\newcommand{\xbD}{\Delta}

\newcommand{\xbG}{\Gamma}

\newcommand{\xbL}{\Lambda}

\newcommand{\xbP}{\Pi}

\newcommand{\xbS}{\Sigma}

\newcommand{\xba}{\alpha}
\newcommand{\xbb}{\beta}

\newcommand{\xbd}{\delta}
\newcommand{\xbe}{\in}
\newcommand{\xbf}{\phi}
\newcommand{\xbg}{\gamma}

\newcommand{\xbk}{\kappa}
\newcommand{\xbl}{\lambda}
\newcommand{\xbm}{\mu}
\newcommand{\xbn}{\nu}
\newcommand{\xbo}{\omega}

\newcommand{\xbq}{\psi}

\newcommand{\xbs}{\sigma}
\newcommand{\xbt}{\tau}

\newcommand{\xCK}{\times}
\newcommand{\xCL}{\pm}
\newcommand{\xCN}{\neg}
\newcommand{\xCQ}{\emptyset}

\newcommand{\xcA}{\forall}

\newcommand{\xcC}{\not\subseteq}

\newcommand{\xcE}{\exists}

\newcommand{\xcM}{\not\models}
\newcommand{\xcN}{\hspace{0.2em}\not\sim\hspace{-0.9em}\mid\hspace{0.8em}}

\newcommand{\xcP}{\not\rightarrow}
\newcommand{\xcS}{\bigcap}
\newcommand{\xcT}{\bot}
\newcommand{\xcU}{\bigwedge}
\newcommand{\xcV}{\bigcup}

\newcommand{\xca}{\infty}
\newcommand{\xcb}{\subset}
\newcommand{\xcc}{\subseteq}
\newcommand{\xcd}{\supseteq}
\newcommand{\xce}{\not\in}

\newcommand{\xcg}{\geq}
\newcommand{\xch}{\Rightarrow}

\newcommand{\xcj}{\Leftrightarrow}
\newcommand{\xck}{\leq}
\newcommand{\xcl}{\vdash}
\newcommand{\xcm}{\models}
\newcommand{\xcn}{\hspace{0.2em}\sim\hspace{-0.9em}\mid\hspace{0.58em}}

\newcommand{\xco}{\vee}
\newcommand{\xcp}{\rightarrow}
\newcommand{\xcq}{\leftarrow}
\newcommand{\xcr}{\leftrightarrow}
\newcommand{\xcs}{\cap}
\newcommand{\xcu}{\wedge}
\newcommand{\xcv}{\cup}

\newcommand{\xcz}{\Box}

\newcommand{\xdC}{\mbox{\boldmath$C$}}
\newcommand{\xdD}{\mbox{\boldmath$D$}}

\newcommand{\xdl}{{\cal L}}
\newcommand{\xdm}{{\cal M}}

\newcommand{\xdp}{{\cal P}}

\newcommand{\xdu}{{\cal U}}

\newcommand{\xdw}{{\cal W}}
\newcommand{\xdx}{{\cal X}}
\newcommand{\xdy}{{\cal Y}}
\newcommand{\xdz}{{\cal Z}}

\newcommand{\xEH}{ & }
\newcommand{\xEI}{\begin{itemize}}
\newcommand{\xEJ}{\end{itemize}}
\newcommand{\xEP}{ \\ }

\newcommand{\xEd}{\neq}
\newcommand{\xEh}{\begin{enumerate}}
\newcommand{\xEj}{\end{enumerate}}

\newcommand{\xeB}{\not\prec}
\newcommand{\xeC}{\not\preceq}

\newcommand{\xeb}{\prec}
\newcommand{\xec}{\preceq}

\newcommand{\xee}{\succ}

\newcommand{\xej}{\lhd}

\newcommand{\xem}{\rhd}

\newcommand{\xex}{\lceil}

\newcommand{\xFO}{\parallel}

\newcommand{\xfA}{\mid}

\newcommand{\Xl}{\ldots}

\newcommand{\ol}{\overline}
\newcommand{\ul}{\underline}

\newcommand{\wt}{\overbrace}

\newcommand{\bl}{\begin{lemma} \rm}
\newcommand{\el}{\end{lemma}}
\newcommand{\br}{\begin{remark} \rm}
\newcommand{\er}{\end{remark}}
\newcommand{\be}{\begin{example} \rm}
\newcommand{\ee}{\end{example}}
\newcommand{\bco}{\begin{corollary} \rm}
\newcommand{\eco}{\end{corollary}}
\newcommand{\bc}{\begin{claim} \rm}
\newcommand{\ec}{\end{claim}}
\newcommand{\bfa}{\begin{fact} \rm}
\newcommand{\efa}{\end{fact}}
\newcommand{\bp}{\begin{proposition} \rm}
\newcommand{\ep}{\end{proposition}}
\newcommand{\bd}{\begin{definition} \rm}
\newcommand{\ed}{\end{definition}}
\newcommand{\bcs}{\begin{construction} \rm}
\newcommand{\ecs}{\end{construction}}
\newcommand{\bcd}{\begin{condition} \rm}
\newcommand{\ecd}{\end{condition}}
\newcommand{\bt}{\begin{theorem} \rm}
\newcommand{\et}{\end{theorem}}
\newcommand{\bn}{\begin{notation} \rm}
\newcommand{\en}{\end{notation}}
\newcommand{\bfi}{\begin{bild} \rm}
\newcommand{\efi}{\end{bild}}

\newcommand{\bfc}{\begin{figure}[htb] \begin{center}}
\newcommand{\efc}{\end{center} \end{figure}}

\documentstyle[12pt]{article}
\title{
DOMAIN CLOSURE CONDITIONS AND DEFINABILITY PRESERVATION
}

\author{Karl Schlechta \\
Laboratoire d'Informatique Fondamentale de Marseille \\
UMR 6166 \\
and \\
Universit\'{e} de Provence \\
CMI, Technop\^{o}le de Ch\^{a}teau-Gombert \\
39, rue Joliot-Curie \\
F-13453 Marseille Cedex 13, France \\
ks@cmi.univ-mrs.fr, schlechta@free.fr, ks1ab@web.de \\
http://www.cmi.univ-mrs.fr/ $\sim$ ks}

\date{June 7, 2006}

\begin{document}

\newtheorem{lemma}{Lemma}[section]
\newtheorem{theorem}[lemma]{Theorem}
\newtheorem{proposition}[lemma]{Proposition}
\newtheorem{corollary}[lemma]{Corollary}
\newtheorem{claim}[lemma]{Claim}
\newtheorem{fact}[lemma]{Fact}
\newtheorem{remark}[lemma]{Remark}
\newtheorem{definition}{Definition}[section]
\newtheorem{construction}{Construction}[section]
\newtheorem{condition}{Condition}[section]
\newtheorem{example}{Example}[section]
\newtheorem{notation}{Notation}[section]
\newtheorem{bild}{Figure}[section]

\maketitle

\pagebreak

\tableofcontents

\pagebreak

\renewcommand{\labelenumi}
  {(\arabic{enumi})}
\renewcommand{\labelenumii}
  {(\arabic{enumi}.\arabic{enumii})}
\renewcommand{\labelenumiii}
  {(\arabic{enumi}.\arabic{enumii}.\arabic{enumiii})}
\renewcommand{\labelenumiv}
  {(\arabic{enumi}.\arabic{enumii}.\arabic{enumiii}.\arabic{enumiv})}

\section*{
ABSTRACT
}

% (+++*** Orig.:   (0)  ABSTRACT  (-> DCB 1) )

We show in Section (2) the importance of closure of the domain under
finite
unions, in particular for Cumulativity, and representation results. We see
that in the absence of this closure, Cumulativity fans out to an infinity
of
different conditions.

We introduce in Section (3) the concept of an algebraic limit, and discuss
its
importance. We then present a representation result for a new concept of
revision, introduced by Booth et al., using approximation by formulas.

We analyse in Section (4) definability preservation problems, and show
that
intersection is the crucial step. We simplify older proofs for the
non-definability cases, and add a new result for ranked structures.

The (very brief) Section 5 puts our considerations in a larger
perspective,
and gives an outlook on open questions.

\section{
INTRODUCTION
}

% (+++*** Orig.:   (1)  INTRODUCTION )

\subsection{
Overview
}

% (+++*** Orig.:   (1.1)  Overview  (-> DCB 1) )

We use and go beyond [Sch04] here. Thus, we will use some results shown
there,
and put them into a more general perspective. We will also re-prove some
results of [Sch04] by more general means. And there are also a number of new
results in the present text, which complement those of [Sch04], for instance,
solving new, or more general cases.

The subject are closure conditions (and related themes) of the domain. We
consider two types, first closure of the domain under simple set-theoretic
operations, particularly under finite unions, second, whether
the algebraic choice functions (or similar objects in the limit case) goes
from
definable sets to definable sets, or to arbitrary subsets, i.e. whether
the
domain is closed under these functions, in other words whether these
functions
preserve definability.

In particular, we will show that the absence of closure under finite
unions has
important consequences for a property called cumulativity. Without this
closure, there is an infinity of different versions of cumulativity, which
all
collaps to one condition in the presence of closure.

This is motivated by the following:
The sets of theory definable propositional model sets - i.e. of the type
$M(T):=\{m:m \xcm T\}$ - are, in the infinite case, closed under arbitrary
intersections, finite unions, but not set difference. When we consider
logics
defined e.g. by suitable sequent calculi, closure under finite unions need
not
hold any longer. This has far reaching consequences for representation
problems,
as the author first noted in [Leh92a]. We investigate the problem and
solutions now
systematically.

Domain closure problems interfere also in the limit variant of several
logics,
and tend to complicate the already relatively complex picture. As an
example of
the limit variant, we take preferential logics. $ \xbm (T)$ (the set of
minimal
models of $T,$ considered in the ``minimal variant'') might well be
empty, even if
$M(T)$ is not empty, e.g. due to infinite
$<-$descending chains of models. This leads to a degenerate case, which
can be
avoided by considering those formulas as consequences, which hold
``finally'' or
``in the limit'', i.e. from a certain level onward, even if there are no
ideal
(minimal) elements.

This limit approach (to various structures and logics) is particularly
recalcitrant, as algebraic and definability preserving problems may occur
together or separately. To distinguish between the two, we introduce the
concept of an algebraic limit, and use this to re-prove as simple
corollaries
of more general situations previous trivialization results, i.e. cases
where the
more complicated limit case can be reduced to the simpler minimal case.

Thus, conceptually, we distinguish between the logical, the algebraic, and
the
structural situation (e.g., preferential structures create certain model
choice
functions, which create certain logics - we made this distinction already
in our
[Sch90] paper), and distinguish now also between the structural and the
algebraic limit. We argue that a ``good'' limit should not only capture the
idea
of a structural limit, but it should also be an algebraic limit, i.e.
capture
the essential algebraic properties of the minimal choice functions. There
might still be definability problems, but our approach allows to
distinguish
them now clearly. The latter, definability, problems are closure
questions,
common to the limit and the minimal variant, and can thus be treated
together
or in a similar way.

Note that these clear distinctions have some philosophical importance,
too.
The structures need an intuitive or philosophical justification, why do we
describe preference by transitive relations, why do we admit copies, etc.?
The resulting algebraic choice functions are void of such questions.

We summarize the distinction:

For e.g. preferential structures, we have:

- logical laws or descriptions like $ \xba \xcn \xba $ - they are the
(imperfect - by
definability preservation problems) reflection of the abstract
description,

- abstract or algebraic semantics, like $ \xbm (A) \xcc A$ - they are the
abstract
description of the foundation,

- structural semantics - they are the intuitive foundation.

Likewise, for the limit situation, we have:

- structural limits - they are again the foundation,

- resulting abstract behaviour, which, again, has to be an abstract or
algebraic
limit, resulting from the structural limit,

- a logical limit, which reflects the abstract limit, and may be plagued
by
definability preservation problems etc. when going from the model to the
logics side.

\subsection{
Basic definitions
}

% (+++*** Orig.:   (1.2)  Basic definitions  (most or all already in CS) )

We summarize now frequently used logical and algebraic properties in the
following table. The left hand column presents the single formula version,
the
center column
the theory version (a theory is, for us, an arbitrary set of formulas),
the
right hand column the algebraic version, describing the choice function on
the
model set, e.g. $f(X) \xcc X$ corresponds to the rule $ \xbf \xcn \xbq $
implies $ \xbf \xcn \xbq $ in the
formula version, and to $ \ol{T} \xcc \ol{ \ol{T} }$ in the theory
version. A short discussion of some
of the properties follows the table.

\bd

$\hspace{0.5em}$

% (*** Orig. No.:  Definition 2.1: )

{\small

\begin{tabular*}{11.5cm}{c@{\extracolsep\fill}|c|c}

(AND) \xEH (AND) \xEH \xEP

$ \xbf \xcn \xbq,  \xbf \xcn \xbq '   \xch $ \xEH
$ T \xcn \xbq, T \xcn \xbq '   \xch $ \xEH
\xEP

$ \xbf \xcn \xbq \xcu \xbq ' $ \xEH
$ T \xcn \xbq \xcu \xbq ' $ \xEH
\xEP

\hline

(OR) \xEH (OR) \xEH $( \xbm \xcv w)$ - $w$ for weak \xEP

$ \xbf \xcn \xbq,  \xbf ' \xcn \xbq   \xch $ \xEH
$T \xcn \xbq, T' \xcn \xbq   \xch $ \xEH
$f(A \xcv B) \xcc f(A) \xcv f(B)$
\xEP

$ \xbf \xco \xbf ' \xcn \xbq $ \xEH
$T \xco T' \xcn \xbq $ \xEH
\xEP

\hline

(LLE) or \xEH (LLE) \xEH \xEP

Left Logical Equivalence \xEH \xEH \xEP

$ \xcl \xbf \xcr \xbf ',  \xbf \xcn \xbq   \xch $ \xEH
$ \ol{T}= \ol{T' }  \xch   \ol{ \ol{T} }= \ol{ \ol{T' } }$ \xEH \xEP

$ \xbf ' \xcn \xbq $ \xEH \xEH \xEP

\hline

(RW) or Right Weakening \xEH (RW) \xEH \xEP

$ \xbf \xcn \xbq,  \xcl \xbq \xcp \xbq '   \xch $ \xEH
$ T \xcn \xbq,  \xcl \xbq \xcp \xbq '   \xch $ \xEH
\xEP

$ \xbf \xcn \xbq ' $ \xEH
$T \xcn \xbq ' $ \xEH
\xEP

\hline

(CCL) or Classical Closure \xEH (CCL) \xEH \xEP

\xEH $ \ol{ \ol{T} }$ is classically \xEH \xEP

\xEH closed \xEH \xEP

\hline

(SC) or Supraclassicality \xEH (SC) \xEH $( \xbm \xcc )$ \xEP

$ \xbf \xcl \xbq $ $ \xch $ $ \xbf \xcn \xbq $ \xEH $ \ol{T} \xcc \ol{
\ol{T} }$ \xEH $f(X) \xcc X$ \xEP
\hline

(CP) or \xEH (CP) \xEH $( \xbm \xCQ )$ \xEP

Consistency Preservation \xEH \xEH \xEP

$ \xbf \xcn \xcT $ $ \xch $ $ \xbf \xcl \xcT $ \xEH $T \xcn \xcT $ $ \xch
$ $T \xcl \xcT $ \xEH $f(X)= \xCQ $ $ \xch $ $X= \xCQ $ \xEP

\hline

(RM) or Rational Monotony \xEH (RM) \xEH $( \xbm =)$ \xEP

$ \xbf \xcn \xbq,  \xbf \xcN \xbq '   \xch $ \xEH
$T \xcn \xbq, T \xcN \xbq '   \xch $ \xEH
$X \xcc Y, Y \xcs f(X) \xEd \xCQ   \xch $
\xEP

$ \xbf \xcu \xbq ' \xcn \xbq $ \xEH
$T \xcv \{ \xbq ' \} \xcn \xbq $ \xEH
$f(X)=f(Y) \xcs X$ \xEP

\hline

(CM) or Cautious Monotony \xEH (CM) \xEH \xEP

$ \xbf \xcn \xbq,  \xbf \xcn \xbq '   \xch $ \xEH
$T \xcc \ol{T' } \xcc \ol{ \ol{T} }  \xch $ \xEH
$f(X) \xcc Y \xcc X  \xch $
\xEP

$ \xbf \xcu \xbq \xcn \xbq ' $ \xEH
$ \ol{ \ol{T} } \xcc \ol{ \ol{T' } }$ \xEH
$f(Y) \xcc f(X)$
\xEP

\hline

(CUM) or Cumulativity \xEH (CUM) \xEH $( \xbm CUM)$ \xEP

$ \xbf \xcn \xbq   \xch $ \xEH
$T \xcc \ol{T' } \xcc \ol{ \ol{T} }  \xch $ \xEH
$f(X) \xcc Y \xcc X  \xch $
\xEP

$( \xbf \xcn \xbq '   \xcj   \xbf \xcu \xbq \xcn \xbq ' )$ \xEH
$ \ol{ \ol{T} }= \ol{ \ol{T' } }$ \xEH
$f(Y)=f(X)$ \xEP

\hline

\xEH (PR) \xEH $( \xbm PR)$ \xEP

$ \ol{ \ol{ \xbf \xcu \xbf ' } }$ $ \xcc $ $ \ol{ \ol{ \ol{ \xbf } } \xcv
\{ \xbf ' \}}$ \xEH
$ \ol{ \ol{T \xcv T' } }$ $ \xcc $ $ \ol{ \ol{ \ol{T} } \xcv T' }$ \xEH
$X \xcc Y$ $ \xch $
\xEP

\xEH \xEH $f(Y) \xcs X \xcc f(X)$
\xEP

\end{tabular*}

}

\ed

(PR) is also called infinite conditionalization - we choose the name for
its
central role for preferential structures (PR) or $( \xbm PR).$ Note that
in the
presence of $( \xbm \xcc ),$ and if $ \xdy $ is closed under finite
intersections, $( \xbm PR)$
is equivalent to

$( \xbm PR' )$ $f(X) \xcs Y \xcc f(X \xcs Y).$

The system of rules (AND), (OR), (LLE), (RW), (SC), (CP), (CM), (CUM) is
also
called system $P$ (for preferential), adding (RM) gives the system $R$
(for
rationality or rankedness).

(AND) is obviously closely related to filters, as we saw already in
Section 1.
(LLE), (RW), (CCL) will all hold automatically, whenever we work with
fixed
model sets. (SC) corresponds to the choice of a subset. (CP) is somewhat
delicate, as it presupposes that the chosen model set is non-empty. This
might
fail in the presence of ever better choices, without ideal ones; the
problem
is addressed by the limit versions - see below in Section 3.4.
(PR) is an inifinitary version of one half of the deduction theorem: Let
$T$ stand
for $ \xbf,$ $T' $ for $ \xbq,$ and $ \xbf \xcu \xbq \xcn \xbs,$ so $
\xbf \xcn \xbq \xcp \xbs,$ but $( \xbq \xcp \xbs ) \xcu \xbq \xcl \xbs.$
(CUM) (whose most interesting half in our context is (CM)) may best be
seen as
normal use of lemmas: We have worked hard and found some lemmas. Now
we can take a rest, and come back again with our new lemmas. Adding them
to the
axioms will neither add new theorems, nor prevent old ones to hold.
(RM) is perhaps best understood by looking at big and small subsets. If
the
set of $ \xbf \xcu \xbq -$models is a big subset of the set of $ \xbf
-$models, and the
set of $ \xbf \xcu \xCN \xbq ' -$models is a not a small subset of the set
of $ \xbf -$models (i.e.
big or of medium size), then the
set of $ \xbf \xcu \xbq \xcu \xbq ' -$models is a big subset of the set of
$ \xbf \xcu \xbq ' -$models.

The following two definitions make preferential structures precise. We
first
give the algebraic definition, and then the definition of the consequence
relation generated by an preferential structure. In the algebraic
definition,
the set $U$ is an arbitrary set, in the application to logic, this will be
the
set of classical models of the underlying propositional language.

In both cases, we first present the simpler variant without copies, and
then
the one with copies. (Note that e.g. [KLM90], [LM92] use labelling
functions
instead, the
version without copies corresponds to injective labelling functions, the
one
with copies to the general case. These are just different ways of
speaking.) We
will discuss the difference between the
version without and the version with copies below, where we show that the
version with copies is strictly more expressive than the version without
copies,
and that transitivity of the relation adds new properties in the case
without
copies. When we summarize our own results below (see Section 2.2.3), we
will
mention that, in the general case with copies, transitivity can be added
without changing properties.

We give here the ``minimal version'', the much more complicated ``limit
version''
is presented and discussed in Section 3. Recall the intuition that the
relation
$ \xeb $ expresses ``normality'' or ``importance'' - the $ \xeb -$smaller, the
more normal or
important. The smallest elements are those which count.

\bd

$\hspace{0.01em}$

% (+++*** Orig. No.:  Definition 1.2: )

Fix $U \xEd \xCQ,$ and consider arbitrary $X.$
Note that this $X$ has not necessarily anything to do with $U,$ or $ \xdu
$ below.
Thus, the functions $ \xbm_{ \xdm }$ below are in principle functions from
$V$ to $V$ - where $V$
is the set theoretical universe we work in.

(A) Preferential models or structures.

(1) The version without copies:

A pair $ \xdm:=<U, \xeb >$ with $U$ an arbitrary set, and $ \xeb $ an
arbitrary binary relation
is called a preferential model or structure.

(2) The version with copies:

A pair $ \xdm:=< \xdu, \xeb >$ with $ \xdu $ an arbitrary set of pairs,
and $ \xeb $ an arbitrary binary
relation is called a preferential model or structure.

If $<x,i> \xbe \xdu,$ then $x$ is intended to be an element of $U,$ and
$i$ the index of the
copy.

(B) Minimal elements, the functions $ \xbm_{ \xdm }$

(1) The version without copies:

Let $ \xdm:=<U, \xeb >,$ and define

$ \xbm_{ \xdm }(X)$ $:=$ $\{x \xbe X:$ $x \xbe U$ $ \xcu $ $ \xCN \xcE x'
\xbe X \xcs U.x' \xeb x\}.$

$ \xbm_{ \xdm }(X)$ is called the set of minimal elements of $X$ (in $
\xdm ).$

(2) The version with copies:

Let $ \xdm:=< \xdu, \xeb >$ be as above. Define

$ \xbm_{ \xdm }(X)$ $:=$ $\{x \xbe X:$ $ \xcE <x,i> \xbe \xdu. \xCN \xcE
<x',i' > \xbe \xdu (x' \xbe X$ $ \xcu $ $<x',i' >' \xeb <x,i>)\}.$

Again, by abuse of language, we say that $ \xbm_{ \xdm }(X)$ is the set of
minimal elements
of $X$ in the structure. If the context is clear, we will also write just
$ \xbm.$

We sometimes say that $<x,i>$ ``kills'' or ``minimizes'' $<y,j>$ if
$<x,i> \xeb <y,j>.$ By abuse of language we also say a set $X$ kills or
minimizes a set
$Y$ if for all $<y,j> \xbe \xdu,$ $y \xbe Y$ there is $<x,i> \xbe \xdu,$
$x \xbe X$ s.t. $<x,i> \xeb <y,j>.$

$ \xdm $ is also called injective or 1-copy, iff there is always at most
one copy
$<x,i>$ for each $x.$

We say that $ \xdm $ is transitive, irreflexive, etc., iff $ \xeb $ is.

\ed

Recall that $ \xbm (X)$ might well be empty, even if $X$ is not.

\bd

$\hspace{0.01em}$

% (+++*** Orig. No.:  Definition 1.3: )

We define the consequence relation of a preferential structure for a
given propositional language $ \xdl.$

(A)

(1) If $m$ is a classical model of a language $ \xdl,$ we say by abuse
of language

$<m,i> \xcm \xbf $ iff $m \xcm \xbf,$

and if $X$ is a set of such pairs, that

$X \xcm \xbf $ iff for all $<m,i> \xbe X$ $m \xcm \xbf.$

(2) If $ \xdm $ is a preferential structure, and $X$ is a set of $ \xdl
-$models for a
classical propositional language $ \xdl,$ or a set of pairs $<m,i>,$
where the $m$ are
such models, we call $ \xdm $ a classical preferential structure or model.

(B)

Validity in a preferential structure, or the semantical consequence
relation
defined by such a structure:

Let $ \xdm $ be as above.

We define:

$T \xcm_{ \xdm } \xbf $ iff $ \xbm_{ \xdm }(M(T)) \xcm \xbf,$ i.e. $
\xbm_{ \xdm }(M(T)) \xcc M( \xbf ).$

$ \xdm $ will be called definability preserving iff for all $X \xbe \xdD_{
\xdl }$ $ \xbm_{ \xdm }(X) \xbe \xdD_{ \xdl }.$

As $ \xbm_{ \xdm }$ is defined on $ \xdD_{ \xdl },$ but need by no means
always result in some new
definable set, this is (and reveals itself as a quite strong) additional
property.

\ed

We define now two additional properties of the relation, smoothness and
rankedness.

The first condition says that if $x \xbe X$ is not a minimal element
of $X,$ then there is $x' \xbe \xbm (X)$ $x' \xeb x.$
In the finite case without copies, smoothness is a trivial consequence of
transitivity and lack of cycles. But note that in the other cases infinite
descending chains might still exist, even if the smoothness condition
holds,
they are just ``short-circuited'': we might have such chains, but below
every
element in the chain is a minimal element. In the author's opinion,
smoothness
is difficult to justify as a structural property (or, in a more
philosophical
spirit, as a property of the world): why should we always have such
minimal
elements below non-minimal ones? Smoothness has, however, a justification
from
its consequences. Its attractiveness comes from two sides:

First, it generates a very valuable logical property, cumulativity (CUM):
If $ \xdm $ is smooth, and
$ \ol{ \ol{T} }$ is the set of $ \xcm_{ \xdm }$-consequences,
then for $T \xcc \ol{T' } \xcc \ol{ \ol{T} }$ $ \xch $ $ \ol{ \ol{T}
}= \ol{ \ol{T' } }.$

Second, for certain approaches, it facilitates completeness proofs, as we
can
look directly at ``ideal'' elements, without having to bother about
intermediate
stages. See in particular the work by Lehmann and his co-authors, [KLM90],
[LM92].

``Smoothness'', or, as it is also called,
``stopperedness'' seems - in the author's opinion - a misnamer. $I$ think it
should
better be called something like ``weak transitivity'': consider the case
where
$a \xee b \xee c,$ but $c \xeB a,$ with $c \xbe \xbm (X).$ It is then not
necessarily the case that
$a \xee c,$ but there is $c' $ ``sufficiently close to $c$'', i.e. in $
\xbm (X),$ s.t. $a \xee c'.$
Results and proof techniques underline this idea. First, in the general
case
with copies, and in the smooth case (in the presence of $( \xcv )!),$
transitivity
does not add new properties, it is ``already present'', second, the
construction
of smoothness by sequences $ \xbs $ (see below in Section 2.3) is very
close in
spirit to a transitive construction.

The second condition, rankedness, seems easier to justify already as a
property
of the structure. It says that, essentially, the elements are ordered in
layers:
If a and $b$ are not comparable, then they are in the same layer. So, if
$c$ is
above (below) $a,$ it will also be above (below) $b$ - like pancakes or
geological
strata. Apart from the triangle
inequality (and leaving aside cardinality questions), this is then just a
distance from some imaginary, ideal point. Again, this property has
important
consequences on the resulting model choice functions and consequence
relations,
making proof techniques for the non-ranked and the ranked case very
different.

\bd

$\hspace{0.01em}$

% (+++*** Orig. No.:  Definition 1.4: )

Let $ \xdz \xcc \xdp (U).$ (In applications to logic, $ \xdz $ will be $
\xdD_{ \xdl }.)$

A preferential structure $ \xdm $ is called $ \xdz -$smooth iff in every
$X \xbe \xdz $ every element
$x \xbe X$ is either minimal in $X$ or above an element, which is minimal
in $X.$ More
precisely:

(1) The version without copies:

If $x \xbe X \xbe \xdz,$ then either $x \xbe \xbm (X)$ or there is $x'
\xbe \xbm (X).x' \xeb x.$

(2) The version with copies:

If $x \xbe X \xbe \xdz,$ and $<x,i> \xbe \xdu,$ then either there is no
$<x',i' > \xbe \xdu,$ $x' \xbe X,$
$<x',i' > \xeb <x,i>$ or there is $<x',i' > \xbe \xdu,$ $<x',i' > \xeb
<x,i>,$ $x' \xbe X,$ s.t. there is
no $<x'',i'' > \xbe \xdu,$ $x'' \xbe X,$ with $<x'',i'' > \xeb <x',i'
>.$

When considering the models of a language $ \xdl,$ $ \xdm $ will be
called smooth iff
it is $ \xdD_{ \xdl }-$smooth; $ \xdD_{ \xdl }$ is the default.

Obviously, the richer the set $ \xdz $ is, the stronger the condition $
\xdz -$smoothness
will be.

\ed

\bd

$\hspace{0.01em}$

% (+++*** Orig. No.:  Definition 1.5: )

A relation $ \xeb_{U}$ on $U$ is called ranked iff
there is an order-preserving function from $U$ to a total order $O,$ $f:U
\xcp O,$ with
$u \xeb_{U}u' $ iff $f(u) \xeb_{O}f(u' ),$ equivalently, if $x$ and $x' $
are $ \xeb_{U}-$incomparable,

then $(y \xeb_{U}x$ iff $y \xeb_{U}x' )$ and $(y \xee_{U}x$ iff $y
\xee_{U}x' )$ for all $y.$

\ed

We conclude with the following standard example:

\be

$\hspace{0.01em}$

% (+++*** Orig. No.:  Example 1.1: )

If $v( \xdl )$ is infinite, and $m$ any model for $ \xdl,$ then $M:=M_{
\xdl }-\{m\}$ is not definable
by any theory $T.$ (Proof: Suppose it were, and let $ \xbf $ hold in $M'
,$
but not in $m,$ so in $m$ $ \xCN \xbf $ holds, but as $ \xbf $ is finite,
there is a model $m' $ in
$M' $ which coincides on all propositional variables of $ \xbf $ with $m,$
so in $m' $ $ \xCN \xbf $
holds, too, a contradiction.)

\ee

\section{
UNIONS AND CUMULATIVITY
}

% (+++*** Orig.:   (2)  UNIONS AND CUMULATIVITY )

\subsection{
Introduction
}

% (+++*** Orig.:   (2.1)  Introduction  (-> DCB 3.3.0) )

This section was motivated by Lehmann's Plausibility Logic, and
re-motivated
by the work of Arieli and Avron, see [AA00]. In both cases, the language does
not have a
built-in ``or'' - resulting in absence $( \xcv )$ of the domain. It is thus
an essay on
the enormous strength of closure of the domain under finite unions, and,
more
generally, on the importance of domain closure conditions.

In the resulting completeness proofs again, a strategy of ``divide and
conquer''
is useful. This helps us to unify
(or extend) our past completeness proofs for the smooth case
in the following way: We will identify more clearly than in the past a
more
or less simple algebraic property - (HU), (HUx) etc. - which allows us to
split
the proofs into two parts. The first part shows validity of the property,
and
this demonstration depends on closure properties, the second part shows
how to
construct a representing structure using the algebraic property. This
second
part will be totally independent from closure properties, and is
essentially an
``administrative'' way to use the property for a construction. This split
approach
allows us thus to isolate the demonstration of the used property from the
construction itself, bringing both parts clearer to light, and simplifying
the proofs, by using common parts.

The reader will see that the successively more complicated conditions
(HU),
(HUx), $( \xbm \xbt )$ reflect well the successively more complicated
situations of
representation:

(HU): smooth (and transitive) structures in the presence of $( \xcv ),$

(HUx): smooth structures in the absence of $( \xcv ),$

$( \xbm \xbt ):$ smooth and transitive structures in the absence of $(
\xcv ).$

This comparison becomes clearer when we see that in the final, most
complicated
case, we will have to carry around all the history of minimization,
$<Y_{0}, \Xl,Y_{n}>,$
necessary for transitivity, which could be summarized in the first case
with
to finite unions. Thus, from an abstract point of view, it is a very
natural
development.

\subsection{
Basic definitions and results
}

% (+++*** Orig.:   (2.2)  Basic definitions and results  (-> DCB 3.3.0) )

We show that, without sufficient closure
properties, there is an infinity of versions of cumulativity, which
collaps
to usual cumulativity when the domain is closed under finite unions.
Closure properties thus reveal themselves as a powerful tool to show
independence of properties.

We work in some fixed arbitrary set $Z,$ all sets considered will be
subsets of $Z.$

We recall or define

\bd

$\hspace{0.01em}$

% (+++*** Orig. No.:  Definition 2.1: )

$( \xbm PR)$ $X \xcc Y$ $ \xcp $ $ \xbm (Y) \xcs X \xcc \xbm (X)$

$( \xbm \xcc )$ $ \xbm (X) \xcc X$

$( \xbm Cum)$ $ \xbm (X) \xcc Y \xcc X$ $ \xcp $ $ \xbm (X)= \xbm (Y)$

$( \xcv )$ is closure of the domain under finite unions.

$( \xcs )$ is closure of the domain under finite intersections.

$( \xcS )$ is closure of the domain under arbitrary intersections.

\ed

We use without further mentioning $( \xbm PR)$ and $( \xbm \xcc ).$

\bd

$\hspace{0.01em}$

% (+++*** Orig. No.:  Definition 2.2: )

For any ordinal $ \xba,$ we define

$( \xbm Cum \xba ):$

If for all $ \xbb \xck \xba $ $ \xbm (X_{ \xbb }) \xcc U \xcv \xcV \{X_{
\xbg }: \xbg < \xbb \}$ hold, then so does
$ \xcS \{X_{ \xbg }: \xbg \xck \xba \} \xcs \xbm (U) \xcc \xbm (X_{ \xba
}).$

$( \xbm Cumt \xba ):$

If for all $ \xbb \xck \xba $ $ \xbm (X_{ \xbb }) \xcc U \xcv \xcV \{X_{
\xbg }: \xbg < \xbb \}$ hold, then so does
$X_{ \xba } \xcs \xbm (U) \xcc \xbm (X_{ \xba }).$

(``t'' stands for transitive, see Fact 2.1, (2.2) below.)

$( \xbm Cum \xca )$ and $( \xbm Cumt \xca )$ will be the class of all $(
\xbm Cum \xba )$ or $( \xbm Cumt \xba )$ -
read their ``conjunction'', i.e. if we say that $( \xbm Cum \xca )$ holds,
we mean that
all $( \xbm Cum \xba )$ hold.

The following inductive definition of $H(U,x)$ and of the property (HUx)
concerns closure under $( \xbm Cum \xca ),$ its main
property is formulated in Fact 2.1, (8), its main interest is its use in
the
proof of Proposition 2.2.

$H(U,x)_{0}$ $:=$ $U,$

$H(U,x)_{ \xba +1}$ $:=$ $H(U,x)_{ \xba }$ $ \xcv $ $ \xcV \{X:$ $x \xbe
X$ $ \xcu $ $ \xbm (X) \xcc H(U,x)_{ \xba }\},$

$H(U,x)_{ \xbl }$ $:=$ $ \xcV \{H(U,x)_{ \xba }: \xba < \xbl \}$ for
$limit( \xbl ),$

$H(U,x)$ $:=$ $ \xcV \{H(U,x)_{ \xba }: \xba < \xbk \}$ for $ \xbk $
sufficiently big $(card(Z)$ suffices, as

the procedure trivializes, when we cannot add any new elements).

(HUx) is the property:

$x \xbe \xbm (U),$ $x \xbe Y- \xbm (Y)$ $ \xcp $ $ \xbm (Y) \xcC H(U,x)$ -
of course for all $x$ and $U.$
$(U,Y \xbe \xdy ).$

For the case with $( \xcv ),$ we further define, independent of $x,$

$H(U)_{0}$ $:=$ $U,$

$H(U)_{ \xba +1}$ $:=$ $H(U)_{ \xba }$ $ \xcv $ $ \xcV \{X:$ $ \xbm (X)
\xcc H(U)_{ \xba }\},$

$H(U)_{ \xbl }$ $:=$ $ \xcV \{H(U)_{ \xba }: \xba < \xbl \}$ for $limit(
\xbl ),$

$H(U)$ $:=$ $ \xcV \{H(U)_{ \xba }: \xba < \xbk \}$ again for $ \xbk $
sufficiently big

(HU) is the property:

$x \xbe \xbm (U),$ $x \xbe Y- \xbm (Y)$ $ \xcp $ $ \xbm (Y) \xcC H(U)$ -
of course for all $U.$
$(U,Y \xbe \xdy ).$

\ed

Obviously, $H(U,x) \xcc H(U),$ so $(HU) \xcp (HUx).$

\paragraph{
Note:
}

$\hspace{0.01em}$

The first conditions thus have the form:

$( \xbm Cum0)$ $ \xbm (X_{0}) \xcc U$ $ \xcp $ $X_{0} \xcs \xbm (U) \xcc
\xbm (X_{0}),$

$( \xbm Cum1)$ $ \xbm (X_{0}) \xcc U,$ $ \xbm (X_{1}) \xcc U \xcv X_{0}$ $
\xcp $ $X_{0} \xcs X_{1} \xcs \xbm (U) \xcc \xbm (X_{1}),$

$( \xbm Cum2)$ $ \xbm (X_{0}) \xcc U,$ $ \xbm (X_{1}) \xcc U \xcv X_{0},$
$ \xbm (X_{2}) \xcc U \xcv X_{0} \xcv X_{1}$ $ \xcp $ $X_{0} \xcs X_{1}
\xcs X_{2} \xcs \xbm (U) \xcc \xbm (X_{2}).$

$( \xbm Cumt \xba )$ differs from $( \xbm Cum \xba )$ only in the
consequence, the intersection contains
only the last $X_{ \xba }$ - in particular, $( \xbm Cum0)$ and $( \xbm
Cumt0)$ coincide.

Recall that condition $( \xbm Cum1)$ is the crucial condition in [Leh92a],
which
failed, despite $( \xbm Cum),$ but which has to hold in all smooth models.
This
condition $( \xbm Cum1)$ was the starting point of the investigation.

\be

$\hspace{0.01em}$

% (+++*** Orig. No.:  Example 2.1: )

Perhaps the main result of this section is the following example, which
shows
that above conditions are all different in the absence of $( \xcv ),$ in
its
presence they all collaps (see  \Xl. below). More precisely,
the following (class of) $example(s)$ shows that the $( \xbm Cum \xba )$
increase in
strength. For any finite or infinite ordinal $ \xbk >0$ we construct an
example s.t.

(a) $( \xbm PR)$ and $( \xbm \xcc )$ hold

(b) $( \xbm Cum)$ holds

(c) $( \xcS )$ holds

(d) $( \xbm Cumt \xba )$ holds for $ \xba < \xbk $

(e) $( \xbm Cum \xbk )$ fails.

\ee

We define a suitable base set and a non-transitive binary relation $ \xeb
$
on this set, as well as a suitable set $ \xdx $ of subsets, closed under
arbitrary
intersections, but not under finite unions, and define $ \xbm $ on these
subsets
as usual in preferential structures by $ \xeb.$ Thus, $( \xbm PR)$ and $(
\xbm \xcc )$ will hold.
It will be immediate that $( \xbm Cum \xbk )$ fails, and we will show that
$( \xbm Cum)$ and
$( \xbm Cumt \xba )$ for $ \xba < \xbk $ hold by examining the cases.

For simplicity, we first define a set of generators for $ \xdx,$ and
close under
$( \xcS )$ afterwards. The set $U$ will have a special position, it is the
``useful''
starting point to construct chains corresponding to above definitions of
$( \xbm Cum \xba )$ and $( \xbm Cumt \xba ).$

\bn

$\hspace{0.01em}$

% (+++*** Orig. No.:  Notation: )

i,j will be successor ordinals, $ \xbl $ etc. limit ordinals, $ \xba,$ $
\xbb,$ $ \xbk $ any ordinals,
thus e.g. $ \xbl \xck \xbk $ will imply that $ \xbl $ is a limit ordinal $
\xck \xbk,$ etc.

\paragraph{
The base set and the relation $ \xeb $:
}

$\hspace{0.01em}$

$ \xbk >0$ is fixed, but arbitrary. We go up to $ \xbk >0.$

The base set is $\{a,b,c\}$ $ \xcv $ $\{d_{ \xbl }: \xbl \xck \xbk \}$ $
\xcv $ $\{x_{ \xba }: \xba \xck \xbk +1\}$ $ \xcv $ $\{x'_{ \xba }: \xba
\xck \xbk \}.$
$a \xeb b \xeb c,$ $x_{ \xba } \xeb x_{ \xba +1},$ $x_{ \xba } \xeb x'_{
\xba },$ $x'_{0} \xeb x_{ \xbl }$ (for any $ \xbl )$ - $ \xeb $ is NOT
transitive.

\paragraph{
The generators:
}

$\hspace{0.01em}$

$U:=\{a,c,x_{0}\} \xcv \{d_{ \xbl }: \xbl \xck \xbk \}$ - i.e. $ \Xl
.\{d_{ \xbl }:lim( \xbl ) \xcu \xbl \xck \xbk \},$

$X_{i}:=\{c,x_{i},x'_{i},x_{i+1}\}$ $(i< \xbk ),$

$X_{ \xbl }:=\{c,d_{ \xbl },x_{ \xbl },x'_{ \xbl },x_{ \xbl +1}\} \xcv
\{x'_{ \xba }: \xba < \xbl \}$ $( \xbl < \xbk ),$

$X'_{ \xbk }:=\{a,b,c,x_{ \xbk },x'_{ \xbk },x_{ \xbk +1}\}$ if $ \xbk $
is a successor,

$X'_{ \xbk }:=\{a,b,c,d_{ \xbk },x_{ \xbk },x'_{ \xbk },x_{ \xbk +1}\}
\xcv \{x'_{ \xba }: \xba < \xbk \}$ if $ \xbk $ is a limit.

Thus, $X'_{ \xbk }=X_{ \xbk } \xcv \{a,b\}$ if $X_{ \xbk }$ were defined.

Note that there is only one $X'_{ \xbk },$ and $X_{ \xba }$ is defined
only for $ \xba < \xbk,$ so we will
not have $X_{ \xba }$ and $X'_{ \xba }$ at the same time.

Thus, the values of the generators under $ \xbm $ are:

$ \xbm (U)=U,$

$ \xbm (X_{i})=\{c,x_{i}\},$

$ \xbm (X_{ \xbl })=\{c,d_{ \xbl }\} \xcv \{x'_{ \xba }: \xba < \xbl \},$

$ \xbm (X'_{i})=\{a,x_{i}\}$ $(i>0!),$

$ \xbm (X'_{ \xbl })=\{a,d_{ \xbl }\} \xcv \{x'_{ \xba }: \xba < \xbl \}.$

(We do not assume that the domain is closed under $ \xbm.)$

\paragraph{
Intersections:
}

$\hspace{0.01em}$

We consider first pairwise intersections:

(1) $U \xcs X_{0}=\{c,x_{0}\},$

(2) $U \xcs X_{i}=\{c\},$ $i>0,$

(3) $U \xcs X_{ \xbl }=\{c,d_{ \xbl }\},$

(4) $U \xcs X'_{i}=\{a,c\}$ $(i>0!),$

(5) $U \xcs X'_{ \xbl }=\{a,c,d_{ \xbl }\},$

(6) $X_{i} \xcs X_{j}:$

(6.1) $j=i+1$ $\{c,x_{i+1}\},$

(6.2) else $\{c\},$

(7) $X_{i} \xcs X_{ \xbl }:$

(7.1) $i< \xbl $ $\{c,x'_{i}\},$

(7.2) $i= \xbl +1$ $\{c,x_{ \xbl +1}\},$

(7.3) $i> \xbl +1$ $\{c\},$

(8) $X_{ \xbl } \xcs X_{ \xbl ' }:$ $\{c\} \xcv \{x'_{ \xba }: \xba \xck
min( \xbl, \xbl ' )\}.$

As $X'_{ \xbk }$ occurs only once, $X_{ \xba } \xcs X'_{ \xbk }$ etc. give
no new results.

Note that $ \xbm $ is constant on all these pairwise intersections.

Iterated intersections:

As $c$ is an element of all sets, sets of the type $\{c,z\}$ do not give
any
new results. The possible subsets of $\{a,c,d_{ \xbl }\}:$ $\{c\},$
$\{a,c\},$ $\{c,d_{ \xbl }\}$ exist
already. Thus, the only source of new sets via iterated intersections is
$X_{ \xbl } \xcs X_{ \xbl ' }=\{c\} \xcv \{x'_{ \xba }: \xba \xck min(
\xbl, \xbl ' )\}.$ But, to intersect them, or with some
old sets, will not generate any new sets either. Consequently, the example
satisfies $( \xcS )$ for $ \xdx $ defined by $U,$ $X_{i}$ $(i< \xbk ),$
$X_{ \xbl }$ $( \xbl < \xbk ),$ $X'_{ \xbk },$ and above
paiwise intersections.

We will now verify the positive properties. This is tedious, but
straightforward, checking the different cases.

\paragraph{
Validity of $( \xbm Cum)$:
}

$\hspace{0.01em}$

Consider the prerequisite $ \xbm (X) \xcc Y \xcc X.$ If $ \xbm (X)=X$ or
if $X- \xbm (X)$ is a singleton,
$X$ cannot give a violation of $( \xbm Cum).$ So we are left with the
following
candidates for $X:$

(1) $X_{i}:=\{c,x_{i},x'_{i},x_{i+1}\},$ $ \xbm (X_{i})=\{c,x_{i}\}$

Interesting candidates for $Y$ will have 3 elements, but they will all
contain
$a.$ (If $ \xbk < \xbo:$ $U=\{a,c,x_{0}\}.)$

(2) $X_{ \xbl }:=\{c,d_{ \xbl },x_{ \xbl },x'_{ \xbl },x_{ \xbl +1}\} \xcv
\{x'_{ \xba }: \xba < \xbl \},$ $ \xbm (X_{ \xbl })=\{c,d_{ \xbl }\} \xcv
\{x'_{ \xba }: \xba < \xbl \}$

The only sets to contain $d_{ \xbl }$ are $X_{ \xbl },$ $U,$ $U \xcs X_{
\xbl }.$ But $a \xbe U,$ and
$U \xcs X_{ \xbl }$ ist finite. $(X_{ \xbl }$ and $X'_{ \xbl }$ cannot be
present at the same time.)

(3) $X'_{i}:=\{a,b,c,x_{i},x'_{i},x_{i+1}\},$ $ \xbm (X'_{i})=\{a,x_{i}\}$

a is only in $U,$ $X'_{i},$ $U \xcs X'_{i}=\{a,c\},$ but $x_{i} \xce U,$
as $i>0.$

(4) $X'_{ \xbl }:=\{a,b,c,d_{ \xbl },x_{ \xbl },x'_{ \xbl },x_{ \xbl +1}\}
\xcv \{x'_{ \xba }: \xba < \xbl \},$ $ \xbm (X'_{ \xbl })=\{a,d_{ \xbl }\}
\xcv \{x'_{ \xba }: \xba < \xbl \}$

$d_{ \xbl }$ is only in $X'_{ \xbl }$ and $U,$ but $U$ contains no $x'_{
\xba }.$

Thus, $( \xbm Cum)$ holds trivially.

\paragraph{
$( \xbm Cumt \xba )$ hold for $ \xba < \xbk $:
}

$\hspace{0.01em}$

To simplify language, we say that we reach $Y$ from $X$ iff $X \xEd Y$ and
there is a
sequence $X_{ \xbb },$ $ \xbb \xck \xba $ and $ \xbm (X_{ \xbb }) \xcc X
\xcv \xcV \{X_{ \xbg }: \xbg < \xbb \},$ and $X_{ \xba }=Y,$ $X_{0}=X.$
Failure of $( \xbm Cumt \xba )$ would then mean that there are $X$ and
$Y,$ we can reach
$Y$ from $X,$ and $x \xbe ( \xbm (X) \xcs Y)- \xbm (Y).$ Thus, in a
counterexample, $Y= \xbm (Y)$ is
impossible, so none of the intersections can be such $Y.$

To reach $Y$ from $X,$ we have to get started from $X,$ i.e. there must be
$Z$ s.t.
$ \xbm (Z) \xcc X,$ $Z \xcC X$ (so $ \xbm (Z) \xEd Z).$ Inspection of the
different cases shows that
we cannot reach any set $Y$ from any case of the intersections, except
from
(1), (6.1), (7.2).

If $Y$ contains a globally minimal element (i.e. there is no smaller
element in
any set), it can only be reached from any $X$ which already contains this
element. The globally minimal elements are $a,$ $x_{0},$ and the $d_{ \xbl
},$ $ \xbl \xck \xbk.$

By these observations, we see that $X_{ \xbl }$ and $X'_{ \xbk }$ can only
be reached from $U.$
From no $X_{ \xba }$ $U$ can be reached, as the globally minimal a is
missing. But
$U$ cannot be reached from $X'_{ \xbk }$ either, as the globally minimal
$x_{0}$ is missing.

When we look at the relation $ \xeb $ defining $ \xbm,$ we see that we
can reach $Y$ from $X$
only by going upwards, adding bigger elements. Thus, from $X_{ \xba },$ we
cannot reach
any $X_{ \xbb },$ $ \xbb < \xba,$ the same holds for $X'_{ \xbk }$ and
$X_{ \xbb },$ $ \xbb < \xbk.$ Thus, from $X'_{ \xbk },$ we
cannot go anywhere interesting (recall that the intersections are not
candidates
for a $Y$ giving a contradiction).

Consider now $X_{ \xba }.$ We can go up to any $X_{ \xba +n},$ but not to
any $X_{ \xbl },$ $ \xba < \xbl,$ as
$d_{ \xbl }$ is missing, neither to $X'_{ \xbk },$ as a is missing. And we
will be stopped by
the first $ \xbl > \xba,$ as $x_{ \xbl }$ will be missing to go beyond
$X_{ \xbl }.$ Analogous observations
hold for the remaining intersections (1), (6.1), (7.2). But in all these
sets we
can reach, we will not destroy minimality of any element of $X_{ \xba }$
(or of the
intersections).

Consequently, the only candidates for failure will all start with $U.$ As
the only
element of $U$ not globally minimal is $c,$ such failure has to have $c
\xbe Y- \xbm (Y),$ so
$Y$ has to be $X'_{ \xbk }.$ Suppose we omit one of the $X_{ \xba }$ in
the sequence going up to
$X'_{ \xbk }.$ If $ \xbk \xcg \xbl > \xba,$ we cannot reach $X_{ \xbl }$
and beyond, as $x'_{ \xba }$ will be missing.
But we cannot go to $X_{ \xba +n}$ either, as $x_{ \xba +1}$ is missing.
So we will be stopped
at $X_{ \xba }.$ Thus, to see failure, we need the full sequence
$U=X_{0},$ $X'_{ \xbk }=Y_{ \xbk },$
$Y_{ \xba }=X_{ \xba }$ for $0< \xba < \xbk.$

\paragraph{
$( \xbm Cum \xbk )$ fails:
}

$\hspace{0.01em}$

The full sequence $U=X_{0},$ $X'_{ \xbk }=Y_{ \xbk },$ $Y_{ \xba }=X_{
\xba }$ for $0< \xba < \xbk $ shows this, as
$c \xbe \xbm (U) \xcs X'_{ \xbk },$ but $c \xce \xbm (X'_{ \xbk }).$

Consequently, the example satisfies $( \xcS ),$ $( \xbm Cum),$ $( \xbm
Cumt \xba )$ for $ \xba < \xbk,$ and
$( \xbm Cum \xbk )$ fails.

$ \xcz $
\\[3ex]

\en

\bfa

$\hspace{0.01em}$

% (+++*** Orig. No.:  Fact 2.1: )

We summarize some properties of $( \xbm Cum \xba )$ and $( \xbm Cumt \xba
)$ - sometimes with some
redundancy. Unless said otherwise, $ \xba,$ $ \xbb $ etc. will be
arbitrary ordinals.

For (1) to (6) $( \xbm PR)$ and $( \xbm \xcc )$ are assumed to hold, for
(7) and (8) only
$( \xbm \xcc ).$

(1) Downward:

(1.1) $( \xbm Cum \xba )$ $ \xcp $ $( \xbm Cum \xbb )$ for all $ \xbb \xck
\xba $

(1.2) $( \xbm Cumt \xba )$ $ \xcp $ $( \xbm Cumt \xbb )$ for all $ \xbb
\xck \xba $

\efa

(2) Validity of $( \xbm Cum \xba )$ and $( \xbm Cumt \xba )$:

(2.1) All $( \xbm Cum \xba )$ hold in smooth preferential structures

(2.2) All $( \xbm Cumt \xba )$ hold in transitive smooth preferential
structures

(2.3) $( \xbm Cumt \xba )$ for $0< \xba $ do not necessarily hold in
smooth structures without
transitivity, even in the presence of $( \xcS )$

(3) Upward:

(3.1) $( \xbm Cum \xbb )$ $+$ $( \xcv )$ $ \xcp $ $( \xbm Cum \xba )$ for
all $ \xbb \xck \xba $

(3.2) $( \xbm Cumt \xbb )$ $+$ $( \xcv )$ $ \xcp $ $( \xbm Cumt \xba )$
for all $ \xbb \xck \xba $

(3.3) $\{( \xbm Cumt \xbb ): \xbb < \xba \}$ $+$ $( \xbm Cum)$ $+$ $( \xcS
)$ $ \xcP $ $( \xbm Cum \xba )$ for $ \xba >0.$

(4) Connection $( \xbm Cum \xba )/( \xbm Cumt \xba )$:

(4.1) $( \xbm Cumt \xba )$ $ \xcp $ $( \xbm Cum \xba )$

(4.2) $( \xbm Cum \xba )$ $+$ $( \xcS )$ $ \xcP $ $( \xbm Cumt \xba )$

(4.3) $( \xbm Cum \xba )$ $+$ $( \xcv )$ $ \xcp $ $( \xbm Cumt \xba )$

(5) $( \xbm Cum)$ and $( \xbm Cumi)$:

(5.1) $( \xbm Cum)$ $+$ $( \xcv )$ entail:

(5.1.1) $ \xbm (A) \xcc B$ $ \xcp $ $ \xbm (A \xcv B)= \xbm (B)$

(5.1.2) $ \xbm (X) \xcc U,$ $U \xcc Y$ $ \xcp $ $ \xbm (Y \xcv X)= \xbm
(Y)$

(5.1.3) $ \xbm (X) \xcc U,$ $U \xcc Y$ $ \xcp $ $ \xbm (Y) \xcs X \xcc
\xbm (U)$

(5.2) $( \xbm Cum \xba )$ $ \xcp $ $( \xbm Cum)$ for all $ \xba $

(5.3) $( \xbm Cum)$ $+$ $( \xcv )$ $ \xcp $ $( \xbm Cum \xba )$ for all $
\xba $

(5.4) $( \xbm Cum)$ $+$ $( \xcs )$ $ \xcp $ $( \xbm Cum0)$

(6) $( \xbm Cum)$ and $( \xbm Cumt \xba )$:

(6.1) $( \xbm Cumt \xba )$ $ \xcp $ $( \xbm Cum)$ for all $ \xba $

(6.2) $( \xbm Cum)$ $+$ $( \xcv )$ $ \xcp $ $( \xbm Cumt \xba )$ for all $
\xba $

(6.3) $( \xbm Cum)$ $ \xcP $ $( \xbm Cumt \xba )$ for all $ \xba >0$

(7) $( \xbm Cum0)$ $ \xcp $ $( \xbm PR)$

(8) $( \xbm Cum \xca )$ and (HUx):

(8.1) $x \xbe \xbm (Y),$ $ \xbm (Y) \xcc H(U,x)$ $ \xcp $ $Y \xcc H(U,x)$

(8.2) $( \xbm Cum \xca )$ $ \xcp $ (HUx)

(8.3) (HUx) $ \xcp $ $( \xbm Cum \xca )$

\paragraph{
Proof:
}

$\hspace{0.01em}$

We prove these facts in a different order: (1), (2), (5.1), (5.2), (4.1),
(6.1),
(6.2), (5.3), (3.1), (3.2), (4.2), (4.3), (5.4), (3.3), (6.3), (7), (8).

(1.1)

For $ \xbb < \xbg \xck \xba $ set $X_{ \xbg }:=X_{ \xbb }.$ Let the
prerequisites of $( \xbm Cum \xbb )$ hold. Then for
$ \xbg $ with $ \xbb < \xbg \xck \xba $ $ \xbm (X_{ \xbg }) \xcc X_{ \xbb
}$ by $( \xbm \xcc ),$ so the prerequisites
of $( \xbm Cum \xba )$ hold, too, so by $( \xbm Cum \xba )$ $ \xcS \{X_{
\xbd }: \xbd \xck \xbb \} \xcs \xbm (U)$ $=$
$ \xcS \{X_{ \xbd }: \xbd \xck \xba \} \xcs \xbm (U)$ $ \xcc $ $ \xbm (X_{
\xba })$ $=$ $ \xbm (X_{ \xbb }).$

(1.2)

Analogous.

(2.1)

Proof by induction.

$( \xbm Cum0)$ Let $ \xbm (X_{0}) \xcc U,$ suppose there is $x \xbe \xbm
(U) \xcs (X_{0}- \xbm (X_{0})).$ By smoothness,
there is $y \xeb x,$ $y \xbe \xbm (X_{0}) \xcc U,$ $contradiction$ (The
same arguments works for copies: all
copies of $x$ must be minimized by some $y \xbe \xbm (X_{0}),$ but at
least one copy of $x$
has to be minimal in U.)

Suppose $( \xbm Cum \xbb )$ hold for all $ \xbb < \xba.$ We show $( \xbm
Cum \xba ).$ Let the prerequisites
of $( \xbm Cum \xba )$ hold, then those for $( \xbm Cum \xbb ),$ $ \xbb <
\xba $ hold, too. Suppose there is
$x \xbe \xbm (U) \xcs \xcS \{X_{ \xbg }: \xbg \xck \xba \}- \xbm (X_{ \xba
}).$ So by $( \xbm Cum \xbb )$ for $ \xbb < \xba $ $x \xbe \xbm (X_{ \xbb
})$
moreover $x \xbe \xbm (U).$ By smoothness, there is $y \xbe \xbm (X_{ \xba
}) \xcc U \xcv \xcV \{X_{ \xbb }: \xbb < \xba \},$ $y \xeb x,$
but this is a contradiction. The same argument works again for copies.

(2.2)

We use the following Fact:
Let, in a smooth transitive structure, $ \xbm (X_{ \xbb })$ $ \xcc $ $U
\xcv \xcV \{X_{ \xbg }: \xbg < \xbb \}$ for all $ \xbb \xck \xba,$
and let $x \xbe \xbm (U).$ Then there is no $y \xeb x,$ $y \xbe U \xcv
\xcV \{X_{ \xbg }: \xbg \xck \xba \}.$

Proof of the Fact by induction:
$ \xba =0:$ $y \xbe U$ is impossible: if $y \xbe X_{0},$ then if $y \xbe
\xbm (X_{0}) \xcc U,$ which is impossible,
or there is $z \xbe \xbm (X_{0}),$ $z \xeb y,$ so $z \xeb x$ by
transitivity, but $ \xbm (X_{0}) \xcc U.$
Let the result hold for all $ \xbb < \xba,$ but fail for $ \xba,$
so $ \xCN \xcE y \xeb x.y \xbe U \xcv \xcV \{X_{ \xbg }: \xbg < \xba \},$
but $ \xcE y \xeb x.y \xbe U \xcv \xcV \{X_{ \xbg }: \xbg \xck \xba \},$
so $y \xbe X_{ \xba }.$
If $y \xbe \xbm (X_{ \xba }),$ then $y \xbe U \xcv \xcV \{X_{ \xbg }: \xbg
< \xba \},$ but this is impossible, so
$y \xbe X_{ \xba }- \xbm (X_{ \xba }),$ let by smoothness $z \xeb y,$ $z
\xbe \xbm (X_{ \xba }),$ so by transitivity $z \xeb x,$ $contradiction.$
The result is easily modified for the case with copies.

Let the prerequisites of $( \xbm Cumt \xba )$ hold, then those of the Fact
will hold,
too. Let now $x \xbe \xbm (U) \xcs (X_{ \xba }- \xbm (X_{ \xba })),$ by
smoothness, there must be $y \xeb x,$
$y \xbe \xbm (X_{ \xba }) \xcc U \xcv \xcV \{X_{ \xbg }: \xbg < \xba \},$
contradicting the Fact.

(2.3)

Let $ \xba >0,$ and consider the following structure over $\{a,b,c\}:$
$U:=\{a,c\},$
$X_{0}:=\{b,c\},$ $X_{ \xba }:= \Xl:=X_{1}:=\{a,b\},$ and their
intersections, $\{a\},$ $\{b\},$ $\{c\},$ $ \xCQ $ with
the order $c \xeb b \xeb a$ (without transitivity). This is preferential,
so $( \xbm PR)$ and
$( \xbm \xcc )$ hold.
The structure is smooth for $U,$ all $X_{ \xbb },$ and their
intersections.
We have $ \xbm (X_{0}) \xcc U,$ $ \xbm (X_{ \xbb }) \xcc U \xcv X_{0}$ for
all $ \xbb \xck \xba,$ so $ \xbm (X_{ \xbb }) \xcc U \xcv \xcV \{X_{ \xbg
}: \xbg < \xbb \}$
for all $ \xbb \xck \xba $ but $X_{ \xba } \xcs \xbm (U)=\{a\} \xcC \{b\}=
\xbm (X_{ \xba })$ for $ \xba >0.$

(5.1)

(5.1.1) $ \xbm (A) \xcc B$ $ \xcp $ $ \xbm (A \xcv B) \xcc \xbm (A) \xcv
\xbm (B) \xcc B$ $ \xcp_{( \xbm Cum)}$ $ \xbm (B)= \xbm (A \xcv B).$

(5.1.2) $ \xbm (X) \xcc U \xcc Y$ $ \xcp $ (by (1)) $ \xbm (Y \xcv X)=
\xbm (Y).$

(5.1.3) $ \xbm (Y) \xcs X$ $=$ (by (2)) $ \xbm (Y \xcv X) \xcs X$ $ \xcc $
$ \xbm (Y \xcv X) \xcs (X \xcv U)$ $ \xcc $ (by $( \xbm PR))$
$ \xbm (X \xcv U)$ $=$ (by (1)) $ \xbm (U).$

(5.2)

Using (1.1), it suffices to show $( \xbm Cum0)$ $ \xcp $ $( \xbm Cum).$
Let $ \xbm (X) \xcc U \xcc X.$ By $( \xbm Cum0)$ $X \xcs \xbm (U) \xcc
\xbm (X),$ so by $ \xbm (U) \xcc U \xcc X$ $ \xcp $ $ \xbm (U) \xcc \xbm
(X).$
$U \xcc X$ $ \xcp $ $ \xbm (X) \xcs U \xcc \xbm (U),$ but also $ \xbm (X)
\xcc U,$ so $ \xbm (X) \xcc \xbm (U).$

(4.1)

Trivial.

(6.1)

Follows from (4.1) and (5.2).

(6.2)

Let the prerequisites of $( \xbm Cumt \xba )$ hold.

We first show by induction $ \xbm (X_{ \xba } \xcv U) \xcc \xbm (U).$

Proof:

$ \xba =0:$ $ \xbm (X_{0}) \xcc U$ $ \xcp $ $ \xbm (X_{0} \xcv U)= \xbm
(U)$ by (5.1.1).
Let for all $ \xbb < \xba $ $ \xbm (X_{ \xbb } \xcv U) \xcc \xbm (U) \xcc
U.$ By prerequisite,
$ \xbm (X_{ \xba }) \xcc U \xcv \xcV \{X_{ \xbb }: \xbb < \xba \},$ thus $
\xbm (X_{ \xba } \xcv U)$ $ \xcc $ $ \xbm (X_{ \xba }) \xcv \xbm (U)$ $
\xcc $ $ \xcV \{U \xcv X_{ \xbb }: \xbb < \xba \},$

so $ \xcA \xbb < \xba $ $ \xbm (X_{ \xba } \xcv U) \xcs (U \xcv X_{ \xbb
})$ $ \xcc $ $ \xbm (U)$ by (5.1.3), thus $ \xbm (X_{ \xba } \xcv U) \xcc
\xbm (U).$

Consequently, under the above prerequisites, we have $ \xbm (X_{ \xba }
\xcv U)$ $ \xcc $ $ \xbm (U)$ $ \xcc $
$U$ $ \xcc $ $U \xcv X_{ \xba },$ so by $( \xbm Cum)$ $ \xbm (U)= \xbm
(X_{ \xba } \xcv U),$ and, finally,
$ \xbm (U) \xcs X_{ \xba }= \xbm (X_{ \xba } \xcv U) \xcs X_{ \xba } \xcc
\xbm (X_{ \xba })$ by $( \xbm PR).$

Note that finite unions take us over the limit step, essentially, as all
steps collaps, and $ \xbm (X_{ \xba } \xcv U)$ will always be $ \xbm (U),$
so there are no real
changes.

(5.3)

Follows from (6.2) and (4.1).

(3.1)

Follows from (5.2) and (5.3).

(3.2)

Follows from (6.1) and (6.2).

(4.2)

Follows from (2.3) and (2.1).

(4.3)

Follows from (5.2) and (6.2).

(5.4)

$ \xbm (X) \xcc U$ $ \xcp $ $ \xbm (X) \xcc U \xcs X \xcc X$ $ \xcp $ $
\xbm (X \xcs U)= \xbm (X)$ $ \xcp $
$X \xcs \xbm (U)=(X \xcs U) \xcs \xbm (U) \xcc \xbm (X \xcs U)= \xbm (X)$

(3.3)

See Example 2.1.

(6.3)

This is a consequence of (3.3).

(7)

Trivial. Let $X \xcc Y,$ so by $( \xbm \xcc )$ $ \xbm (X) \xcc X \xcc Y,$
so by $( \xbm Cum0)$ $X \xcs \xbm (Y) \xcc \xbm (X).$

(8.1)

Trivial by definition of $H(U,x).$

(8.2)

Let $x \xbe \xbm (U),$ $x \xbe Y,$ $ \xbm (Y) \xcc H(U,x)$ (and thus $Y
\xcc H(U,x)$ by definition).
Thus, we have a sequence $X_{0}:=U,$ $ \xbm (X_{ \xbb }) \xcc U \xcv \xcV
\{X_{ \xbg }: \xbg < \xbb \}$ and $Y=X_{ \xba }$ for some $ \xba $
(after $X_{0},$ enumerate arbitrarily $H(U,x)_{1},$ then $H(U,x)_{2},$
etc., do nothing at
limits). So $x \xbe \xcS \{X_{ \xbg }: \xbg \xck \xba \} \xcs \xbm (U),$
and $x \xbe \xbm (X_{ \xba })= \xbm (Y)$ by $( \xbm Cum \xca ).$
Remark: The same argument shows that we can replace  ``$x \xbe X$''
equivalently by
``$x \xbe \xbm (X)$''  in the definition of $H(U,x)_{ \xba +1},$ as was
done in Definition 3.7.5
in [Sch04].

(8.3)

Suppose $( \xbm Cum \xba )$ fails, we show that then so does (HUx). As $(
\xbm Cum \xba )$ fails, for
all $ \xbb \xck \xba $ $ \xbm (X_{ \xbb }) \xcc U \xcv \xcV \{X_{ \xbg }:
\xbg < \xbb \},$ but there is $x \xbe \xcS \{X_{ \xbg }: \xbg \xck \xba \}
\xcs \xbm (U),$
$x \xce \xbm (X_{ \xba }).$ Thus for all $ \xbb \xck \xba $ $ \xbm (X_{
\xbb }) \xcc X_{ \xbb } \xcc H(U,x),$ moreover $x \xbe \xbm (U),$
$x \xbe X_{ \xba }- \xbm (X_{ \xba }),$ but $ \xbm (X_{ \xba }) \xcc
H(U,x),$ so (HUx) fails.

$ \xcz $
\\[3ex]

We turn to $H(U).$

\bfa

$\hspace{0.01em}$

% (+++*** Orig. No.:  Fact 2.2: )

Let $A,$ $X,$ $U,$ $U',$ $Y$ and all $A_{i}$ be in $ \xdy.$

(1) $( \xbm \xcc )$ and (HU) entail:

(1.1) $( \xbm PR)$

(1.2) $( \xbm Cum)$

(2) (HU) $+$ $( \xcv )$ $ \xcp $ (HUx)

(3) $( \xbm \xcc )$ and $( \xbm PR)$ entail:

(3.1) $A= \xcV \{A_{i}:i \xbe I\}$ $ \xcp $ $ \xbm (A) \xcc \xcV \{ \xbm
(A_{i}):i \xbe I\},$

(3.2) $U \xcc H(U),$ and $U \xcc U' \xcp H(U) \xcc H(U' ),$

(3.3) $ \xbm (U \xcv Y)-H(U) \xcc \xbm (Y)$ - if $ \xbm (U \xcv Y)$ is
defined, in particular, if $( \xcv )$
holds.

(4) $( \xcv ),$ $( \xbm \xcc ),$ $( \xbm PR),$ $( \xbm CUM)$ entail:

(4.1) $H(U)=H_{1}(U)$

(4.2) $U \xcc A,$ $ \xbm (A) \xcc H(U)$ $ \xcp $ $ \xbm (A) \xcc U,$

(4.3) $ \xbm (Y) \xcc H(U)$ $ \xcp $ $Y \xcc H(U)$ and $ \xbm (U \xcv Y)=
\xbm (U),$

(4.4) $x \xbe \xbm (U),$ $x \xbe Y- \xbm (Y)$ $ \xcp $ $Y \xcC H(U)$ (and
thus (HU)),

(4.5) $Y \xcC H(U)$ $ \xcp $ $ \xbm (U \xcv Y) \xcC H(U).$

(5) $( \xcv ),$ $( \xbm \xcc ),$ (HU) entail

(5.1) $H(U)=H_{1}(U)$

(5.2) $U \xcc A,$ $ \xbm (A) \xcc H(U)$ $ \xcp $ $ \xbm (A) \xcc U,$

(5.3) $ \xbm (Y) \xcc H(U)$ $ \xcp $ $Y \xcc H(U)$ and $ \xbm (U \xcv Y)=
\xbm (U),$

(5.4) $x \xbe \xbm (U),$ $x \xbe Y- \xbm (Y)$ $ \xcp $ $Y \xcC H(U),$

(5.5) $Y \xcC H(U)$ $ \xcp $ $ \xbm (U \xcv Y) \xcC H(U).$

\efa

\paragraph{
Proof:
}

$\hspace{0.01em}$

(1.1) By (HU), if $ \xbm (Y) \xcc H(U),$ then $ \xbm (U) \xcs Y \xcc \xbm
(Y).$ But, if $Y \xcc U,$ then
$ \xbm (Y) \xcc H(U)$ by $( \xbm \xcc ).$

(1.2) Let $ \xbm (U) \xcc X \xcc U.$ Then by (1.1) $ \xbm (U)= \xbm (U)
\xcs X \xcc \xbm (X).$ By prerequisite,
$ \xbm (U) \xcc U \xcc H(X),$ so $ \xbm (X)= \xbm (X) \xcs U \xcc \xbm
(U)$ by $( \xbm \xcc ).$

(2) By (1.2), (HU) entails $( \xbm Cum),$ so by $( \xcv )$ and Fact 2.1,
(5.2) $( \xbm Cum \xca )$
holds, so by Fact 2.1, (8.2) (HUx) holds.

(3.1) $ \xbm (A) \xcs A_{j} \xcc \xbm (A_{j}) \xcc \xcV \xbm (A_{i}),$ so
by $ \xbm (A) \xcc A= \xcV A_{i}$ $ \xbm (A) \xcc \xcV \xbm (A_{i}).$

(3.2) trivial.

(3.3) $ \xbm (U \xcv Y)-H(U)$ $ \xcc_{(3.2)}$ $ \xbm (U \xcv Y)-U$ $ \xcc
$ (by $( \xbm \xcc )$ and (3.1))
$ \xbm (U \xcv Y) \xcs Y$ $ \xcc_{( \xbm PR)}$ $ \xbm (Y).$

(4.1) We show that, if $X \xcc H_{2}(U),$ then $X \xcc H_{1}(U),$ more
precisely, if $ \xbm (X) \xcc H_{1}(U),$
then already $X \xcc H_{1}(U),$ so the construction stops already at
$H_{1}(U).$
Suppose then $ \xbm (X) \xcc \xcV \{Y: \xbm (Y) \xcc U\},$ and let $A:=X
\xcv U.$ We show that $ \xbm (A) \xcc U,$ so
$X \xcc A \xcc H_{1}(U).$ Let $a \xbe \xbm (A).$ By (3.1), $ \xbm (A) \xcc
\xbm (X) \xcv \xbm (U).$ If $a \xbe \xbm (U) \xcc U,$ we
are done. If $a \xbe \xbm (X),$ there is $Y$ s.t. $ \xbm (Y) \xcc U$ and
$a \xbe Y,$ so $a \xbe \xbm (A) \xcs Y.$
By Fact 2.1, (5.1.3), we have for $Y$ s.t. $ \xbm (Y) \xcc U$ and $U \xcc
A$ $ \xbm (A) \xcs Y \xcc \xbm (U).$
Thus $a \xbe \xbm (U),$ and we are done again.

(4.2) Let $U \xcc A,$ $ \xbm (A) \xcc H(U)=H_{1}(U)$ by (4.1). So $ \xbm
(A)$ $=$ $ \xcV \{ \xbm (A) \xcs Y: \xbm (Y) \xcc U\}$ $ \xcc $
$ \xbm (U)$ $ \xcc $ $U,$ again by Fact 2.1, (5.1.3).

(4.3) Let $ \xbm (Y) \xcc H(U),$ then by $ \xbm (U) \xcc H(U)$ and (3.1)
$ \xbm (U \xcv Y) \xcc \xbm (U) \xcv \xbm (Y) \xcc H(U),$ so by (4.2) $
\xbm (U \xcv Y) \xcc U$ and $U \xcv Y \xcc H(U).$
Moreover, $ \xbm (U \xcv Y) \xcc U \xcc U \xcv Y$ $ \xcp_{( \xbm CUM)}$ $
\xbm (U \xcv Y)= \xbm (U).$

(4.4) If not, $Y \xcc H(U),$ so $ \xbm (Y) \xcc H(U),$ so $ \xbm (U \xcv
Y)= \xbm (U)$ by (4.3),
but $x \xbe Y- \xbm (Y)$ $ \xcp_{( \xbm PR)}$ $x \xce \xbm (U \xcv Y)=
\xbm (U),$ $contradiction.$

(4.5) $ \xbm (U \xcv Y) \xcc H(U)$ $ \xcp_{(4.3)}$ $U \xcv Y \xcc H(U).$

(5) Trivial by (1) and (4).

$ \xcz $
\\[3ex]

(5) is just noted for the convenience of the reader. It will be used for
the
proof of Fact 2.3.

Thus, in the presence of $( \xcv )$ $H(U,x)$ can be simplified to $H(U),$
which is
constructed in one single step, and is independent from $x.$ Of course,
then
$H(U,x) \xcc H(U).$

\subsection{
Representation by smooth preferential structures
}

% (+++*** Orig.:   (2.3)  Representation by smooth preferential structures )
%  (2.3.1)  The not necessarily transitive case  (most revised -> DCB 3.3.1)
%  (2.3.1)  The not necessarily transitive case  (most revised -> DCB 3.3.1)
% %
% ==========================================================================
\subsubsection{The not necessarily transitive case}

We adapt Proposition 3.7.15 in [Sch04] and its proof. All we need is (HUx)
and
$( \xbm \xcc ).$ We modify the proof of Remark 3.7.13 (1) in [Sch04] (now
Remark 2.4) so we
will not need $( \xcs )$ any more.
We will give the full proof, although its essential elements have already
been
published, for three reasons: First, the new version will need less
prerequisites than the old proof does (closure under finite intersections
is
not needed any more, and replaced by (HUx)). Second, we will more clearly
separate the requirements to do the construction from the construction
itself,
thus splitting the proof neatly into two parts.

We show how to work with $( \xbm \xcc )$ and (HUx) only. Thus, once we
have shown $( \xbm \xcc )$
and (HUx), we have finished the substantial side, and enter the
administrative
part, which will not use any prerequisites about domain closure any more.
At the same time, this gives a uniform proof of the difficult part for the
case
with and without $( \xcv ),$ in the former case we can even work with the
stronger
$H(U).$ The easy direction of the former parts needs a proof of the
stronger
$H(U),$ but this is easy.

Note that, by Fact 2.1, (8.3) and (7), (HUx) entails $( \xbm PR),$ so we
can use it
in our context, where (HUx) will be the central property.

\bfa

$\hspace{0.01em}$

% (+++*** Orig. No.:  Fact 2.3: )

(1) $x \xbe \xbm (Y),$ $ \xbm (Y) \xcc H(U,x)$ $ \xcp $ $Y \xcc H(U,x),$

(2) (HUx) holds in all smooth models.

\efa

\paragraph{
Proof:
}

$\hspace{0.01em}$

(1) Trivial by definition.

(2) Suppose not. So let $x \xbe \xbm (U),$ $x \xbe Y- \xbm (Y),$ $ \xbm
(Y) \xcc H(U,x).$
By smoothness, there is $x_{1} \xbe \xbm (Y),$ $x \xee x_{1},$
and let $ \xbk_{1}$ be the least $ \xbk $ s.t. $x_{1} \xbe H(U,x)_{
\xbk_{1}}.$ $ \xbk_{1}$ is not a
limit, and $x_{1} \xbe U'_{x_{1}}- \xbm (U'_{x_{1}})$ with $x \xbe \xbm
(U'_{x_{1}})$ for some $U'_{x_{1}},$ so as $x_{1} \xce \xbm (U'_{x_{1}}),$
there must be (by smoothness) some other
$x_{2} \xbe \xbm (U'_{x_{1}}) \xcc H(U,x)_{ \xbk_{1}-1}$ with $x \xee
x_{2}.$ Continue with $x_{2},$ we thus construct
a descending chain of ordinals, which cannot be infinite, so there must be
$x_{n} \xbe \xbm (U'_{x_{n}}) \xcc U,$ $x \xee x_{n},$ contradicting
minimality of $x$ in $U.$
(More precisely, this works for all copies of x.)
$ \xcz $
\\[3ex]

We first show two basic facts and then turn to the main result,
Proposition
2.6.

\bd

$\hspace{0.01em}$

% (+++*** Orig. No.:  Definition 2.3: )

For $x \xbe Z,$ let $ \xdw_{x}:=\{ \xbm (Y)$: $Y \xbe \xdy $ $ \xcu $ $x
\xbe Y- \xbm (Y)\},$ $ \xbG_{x}:= \xbP \xdw_{x}$, and $K:=\{x \xbe Z$: $
\xcE X \xbe \xdy.x \xbe \xbm (X)\}.$

\ed

\br

$\hspace{0.01em}$

% (+++*** Orig. No.:  Remark 2.4: )

(1) $x \xbe K$ $ \xcp $ $ \xbG_{x} \xEd \xCQ,$

(2) $g \xbe \xbG_{x}$ $ \xcp $ $ran(g) \xcc K.$

\er

\paragraph{
Proof:
}

$\hspace{0.01em}$

(1) We give two proofs, the first uses $( \xbm Cum0),$ the second the
stronger
(HUx).

(a) We have to show that $Y \xbe \xdy,$ $x \xbe Y- \xbm (Y)$ $ \xcp $ $
\xbm (Y) \xEd \xCQ.$ Suppose then $x \xbe \xbm (X),$
this exists, as $x \xbe K,$ and $ \xbm (Y)= \xCQ,$ so $ \xbm (Y) \xcc X,$
$x \xbe Y,$ so by $( \xbm Cum0)$ $x \xbe \xbm (Y).$

(b) $ \xbm (Y)= \xCQ $ $ \xcp $ $Y \xcc H(U,x),$ contradicting $x \xbe Y-
\xbm (Y).$

(2) By definition, $ \xbm (Y) \xcc K$ for all $Y \xbe \xdy.$
$ \xcz $
\\[3ex]

\bc

$\hspace{0.01em}$

% (+++*** Orig. No.:  Claim 2.5: )

Let $U \xbe \xdy,$ $x \xbe K.$ Then

(1) $x \xbe \xbm (U)$ $ \xcr $ $x \xbe U$ $ \xcu $ $ \xcE f \xbe
\xbG_{x}.ran(f) \xcs U= \xCQ,$

(2) $x \xbe \xbm (U)$ $ \xcr $ $x \xbe U$ $ \xcu $ $ \xcE f \xbe
\xbG_{x}.ran(f) \xcs H(U,x)= \xCQ.$

\ec

\paragraph{
Proof:
}

$\hspace{0.01em}$

(1)

Case 1: $ \xdw_{x}= \xCQ,$ thus $ \xbG_{x}=\{ \xCQ \}.$

``$ \xcp $'': Take $f:= \xCQ.$

``$ \xcq $'': $x \xbe U \xbe \xdy,$ $ \xdw_{x}= \xCQ $ $ \xcp $ $x \xbe
\xbm (U)$ by definition of $ \xdw_{x}.$

Case 2: $ \xdw_{x} \xEd \xCQ.$

``$ \xcp $'': Let $x \xbe \xbm (U) \xcc U.$ By (HUx), if $Y \xbe
\xdw_{x},$ then $ \xbm (Y)-H(U,x) \xEd \xCQ.$

``$ \xcq $'': If $x \xbe U- \xbm (U),$ $ \xbm (U) \xbe \xdw_{x}$,
moreover $ \xbG_{x} \xEd \xCQ $ by Remark 2.4, (1) and
thus $ \xbm (U) \xEd \xCQ,$ so $ \xcA f \xbe \xbG_{x}.ran(f) \xcs U \xEd
\xCQ.$

(2): The proof is verbatim the same as for (1).

$ \xcz $ (Claim 2.5)
\\[3ex]

The following Proposition 2.6 is the main result of Section 2.3.1 and
shows how to characterize smooth structures in the absence of closure
under
finite unions. The strategy of the proof follows closely the proof of
Proposition 3.3.4 in [Sch04].

\bp

$\hspace{0.01em}$

% (+++*** Orig. No.:  Proposition 2.6: )

Let $ \xbm: \xdy \xcp \xdp (Z).$
Then there is a $ \xdy -$smooth preferential structure $ \xdz,$ s.t. for
all $X \xbe \xdy $
$ \xbm (X)= \xbm_{ \xdz }(X)$ iff $ \xbm $ satisfies $( \xbm \xcc )$ and
(HUx) above.

\ep

\paragraph{
Proof:
}

$\hspace{0.01em}$

``$ \xcp $'' (HUx) was shown in Fact 2.3.

Outline of ``$ \xcq $'': We first define a structure $ \xdz $ which
represents $ \xbm,$ but is
not necessarily $ \xdy -$smooth, refine it to $ \xdz ' $ and show that $
\xdz ' $ represents $ \xbm $
too, and that $ \xdz ' $ is $ \xdy -$smooth.

In the structure $ \xdz ',$ all pairs destroying smoothness in $ \xdz $
are successively
repaired, by adding minimal elements: If $<y,j>$ is not minimal, and has
no minimal
$<x,i>$ below it, we just add one such $<x,i>.$ As the repair process
might itself
generate such ``bad'' pairs, the process may have to be repeated infinitely
often.
Of course, one has to take care that the representation property is
preserved.

\bcs

$\hspace{0.01em}$

% (+++*** Orig. No.:  Construction 2.1: )

(Construction of $ \xdz )$

Let $ \xdx $ $:=$ $\{<x,g>$: $x \xbe K,$ $g \xbe \xbG_{x}\},$ $<x',g' >
\xeb <x,g>$ $: \xcr $ $x' \xbe ran(g),$ $ \xdz:=< \xdx, \xeb >.$

\ecs

\bc

$\hspace{0.01em}$

% (+++*** Orig. No.:  Claim 2.7: )

$ \xcA U \xbe \xdy. \xbm (U)= \xbm_{ \xdz }(U)$

\ec

\paragraph{
Proof:
}

$\hspace{0.01em}$

Case 1: $x \xce K.$ Then $x \xce \xbm (U)$ and $x \xce \xbm_{ \xdz }(U).$

Case 2: $x \xbe K.$

By Claim 2.5, (1) it suffices to show that for all $U \xbe \xdy $
$x \xbe \xbm_{ \xdz }(U)$ $ \xcr $ $x \xbe U$ $ \xcu $ $ \xcE f \xbe
\xbG_{x}.ran(f) \xcs U= \xCQ.$
Fix $U \xbe \xdy.$

``$ \xcp $'': $x \xbe \xbm_{ \xdz }(U)$ $ \xcp $ ex. $<x,f>$ minimal in
$ \xdx \xex U,$ thus $x \xbe U$ and there is no
$<x',f' > \xeb <x,f>,$ $x' \xbe U,$ $x' \xbe K.$ But if $x' \xbe K,$ then
by Remark 2.4, (1), $ \xbG_{x' } \xEd \xCQ,$
so we find suitable $f'.$ Thus, $ \xcA x' \xbe ran(f).x' \xce U$ or $x'
\xce K.$ But $ran(f) \xcc K,$ so
$ran(f) \xcs U= \xCQ.$

``$ \xcq $'': If $x \xbe U,$ $f \xbe \xbG_{x}$ s.t. $ran(f) \xcs U= \xCQ
,$ then $<x,f>$ is minimal in $ \xdx \xex U.$
$ \xcz $ (Claim 2.7)
\\[3ex]

We now construct the refined structure $ \xdz '.$

\bcs

$\hspace{0.01em}$

% (+++*** Orig. No.:  Construction 2.2: )

(Construction of $ \xdz ' )$

$ \xbs $ is called x-admissible sequence iff

1. $ \xbs $ is a sequence of length $ \xck \xbo,$ $ \xbs =\{ \xbs_{i}:i
\xbe \xbo \},$

2. $ \xbs_{o} \xbe \xbP \{ \xbm (Y)$: $Y \xbe \xdy $ $ \xcu $ $x \xbe Y-
\xbm (Y)\},$

3. $ \xbs_{i+1} \xbe \xbP \{ \xbm (X)$: $X \xbe \xdy $ $ \xcu $ $x \xbe
\xbm (X)$ $ \xcu $ $ran( \xbs_{i}) \xcs X \xEd \xCQ \}.$

By 2., $ \xbs_{0}$ minimizes $x,$ and by 3., if $x \xbe \xbm (X),$ and
$ran( \xbs_{i}) \xcs X \xEd \xCQ,$ i.e. we
have destroyed minimality of $x$ in $X,$ $x$ will be above some $y$
minimal in $X$ to
preserve smoothness.

Let $ \xbS_{x}$ be the set of x-admissible sequences, for $ \xbs \xbe
\xbS_{x}$ let $ \wt{ \xbs }:= \xcV \{ran( \xbs_{i}):i \xbe \xbo \}.$
Note that by Remark 2.4, (1), $ \xbS_{x} \xEd \xCQ,$ if $x \xbe K$ (this
does $ \xbs_{0},$ $ \xbs_{i+1}$ is
trivial as by prerequisite $ \xbm (X) \xEd \xCQ ).$

Let $ \xdx ' $ $:=$ $\{<x, \xbs >$: $x \xbe K$ $ \xcu $ $ \xbs \xbe
\xbS_{x}\}$ and $<x', \xbs ' > \xeb ' <x, \xbs >$ $: \xcr $ $x' \xbe \wt{
\xbs }$.
Finally, let $ \xdz ':=< \xdx ', \xeb ' >,$ and $ \xbm ':= \xbm_{ \xdz
' }.$

\ecs

It is now easy to show that $ \xdz ' $ represents $ \xbm,$ and that $
\xdz ' $ is smooth.
For $x \xbe \xbm (U),$ we construct a special x-admissible sequence $
\xbs^{x,U}$ using the
properties of $H(U,x)$ as described at the beginning of this section.

\bc

$\hspace{0.01em}$

% (+++*** Orig. No.:  Claim 2.8: )

For all $U \xbe \xdy $ $ \xbm (U)= \xbm_{ \xdz }(U)= \xbm ' (U).$

\ec

\paragraph{
Proof:
}

$\hspace{0.01em}$

If $x \xce K,$ then $x \xce \xbm_{ \xdz }(U),$ and $x \xce \xbm ' (U)$ for
any $U.$ So assume $x \xbe K.$ If $x \xbe U$ and
$x \xce \xbm_{ \xdz }(U),$ then for all $<x,f> \xbe \xdx,$ there is $<x'
,f' > \xbe \xdx $ with $<x',f' > \xeb <x,f>$ and
$x' \xbe U.$ Let now $<x, \xbs > \xbe \xdx ',$ then $<x, \xbs_{0}> \xbe
\xdx,$ and let $<x',f' > \xeb <x, \xbs_{0}>$ in $ \xdz $ with
$x' \xbe U.$ As $x' \xbe K,$ $ \xbS_{x' } \xEd \xCQ,$ let $ \xbs ' \xbe
\xbS_{x' }$. Then $<x', \xbs ' > \xeb ' <x, \xbs >$ in $ \xdz '.$ Thus
$x \xce \xbm ' (U).$
Thus, for all $U \xbe \xdy,$ $ \xbm ' (U) \xcc \xbm_{ \xdz }(U)= \xbm
(U).$

It remains to show $x \xbe \xbm (U) \xcp x \xbe \xbm ' (U).$

Assume $x \xbe \xbm (U)$ (so $x \xbe K),$ $U \xbe \xdy,$ we will
construct minimal $ \xbs,$ i.e. show that
there is $ \xbs^{x,U} \xbe \xbS_{x}$ s.t. $ \wt{ \xbs^{x,U}} \xcs U= \xCQ
.$ We construct this $ \xbs^{x,U}$ inductively, with the
stronger property that $ran( \xbs^{x,U}_{i}) \xcs H(U,x)= \xCQ $ for all
$i \xbe \xbo.$

$ \ul{ \xbs^{x,U}_{0}:}$

$x \xbe \xbm (U),$ $x \xbe Y- \xbm (Y)$ $ \xcp $ $ \xbm (Y)-H(U,x) \xEd
\xCQ $ by (HUx).
Let $ \xbs^{x,U}_{0}$ $ \xbe $ $ \xbP \{ \xbm (Y)-H(U,x):$ $Y \xbe \xdy,$
$x \xbe Y- \xbm (Y)\},$ so $ran( \xbs^{x,U}_{0}) \xcs H(U,x)= \xCQ.$

$ \ul{ \xbs^{x,U}_{i} \xcp \xbs^{x,U}_{i+1}:}$

By the induction hypothesis, $ran( \xbs^{x,U}_{i}) \xcs H(U,x)= \xCQ.$
Let $X \xbe \xdy $ be s.t. $x \xbe \xbm (X),$
$ran( \xbs^{x,U}_{i}) \xcs X \xEd \xCQ.$ Thus $X \xcC H(U,x),$ so $ \xbm
(X)-H(U,x) \xEd \xCQ $ by Fact 2.3, (1).
Let $ \xbs^{x,U}_{i+1}$ $ \xbe $ $ \xbP \{ \xbm (X)-H(U,x):$ $X \xbe \xdy
,$ $x \xbe \xbm (X),$ $ran( \xbs^{x,U}_{i}) \xcs X \xEd \xCQ \},$ so $ran(
\xbs^{x,U}_{i+1}) \xcs H(U,x)= \xCQ.$
The construction satisfies the x-admissibility condition.
$ \xcz $
\\[3ex]

It remains to show:

\bc

$\hspace{0.01em}$

% (+++*** Orig. No.:  Claim 2.9: )

$ \xdz ' $ is $ \xdy -$smooth.

\ec

\paragraph{
Proof:
}

$\hspace{0.01em}$

Let $X \xbe \xdy,$ $<x, \xbs > \xbe \xdx ' \xex X.$

Case 1, $x \xbe X- \xbm (X):$ Then $ran( \xbs_{0}) \xcs \xbm (X) \xEd \xCQ
,$ let $x' \xbe ran( \xbs_{0}) \xcs \xbm (X).$ Moreover,
$ \xbm (X) \xcc K.$ Then for all $<x', \xbs ' > \xbe \xdx ' $ $<x', \xbs
' > \xeb <x, \xbs >.$ But $<x', \xbs^{x',X}>$ as
constructed in the proof of Claim 2.8 is minimal in $ \xdx ' \xex X.$

Case 2, $x \xbe \xbm (X)= \xbm_{ \xdz }(X)= \xbm ' (X):$ If $<x, \xbs >$
is minimal in $ \xdx ' \xex X,$ we are done.
So suppose there is $<x', \xbs ' > \xeb <x, \xbs >,$ $x' \xbe X.$ Thus
$x' \xbe \wt{ \xbs }.$ Let
$x' \xbe ran( \xbs_{i}).$ So $x \xbe \xbm (X)$ and $ran( \xbs_{i}) \xcs X
\xEd \xCQ.$ But
$ \xbs_{i+1} \xbe \xbP \{ \xbm (X' )$: $X' \xbe \xdy $ $ \xcu $ $x \xbe
\xbm (X' )$ $ \xcu $ $ran( \xbs_{i}) \xcs X' \xEd \xCQ \},$ so $X$ is one
of the $X',$
moreover $ \xbm (X) \xcc K,$ so there is $x'' \xbe \xbm (X) \xcs ran(
\xbs_{i+1}) \xcs K,$ so for all $<x'', \xbs '' > \xbe \xdx ' $
$<x'', \xbs '' > \xeb <x, \xbs >.$ But again $<x'', \xbs^{x'',X}>$ as
constructed in the proof of Claim 2.8
is minimal in $ \xdx ' \xex X.$

$ \xcz $ (Claim 2.9 and Proposition 2.6)
\\[3ex]

We conclude this section by showing that we cannot improve substantially.

\bp

$\hspace{0.01em}$

% (+++*** Orig. No.:  Proposition 2.10: )

There is no fixed size characterization of $ \xbm -$functions which are
representable
by smooth structures, if the domain is not closed under finite unions.

\ep

\paragraph{
Proof:
}

$\hspace{0.01em}$

Suppose we have a fixed size characterization, which allows to distinguish
$ \xbm -$functions on domains which are not necessarily closed under
finite unions,
and which can be represented by smooth structures, from those which cannot
be
represented in this way. Let the characterization have $ \xba $ parameters
for sets,
and consider Example 2.1 with $ \xbk = \xbb +1,$ $ \xbb > \xba $ (as a
cardinal). This structure
cannot be represented, as $( \xbm Cum \xbk )$ fails. As we have only $
\xba $ parameters, at
least one of the $X_{ \xbg }$ is not mentioned, say $X_{ \xbd }.$ Wlog, we
may assume that
$ \xbd = \xbd ' +1.$ We change now the structure, and erase one pair of
the relation,
$x_{ \xbd } \xeb x_{ \xbd +1}.$ Thus, $ \xbm (X_{ \xbd })=\{c,x_{ \xbd
},x_{ \xbd +1}\}.$ But now we cannot go any more from
$X_{ \xbd ' }$ to $X_{ \xbd ' +1}=X_{ \xbd },$ as $ \xbm (X_{ \xbd }) \xcC
X_{ \xbd ' }.$ Consequently, the only chain showing
that $( \xbm Cum \xca )$ fails is interrupted - and we have added no new
possibilities,
as inspection of cases shows. $(x_{ \xbd +1}$ is now globally minimal, and
increasing
$ \xbm (X)$ cannot introduce new chains, only interrupt chains.) Thus, $(
\xbm Cum \xca )$ holds
in the modified example, and it is thus representable by a smooth
structure, as
above proposition shows. As we did not touch any of the parameters, the
truth
value of the characterization is unchanged, which was negative. So the
``characterization'' cannot be correct.
$ \xcz $
\\[3ex]
%  (2.3.2)  The transitive case  (-> DCB 3.3.2, part revised)
%  (2.3.2)  The transitive case  (-> DCB 3.3.2, part revised)
% %
% ===========================================================
\subsubsection{The transitive case}

Unfortunately, $( \xbm Cumt \xca )$ is a necessary but not sufficient
condition for
smooth transitive structures, as can be seen in the following example.

\be

$\hspace{0.01em}$

% (+++*** Orig. No.:  Example 2.2: )

We assume no closure whatever.

$U:=\{u_{1},u_{2},u_{3},u_{4}\},$ $ \xbm (U):=\{u_{3},u_{4}\}$

$Y_{1}:=\{u_{4},v_{1},v_{2},v_{3},v_{4}\},$ $ \xbm
(Y_{1}):=\{v_{3},v_{4}\}$

$Y_{2,1}:=\{u_{2},v_{2},v_{4}\},$ $ \xbm (Y_{2,1}):=\{u_{2},v_{2}\}$

$Y_{2,2}:=\{u_{1},v_{1},v_{3}\},$ $ \xbm (Y_{2,2}):=\{u_{1},v_{1}\}$

For no A,B $ \xbm (A) \xcc B$ $(A \xEd B),$ so the prerequisite of $( \xbm
Cumt \xba )$ is false,
and $( \xbm Cumt \xba )$ holds, but there is no smooth transitive
representation possible:
if $u_{4} \xee v_{3},$ then $Y_{2,2}$ makes this impossible, if $u_{4}
\xee v_{4},$ then $Y_{2,1}$ makes this
impossible.

$ \xcz $
\\[3ex]

\ee

\br

$\hspace{0.01em}$

% (+++*** Orig. No.:  Remark 2.11: )

(1)
The situation does not change when we have copies, the same argument will
still
work: There is a U-minimal copy $<u_{4},i>,$ by smoothness and $Y_{1},$
there must be a
$Y_{1}-$minimal copy, e.g. $<v_{3},j> \xeb <u_{4},i>.$ By smoothness and
$Y_{2,2},$ there must be a
$Y_{2,2}-$minimal $<u_{1},k>$ or $<v_{1},l>$ below $<v_{3},j>.$ But
$v_{1}$ is in $Y_{1},$ contradicting
minimality of $<v_{3},j>,$ $u_{1}$ is in $U,$ contadicting minimality of
$<u_{4},i>$ by
transitivity. If we choose $<v_{4},j>$ minimal below $<u_{4},i>,$ we will
work with
$Y_{2,1}$ instead of $Y_{2,2}.$

(2)
We can also close under arbitrary intersections, and the example will
still
work: We have to consider $U \xcs Y_{1},$ $U \xcs Y_{2,1},$ $U \xcs
Y_{2,2},$ $Y_{2,1} \xcs Y_{2,2},$ $Y_{1} \xcs Y_{2,1},$ $Y_{1} \xcs
Y_{2,2},$
there are no further intersections to consider. We may assume $ \xbm
(A)=A$ for all
these intersections (working with copies). But then $ \xbm (A) \xcc B$
implies $ \xbm (A)=A$
for all sets, and all $( \xbm Cumt \xba )$ hold again trivially.

(3) If we had finite unions, we could form $A:=U \xcv Y_{1} \xcv Y_{2,1}
\xcv Y_{2,2},$ then $ \xbm (A)$ would
have to be a subset of $\{u_{3}\}$ by $( \xbm PR),$ so by $( \xbm CUM)$
$u_{4} \xce \xbm (U),$ a contradiction.
Finite unions allow us to ``look ahead'', without $( \xcv ),$ we see
desaster only at
the end - and have to backtrack, i.e. try in our example $Y_{2,1},$ once
we have seen
impossibility via $Y_{2,2},$ and discover impossibility again at the end.

\er

We discuss and define now an analogon to (HU) or (HUx), condition $( \xbm
\xbt ),$
defined in Definition 2.4.
%  The different possible cases
%  The different possible cases
% %
% =============================
\paragraph{The different possible cases}

The problem is to minimize an element, already
minimal in a finite number of sets, in a new set, without destroying
previous
minimality. We have to examine the way the new minimal element is chosen.

\paragraph{
(I) Going forward
}

$\hspace{0.01em}$

Let $Y_{0}, \Xl,Y_{n}$ be treated and $x_{n} \xbe \xbm (Y_{n}).$ (We can
argue without copies, as we may
assume that we have chosen a minimal copy.)

Treating $Y_{n+1}:$

For all $Y_{n+1}$ s.t. $x_{n} \xbe Y_{n+1},$ we have to treat $Y_{n+1}$
and choose $x_{n+1} \xbe \xbm (Y_{n+1}).$

If $x_{n} \xce Y_{n+1},$ there is nothing to do: $Y_{n+1}$ is not to be
considered.

Case 1: $x_{n} \xbe \xbm (Y_{n+1}),$ then we can either

Case 1.1: leave it as it is, i.e. $x_{n+1}:=x_{n},$

Case 1.2: minimize it by another $x_{n+1} \xbe \xbm (Y_{n+1}),$ outside
$Y_{0} \xcv  \Xl  \xcv Y_{n}$ as we do
not want to destroy previous minimality - assuming that $ \xbm (Y_{n+1})
\xcC Y_{0} \xcv  \Xl  \xcv Y_{n}.$

Case 2: $x_{n} \xbe Y_{n+1}- \xbm (Y_{n+1}),$ then we have to

minimize it by another $x_{n+1} \xbe \xbm (Y_{n+1}),$ outside $Y_{0} \xcv
\Xl  \xcv Y_{n}$ as we do not want to
destroy previous minimality - assuming that $ \xbm (Y_{n+1}) \xcC Y_{0}
\xcv  \Xl  \xcv Y_{n}.$

We tentatively write this down as follows:

- $(Y_{n+1},x_{n+1} \xbe \xbm (Y_{n+1}),m)$ - if we modified it, i.e.
$x_{n+1} \xEd x_{n},$ $m$ for modify, and

- $(Y_{n+1},x_{n+1} \xbe \xbm (Y_{n+1}),c)$ - if $x_{n+1}=x_{n},$ $c$ for
constant.

We have to do this for all $Y_{n+1}$ s.t. $x_{n} \xbe Y_{n+1}.$

We continue to go forward one further step.

Treating $Y_{n+2}:$

For all $Y_{n+2}$ s.t. $x_{n+1} \xbe Y_{n+2},$ we treat $Y_{n+2}$ by
choosing $x_{n+2} \xbe \xbm (Y_{n+2}).$

Let $x_{n+1} \xbe Y_{n+2}.$

Case 1: $x_{n+1} \xbe \xbm (Y_{n+2})$

Case 1.1: We leave $x_{n+1}$ as it is, $x_{n+2}:=x_{n+1},$ so the next
element is
$(Y_{n+2},x_{n+2} \xbe \xbm (Y_{n+2}),c).$

Case 1.2: We try to modify.

If $ \xbm (Y_{n+2}) \xcC Y_{0} \xcv  \Xl  \xcv Y_{n} \xcv Y_{n+1},$ then
we can
choose any $x_{n+2} \xbe \xbm (Y_{n+2})-Y_{0} \xcv  \Xl  \xcv Y_{n} \xcv
Y_{n+1},$ and the successor is
$(Y_{n+2},x_{n+2} \xbe \xbm (Y_{n+2}),m).$

If $ \xbm (Y_{n+2}) \xcc Y_{0} \xcv  \Xl  \xcv Y_{n} \xcv Y_{n+1},$ then
this is impossible, and we have
to work with Case 1.1.

Case 2: $x_{n+1} \xbe Y_{n+2}- \xbm (Y_{n+2})$

We have to modify.

If $ \xbm (Y_{n+2}) \xcC Y_{0} \xcv  \Xl  \xcv Y_{n} \xcv Y_{n+1},$ then
we can
choose any $x_{n+2} \xbe \xbm (Y_{n+2})-Y_{0} \xcv  \Xl  \xcv Y_{n} \xcv
Y_{n+1},$ and the successor is
$(Y_{n+2},x_{n+2} \xbe \xbm (Y_{n+2}),m)$

If $ \xbm (Y_{n+2}) \xcc Y_{0} \xcv  \Xl  \xcv Y_{n} \xcv Y_{n+1},$ then
this is impossible, but now we
have no alternative, thus already the choice of $x_{n+1} \xbe Y_{n+2}-
\xbm (Y_{n+2})$
was impossible, as we then have to consider $Y_{n+2}.$

This is the heart of the problem: a subsequent step (considering
$Y_{n+2})$
may show that a previous choice $(x_{n+1})$ was impossible, so we have to
backtrack.

Of course, this does not only eliminate this particular $x_{n+1},$ but any
other $x_{n+1} \xbe Y_{n+2}- \xbm (Y_{n+2}),$ too. But it does not concern
any
$x_{n+1} \xbe \xbm (Y_{n+2}),$ as we then have alternative 1.1 above. We
also note that
it is unimportant how we obtained the previous $x_{n+1}$ - by modification
or
staying constant.

\paragraph{
(II) Backtracking
}

$\hspace{0.01em}$

We discuss now the repercussions of such, in hindsight, impossible
choices.

Treating $Y_{n+1}:$

Let $ \xbm (Y_{n+2}) \xcc Y_{0} \xcv  \Xl  \xcv Y_{n} \xcv Y_{n+1}$ and
$x_{n} \xbe Y_{n+1},$ so we have to treat $Y_{n+1},$ and
choose $x_{n+1} \xbe \xbm (Y_{n+1}).$ If we choose $x_{n+1} \xbe Y_{n+2}-
\xbm (Y_{n+2}),$ then the next step
will show us that this is impossible, so we have to choose $x_{n+1}$
outside of
$Y_{n+2}- \xbm (Y_{n+2}).$ If $x_{n} \xbe \xbm (Y_{n+1}),$ and $x_{n} \xce
Y_{n+2}- \xbm (Y_{n+2}),$ we can choose $x_{n+1}:=x_{n}.$
If $x_{n} \xce \xbm (Y_{n+1}),$ we have to choose $x_{n+1} \xbe \xbm
(Y_{n+1})-(Y_{0} \xcv  \Xl  \xcv Y_{n})-(Y_{n+2}- \xbm (Y_{n+2})).$

So there is an additional problem here: if we can choose $x_{n+1}:=x_{n},$
we have
more liberty, we need not necessarily avoid $Y_{0} \xcv  \Xl  \xcv Y_{n}$
(if we were constant
all the time, we need not avoid any of the previous $Y_{i}).$ For this
reason,
the following simplification will not work: Suppose there is a cover of
$ \xbm (Y_{n+1})-(Y_{0} \xcv  \Xl  \xcv Y_{n})$ by such $Y_{n+2,i}- \xbm
(Y_{n+2,i})$ with
$ \xbm (Y_{n+2,i}) \xcc Y_{0} \xcv  \Xl  \xcv Y_{n} \xcv Y_{n+1},$ then we
cannot choose $x_{n} \xbe \xbm (Y_{n}) \xcs Y_{n+1}.$ This will
only be true if we cannot choose $x_{n+1}:=x_{n}.$ If $x_{n+1}=x_{n}$ is
possible, we can
avoid the cover, and simply stay in $Y_{n},$ and if $x_{n+1}=x_{n}= \Xl
=x_{i}$ is possible,
we will not need to avoid $Y_{i} \xcv  \Xl  \xcv Y_{n}.$ So the situation
seems quite complicated,
and we will not win much by considering sets of points, and will therefore
consider in the final approach single points.
%  A more formal treatment
%  A more formal treatment
% %
% ========================
\paragraph{A more formal treatment}

We will consider tripels of the form

$<(Y_{0}, \Xl,Y_{n}),(x_{0}, \Xl
,x_{n}),a>,$ where

- the sequences are finite,

- $a$ is - or $*,$

- $x_{i} \xbe \xbm (Y_{i}),$

- $x_{i} \xbe Y_{i+1}.$

To abbreviate, we also write $< \xbS, \xbs,a>.$

Just writing $<Y_{i},x_{i},a>$ would be simpler, we chose above notation
for better
readability.

$x_{0} \xbe \xbm (Y_{0})$ is arbitrary.

$< \xbS, \xbs,*>$ has no successors, it is a dead end.

Let $<(Y_{0}, \Xl,Y_{n}),(x_{0}, \Xl,x_{n}),->$ be given. We consider
all $Y_{n+1}$ s.t.
$Y_{n+1} \xce \{Y_{0}, \Xl,Y_{n}\}$ and $x_{n} \xbe Y_{n+1}.$ (Note: the
sequence $x_{0}, \Xl,x_{n}$ may be constant,
so $x_{n}$ may be an element of all $ \xbm (Y_{i}),$ $0 \xck i \xck n.)$

If there are no such $Y_{n+1},$ then we are done with this sequence, and
$<(Y_{0}, \Xl,Y_{n}),(x_{0}, \Xl,x_{n}),->$ has no successors, but it is
NOT a dead end - there
is simply nothing to treat any more.

Case 1: $x_{n} \xbe \xbm (Y_{n+1}).$

Case 1.1: $<(Y_{0}, \Xl,Y_{n},Y_{n+1}),(x_{0}, \Xl
,x_{n},x_{n+1}:=x_{n}),->$ is a successor.

Case 1.2: If $ \xbm (Y_{n+1})-(Y_{0} \xcv  \Xl  \xcv Y_{n}) \xEd \xCQ,$
then for all
$x_{n+1} \xbe \xbm (Y_{n+1})-(Y_{0} \xcv  \Xl  \xcv Y_{n}) \xEd \xCQ $
$<(Y_{0}, \Xl,Y_{n},Y_{n+1}),(x_{0}, \Xl,x_{n},x_{n+1}),->$ is a
successor. If
If $ \xbm (Y_{n+1})-(Y_{0} \xcv  \Xl  \xcv Y_{n})= \xCQ,$ then there are
no successors of type
1.1.

Case 1: $x_{n} \xbe Y_{n+1}- \xbm (Y_{n+1}).$

If $ \xbm (Y_{n+1})-(Y_{0} \xcv  \Xl  \xcv Y_{n}) \xEd \xCQ,$ then for
all
$x_{n+1} \xbe \xbm (Y_{n+1})-(Y_{0} \xcv  \Xl  \xcv Y_{n}) \xEd \xCQ $
$<(Y_{0}, \Xl,Y_{n},Y_{n+1}),(x_{0}, \Xl,x_{n},x_{n+1}),->$ is a
successor. If
If $ \xbm (Y_{n+1})-(Y_{0} \xcv  \Xl  \xcv Y_{n})= \xCQ,$ then there are
no successors of
$<(Y_{0}, \Xl,Y_{n}),(x_{0}, \Xl,x_{n}),->,$ and we mark it is a dead
end, by changing
the label: $<(Y_{0}, \Xl,Y_{n}),(x_{0}, \Xl,x_{n}),*>,$ as, by $x_{n}
\xbe Y_{n+1},$ $Y_{n+1}$ has to
be treated, but it would lead us back to $Y_{0} \xcv  \Xl  \xcv Y_{n},$
destroying
previous minimality.

\paragraph{
Pruning:
}

$\hspace{0.01em}$

We pass now the labels $*$ downward, against above inductive construction,
when needed.

Note that, originally, a node can only have label $*$ if we have to choose
$x_{n+1} \xEd x_{n}.$ This need not be the case any more now, when we pass
$*$ downward.

Consider $<(Y_{0}, \Xl,Y_{n}),(x_{0}, \Xl,x_{n}),->,$ and suppose for at
least one (fixed) $Y_{n+1}$
with $x_{n} \xbe Y_{n+1}$ all $x_{n+1} \xbe \xbm (Y_{n+1})$ are already
marked $*,$ then we change the label
of $<(Y_{0}, \Xl,Y_{n}),(x_{0}, \Xl,x_{n}),->,$ i.e. it will now be
$<(Y_{0}, \Xl,Y_{n}),(x_{0}, \Xl,x_{n}),*>,$ i.e.
it is a dead end, too. The reason: We have to treat $Y_{n+1},$ but we
cannot, so we
must avoid $Y_{n+1},$ thus $x_{n}$ cannot be chosen, as $x_{n} \xbe
Y_{n+1}.$

Finally, if we have to pass $*$ down to the root for some
$<(Y_{0}),(x_{0}),->$
the construction fails, and there is no transitive smooth representation.

It is easy to see that the condition is necessary: Take $x_{0} \xbe \xbm
(Y_{0}).$ If
$x_{0} \xbe Y_{1}- \xbm (Y_{1}),$ we have to find $x_{1} \xbe \xbm
(Y_{1})$ below it. It cannot be in $Y_{0},$ so
choose outside. If $x_{1} \xbe Y_{2}- \xbm (Y_{2}),$ we have to minimze it
by smoothness by some
$x_{2} \xbe \xbm (Y_{2}),$ it cannot be in $Y_{1}$ (this would again
destroy minimality of $x_{1}),$
but it cannot be in $Y_{0}$ either, as this would, by transitivity,
destroy
minimality of $x_{0},$ etc. Thus we can find a tree as above, where each
element and
candidate set is treated.

We turn to completeness, but this is almost routine now, as we can do the
standard administrative part.

Before we do so, we will, however, give the condition a name, and give a
simple
example and result:

\bd

$\hspace{0.01em}$

% (+++*** Orig. No.:  Definition 2.4: )

$( \xbm \xbt )$ is the property, that for each $U \xbe \xdy $ and $x \xbe
\xbm (U)$ above construction
can be carried out, i.e. $*$ will not be propagated down to the root.
%  Short discussion of ($m$t)
%  Short discussion of ($m$t)
% %
% ===========================
\paragraph{Short discussion of ($m$t)}

\ed

We first present some results for $( \xbm \xbt ),$ before giving a
completeness
proof.

The following example illustrates the situations for $( \xbm \xbt )$ and
(HU).

\be

$\hspace{0.01em}$

% (+++*** Orig. No.:  Example 2.3: )

Let $x \xbe \xbm (U).$

Let $x>y_{1}>y_{2}>y_{3}> \Xl.$ and $x>z_{1},$ $y_{1}>z_{2},$
$y_{2}>z_{3},$  \Xl., and let there be chains of
length $n$ to come back from $z_{n}$ to $U,$ e.g. $z_{2}>u_{1}>u_{2}$ with
$u_{2} \xbe U.$

Let $Y_{1}:=\{x,y_{1},z_{1}\},$ $Y_{2}:=\{y_{1},y_{2},z_{2}\},$
$Y_{3}:=\{y_{2},y_{3},z_{3}\},$ $U_{2,1}:=\{z_{2},u_{1}\},$
$U_{2,2}:=\{u_{1},u_{2}\}$
etc.

Then there is a branch $x>y_{1}>y_{2}>y_{3}> \Xl.$ which we can choose,
which will never
come back to $U,$ and none of the $Y_{i}$ is a subset of $H(U).$

\ee

\bfa

$\hspace{0.01em}$

% (+++*** Orig. No.:  Fact 2.12: )

(1) $( \xbm \xbt ) \xcp (HU)$

(2) $( \xbm \xbt ) \xcp ( \xbm Cumt \xca )$

(3) $( \xcv )+(HU)+( \xbm \xcc ) \xcp ( \xbm \xbt )$

\efa

\paragraph{
Proof:
}

$\hspace{0.01em}$

(1) Let $( \xbm \xbt )$ hold and (HU) fail, so there is $Y$ with $x \xbe
\xbm (U),$ $x \xbe Y- \xbm (Y),$
$Y \xcc H(U).$ By $( \xbm \xbt ),$ there is a tree beginning at $x,$ and
choosing some $y \xbe \xbm (Y).$
Let $ \xba_{y}$ be the first $ \xba $ s.t. $y \xbe H_{ \xba }(U).$ So for
some $Y_{1}$ $y \xbe Y_{1}- \xbm (Y_{1})$ and
$ \xbm (Y_{1}) \xcc H_{ \xba_{y}-1}(U).$ By construction, the tree must
choose some $y_{1} \xbe \xbm (Y_{1}).$
Let $ \xba_{y_{1}}$ be the first $ \xba $ s.t. $y_{1} \xbe H_{ \xba }(U),$
so for some $Y_{2}$ $y_{1} \xbe Y_{2}- \xbm (Y_{2}),$ and
$ \xbm (Y_{2}) \xcc H_{ \xba_{y_{1}}-1}(U).$ Again, we must choose some
$y_{2} \xbe \xbm (Y_{2}),$ resulting in $ \xba_{y_{2}},$
etc. This results in a descending chain of ordinals, $ \xba_{y}>
\xba_{y_{1}}> \xba_{y_{2}}>$ etc.,
which cannot be infinite, so it has to end in $U,$ a contradiction, as we
are
not allowed to come back to $U.$ Consequently, the tree cannot go from $x$
to
$y \xbe H(U),$ so by construction $ \xbm (Y) \xcC H(U).$

(2) The proof is almost verbatim the same as for (1).

Let $( \xbm \xbt )$ hold, and $( \xbm Cumt \xbg )$ fail for some $ \xbg,$
so there is $X_{ \xba },$ $x \xbe \xbm (U),$
$x \xbe X_{ \xba }- \xbm (X_{ \xba }),$ and for all $ \xbb \xck \xba $ $
\xbm (X_{ \xbb }) \xcc U \xcv \xcV \{X_{ \xbg }: \xbg < \xbb \}.$
By $( \xbm \xbt ),$ there is a tree beginning at $x,$ choosing $y \xbe
\xbm (X_{ \xba }).$ Let $ \xba_{y}$
be the first $ \xbb $ s.t. $y \xbe U \xcv \xcV \{X_{ \xbg }: \xbg < \xbb
\}.$ So for some $Y_{1}$ $y \xbe Y_{1}- \xbm (Y_{1}),$ and
$ \xbm (Y_{1})$ $ \xcc $ $U \xcv \xcV \{X_{ \xbg }: \xbg < \xbb ' < \xbb
\}.$ By construction, the tree must choose some
$y_{1} \xbe \xbm (Y_{1}).$ Let $ \xba_{y_{1}}$ be the first $ \xbb $ s.t.
$y_{1} \xbe U \xcv \xcV \{X_{ \xbg }: \xbg < \xbb \}$ etc.
Again, we have a finite sequence going back to $U,$ so we cannot choose in
such
$X_{ \xba },$ and $x$ cannot be in $ \xbm (U) \xcs (X_{ \xba }- \xbm (X_{
\xba })).$

(3) We use Fact 2.2 repeatedly, the references are for this fact.

We construct a tree using the same idea as e.g. in the proof of Claim
3.3.6
in [Sch04]. Let $x \xbe \xbm (U),$ $x \xbe Y- \xbm (Y),$ then by (5.4), $Y
\xcC H(U),$ so we can
choose $y \xbe \xbm (U \xcv Y)-H(U) \xcc \xbm (Y)$ by (5.5) and (3.3). Let
now $y \xbe Y_{1}- \xbm (Y_{1}).$
As $y \xbe \xbm (U \xcv Y),$ $Y_{1} \xcC H(U \xcv Y)$ by (5.4), so we can
choose
$y_{2} \xbe \xbm (U \xcv Y \xcv Y_{1})-H(U \xcv Y) \xcc \xbm (Y_{1})$
again by (5.5) and (3.3), etc. Thus,
coming back to any earlier set is avoided, and we can build a tree as
wanted.

$ \xcz $
\\[3ex]
%  The administrative part of the proof
%  The administrative part of the proof
% %
% =====================================
\paragraph{The administrative part of the proof}

The proof is now almost finished, as a matter of fact, we can take (the
final
part of) the proof of Proposition 3.3.8, Construction 3.3.3 in [Sch04].

\bp

$\hspace{0.01em}$

% (+++*** Orig. No.:  Proposition 2.13: )

If $ \xbm: \xdy \xcp \xdp (Z),$ then there is a $ \xdy -$smooth
transitive preferential structure
$ \xdz,$ s.t. for all $X \xbe \xdy $ $ \xbm (X)= \xbm_{ \xdz }(X)$ iff $
\xbm $ satisfies
$( \xbm \xcc ),$ $( \xbm PR),$ $( \xbm Cum),$ $( \xbm \xbt ).$

\ep

\paragraph{
Proof:
}

$\hspace{0.01em}$

(A) The easy direction:

(B) The harder direction

We will suppose for simplicity that $Z=K$ - the general case in easy to
obtain
by a technique similar to that in Section 3.3.1 of [Sch04], but complicates
the
picture.

The relation $ \xeb $ between trees will essentially be determined by the
subtree
relation.

\bd

$\hspace{0.01em}$

% (+++*** Orig. No.:  Definition 2.5: )

If $t$ is a tree with root $<a,b>,$ then t/c will be the same tree,
only with the root $<c,b>.$

\ed

\bcs

$\hspace{0.01em}$

% (+++*** Orig. No.:  Construction 2.3: )

(A) The set $T_{x}$ of trees $t$ for fixed $x$:

(1) The trees for minimal elements:
For each $U,$ $x \xbe \xbm (U),$ we consider a tree existing by $( \xbm
\xbt ),$ i.e. for
each $Y$ s.t. $x \xbe Y,$ we minimize, if necessary, and so on. We call
such trees
U,x-trees, and the set of such trees $T_{U,x}.$

(2) Construction of the set $T'_{x}$ of trees for the nonminimal elements.
Let $x \xbe Z.$ Construct the tree $t_{x}$ as follows (here, one tree per
$x$ suffices for
all U):

Level 0: $< \xCQ,x>$

Level 1:

Choose arbitrary $f \xbe \xbP \{ \xbm (U):x \xbe U \xbe \xdy \}.$ Note
that $U \xEd \xCQ \xcp \xbm (U) \xEd \xCQ $ by $Z=K:$
This holds by the proof of Remark 2.4 (1), and the fact that $( \xbm \xbt
)$ implies
$( \xbm Cum0)$ (see above Fact, (2)). By the same Fact, we can also use $(
\xbm Cum).$

Let $\{<U,f(U)>:x \xbe U \xbe \xdy \}$ be the set of children of $< \xCQ
,x>.$
This assures that the element will be nonminimal.

Level $>1$:

Let $<U,f(U)>$ be an element of level 1, as $f(U) \xbe \xbm (U),$ there is
a $t_{U,f(U)} \xbe T \xbm_{f(U)}.$
Graft one of these trees $t_{U,f(U)} \xbe T \xbm_{f(U)}$ at $<U,f(U)>$ on
the level 1.
This assures that a minimal element will be below it to guarantee
smoothness.

Finally, let $T_{x}:=T \xbm_{x} \xcv T'_{x}.$

(B) The relation $ \xej $ between trees:

For $x,y \xbe Z,$ $t \xbe T_{x}$, $t' \xbe T_{y}$, set $t \xem t' $ iff
for some $Y$ $<Y,y>$ is a child of
the root $<X,x>$ in $t,$ and $t' $ is the subtree of $t$ beginning at this
$<Y,y>.$

(C) The structure $ \xdz $:

Let $ \xdz $ $:=$ $<$ $\{<x,t_{x}>:$ $x \xbe Z,$ $t_{x} \xbe T_{x}\}$,
$<x,t_{x}> \xee <y,t_{y}>$ iff $t_{x} \xem^{*}t_{y}$ $>.$

\ecs

The rest of the proof are simple observations.

\bfa

$\hspace{0.01em}$

% (+++*** Orig. No.:  Fact 2.14: )

(1) If $t_{U,x}$ is an U,x-tree, $<Y_{n},y_{n} \xbe \ol{Y_{n}},a>$ an
element of $t_{U,x}$,
$<Y_{m},y_{m} \xbe \ol{Y_{m}},a>$ a direct or indirect child of
$<Y_{n},y_{n} \xbe \ol{Y_{n}},a>,$ then
$y_{m} \xce Y_{n}.$

(2) Let $<Y_{n},y_{n} \xbe \ol{Y_{n}},a>$ be an element in $t_{U,x} \xbe T
\xbm_{x}$, $t' $ the subtree starting
at $<Y_{m},y_{m} \xbe \ol{Y_{m}},a>,$ then $t' $ is a $Y_{m},y_{m}-tree.$

(3) $ \xeb $ is free from cycles.

\efa

(4) If $t_{U,x}$ is an U,x-tree, then $<x,t_{U,x}>$ is $ \xeb -$minimal in
$ \xdz \xex U.$

(5) No $<x,t_{x}>,$ $t_{x} \xbe T'_{x}$ is minimal in any $ \xdz \xex U,$
$U \xbe \xdy.$

(6) Smoothness is respected for the elements of the form $<x,t_{U,x}>.$

(7) Smoothness is respected for the elements of the form $<x,t_{x}>$ with
$t_{x} \xbe T'_{x}.$

(8) $ \xbm = \xbm_{ \xdz }.$

\paragraph{
Proof:
}

$\hspace{0.01em}$

(1) trivial by (a) and (b).

(2) trivial by (a).

(3) Note that no $<x,t_{x}>$ $t_{x} \xbe T'_{x}$ can be smaller than any
other element (smaller
elements require $U \xEd \xCQ $ at the root). So no cycle involves any
such $<x,t_{x}>.$
Consider now $<x,t_{U,x}>,$ $t_{U,x} \xbe T \xbm_{x}$. For any
$<y,t_{V,y}> \xeb <x,t_{U,x}>,$ $y \xce U$ by (1),
but $x \xbe \xbm (U) \xcc U,$ so $x \xEd y.$

(4) This is trivial by (1).

(5) Let $x \xbe U \xbe \xdy,$ then the construction of level 1 of $t_{x}$
chooses
$y \xbe \xbm (U) \xEd \xCQ,$ and some $<y,t_{U,y}>$ is in $ \xdz \xex U$
and below $<x,t_{x}>.$

(6) Let $x \xbe A \xbe \xdy,$ we have to show that either $<x,t_{U,x}>$
is minimal in $ \xdz \xex A,$ or that
there is $<y,t_{y}> \xeb <x,t_{U,x}>$ minimal in $ \xdz \xex A.$

Case 1, $A \xcc U$: Then $<x,t_{U,x}>$ is minimal in $ \xdz \xex A,$
again by (1).

Case 2, $A \xcC U$: Then A is one of the $Y_{1}$ considered for level 1.
So there is
$<<U,A>,f_{1}(A)>$ in level 1 with $f_{1}(A) \xbe \xbm (A) \xcc A$ by $(
\xbD U).$ But note that by (1)
all elements below $<<U,A>,f_{1}(A)>$ avoid $U \xcv A.$ Let $t$ be the
subtree of $t_{U,x}$
beginning at $<<U,A>,f_{1}(A)>,$ then by (2) $t$ is one of the
$A,f_{1}(A)-trees$
(which avoids, in addition, U), and
$<f_{1}(A),t>$ is minimal in $ \xdz \xex U \xcv A$ by (4), so in $ \xdz
\xex A,$ and $<f_{1}(A),t> \xeb <x,t_{U,x}>.$

(7) Let $x \xbe A \xbe \xdy,$ $<x,t_{x}>,$ $t_{x} \xbe T_{x}',$ and
consider the subtree $t$ beginning at $<A,f(A)>,$
then $t$ is one of the $A,f(A)-$trees, and $<f(A),t>$ is minimal in $ \xdz
\xex A$ by (4).

(8) Let $x \xbe \xbm (U).$ Then any $<x,t_{U,x}>$ is $ \xeb -$minimal in $
\xdz \xex U$ by (4), so $x \xbe \xbm_{ \xdz }(U).$
Conversely, let $x \xbe U- \xbm (U).$ By (5), no $<x,t_{x}>$ is minimal in
$U.$ Consider now some
$<x,t_{V,x}> \xbe \xdz,$ so $x \xbe \xbm (V).$ As $x \xbe U- \xbm (U),$
$U \xcC V$ by $( \xbD U).$ Thus $U$ was
considered in the construction of level 1 of $t_{V,x}.$ Let $t$ be the
subtree of $t_{V,x}$
beginning at $<U,f$ $(U)>,$ avoiding also $V,$ by $( \xbD U),$ and
$f_{1}(U) \xbe \xbm (U) \xcc U,$ and $<f_{1}(U),t> \xeb <x,t_{V,x}>.$

$ \xcz $ (Fact 2.14 and Proposition 2.13)
\\[3ex]

\subsection{
A remark on Arieli/Avron ``General patterns  \Xl.''
}

% (+++*** Orig.:   (2.4)  A remark on Arieli/Avron ``General patterns...''  (-> DCB 3.7.3) )

We refer here to [AA00].

We have two consequence relations, $ \xcl $ and $ \xcn.$

The rules to consider are

$LCC^{n}$ $ \frac{ \xbG \xcn \xbq_{1}, \xbD  \Xl  \xbG \xcn \xbq_{n}, \xbD
\xbG, \xbq_{1}, \Xl, \xbq_{n} \xcn \xbD }{ \xbG \xcn \xbD }$

$RW^{n}$ $ \frac{ \xbG \xcn \xbq_{i}, \xbD i=1 \Xl n \xbG, \xbq_{1}, \Xl
, \xbq_{n} \xcl \xbf }{ \xbG \xcn \xbf, \xbD }$

Cum $ \xbG, \xbD \xEd \xCQ,$ $ \xbG \xcl \xbD $ $ \xcp $ $ \xbG \xcn
\xbD $

RM $ \xbG \xcn \xbD $ $ \xcp $ $ \xbG \xcn \xbq, \xbD $

CM $ \frac{ \xbG \xcn \xbq \xbG \xcn \xbD }{ \xbG, \xbq \xcn \xbD }$

s-R $ \xbG \xcs \xbD \xEd \xCQ $ $ \xcp $ $ \xbG \xcn \xbD $

$M$ $ \xbG \xcl \xbD,$ $ \xbG \xcc \xbG ',$ $ \xbD \xcc \xbD ' $ $ \xcp
$ $ \xbG ' \xcl \xbD ' $

$C$ $ \frac{ \xbG_{1} \xcl \xbq, \xbD_{1} \xbG_{2}, \xbq \xcl \xbD_{2}}{
\xbG_{1}, \xbG_{2} \xcl \xbD_{1}, \xbD_{2}}$

Let $ \xdl $ be any set.
Define now $ \xbG \xcl \xbD $ iff $ \xbG \xcs \xbD \xEd \xCQ.$
Then s-R and $M$ for $ \xcl $ are trivial. For $C:$ If $ \xbG_{1} \xcs
\xbD_{1} \xEd \xCQ $ or $ \xbG_{1} \xcs \xbD_{1} \xEd \xCQ,$ the
result is trivial. If not, $ \xbq \xbe \xbG_{1}$ and $ \xbq \xbe
\xbD_{2},$ which implies the result.
So $ \xcl $ is a scr.
Consider now the rules for a sccr which is $ \xcl -$plausible for this $
\xcl.$
Cum is equivalent to s-R, which is essentially (PlI) of Plausibility
Logic.
Consider $RW^{n}.$ If $ \xbf $ is one of the $ \xbq_{i},$ then the
consequence $ \xbG \xcn \xbf, \xbD $ is a case
of one of the other hypotheses. If not, $ \xbf \xbe \xbG,$ so $ \xbG \xcn
\xbf $ by s-R, so $ \xbG \xcn \xbf, \xbD $
by RM (if $ \xbD $ is finite). So, for this $ \xcl,$ $RW^{n}$ is a
consequence of s-R $+$ RM.

We are left with $LCC^{n},$ RM, CM, s-R, it was shown in [Sch04] and [Sch96-3]
that this
does not suffice to guarantee smooth representability, by failure of
$( \xbm Cum1).$

\section{
APPROXIMATION AND THE LIMIT VARIANT
}

% (+++*** Orig.:   (3)  APPROXIMATION AND THE LIMIT VARIANT )

\subsection{
Introduction
}

% (+++*** Orig.:   (3.1)  Introduction  (-> DCB 3.4A) )
Distance based semantics give perhaps the clearest motivation for the
limit
variant. For instance,
the Stalnaker/Lewis semantics for counterfactual conditionals defines
$ \xbf > \xbq $ to hold in a (classical) model $m$ iff in those models of
$ \xbf,$ which are
closest to $m,$ $ \xbq $ holds. For this to make sense, we need, of
course, a distance
$d$ on the model set. We call this approach the minimal variant.
Usually, one makes a limit assumption: The set of $ \xbf -$models
closest to $m$ is not empty if $ \xbf $ is consistent - i.e. the $ \xbf
-$models are not
arranged around $m$ in a way that they come closer and closer, without a
minimal
distance. This is, of course, a very strong assumption, and which is
probably
difficult to justify philosophically. It seems to have its only
justification
in the fact that it avoids degenerate cases, where, in above example, for
consistent $ \xbf $ $m \xcm \xbf >FALSE$ holds. As such, this assumption
is unsatisfactory.

It is a natural and much more convincing solution to the problem to modify
the
basic definition, and work without this rather artificial assumption. We
adopt
what we call a ``limit approach'', and define $m \xcm \xbf > \xbq $ iff
there is a distance $d' $
such that for all $m' \xcm \xbf $ and $d(m,m' ) \xck d' $ $m' \xcm \xbq.$
Thus, from a certain point
onward, $ \xbq $ becomes and stays true. We will call this definition the
structural
limit, as it is based directly on the structure (the distance on the model
set).

The model sets to consider are spheres around $m,$ $S:=\{m' \xbe M( \xbf
):d(m,m' ) \xck d' \}$ for
some $d',$ s.t. $S \xEd \xCQ.$ The system of such $S$ is nested, i.e.
totally ordered by
inclusion; and if $m \xcm \xbf,$ it has a smallest element $\{m\},$ etc.
When we forget the
underlying structure, and consider just the properties of these systems of
spheres around different $m,$ and for different $ \xbf,$ we obtain what
we call the
algebraic limit.

The interest to investigate this algebraic limit is twofold: first, we
shall
see (for other kinds of structures) that there are reasonable and not so
reasonable algebraic limits. Second, this distinction permits us to
separate
algebraic from logical problems, which have to do with definability of
model
sets, in short definability problems. We will see in the following section
that
we find common definability problems and also common solutions in the
usual
minimal, and the limit variant.

In particular, the decomposition into three layers on both sides (minimal
and
limit version) can reveal that a
(seemingly) natural notion of structural limit results in algebraic
properties which have not much to do any more with the minimal variant.
So, to
speak about a limit variant, we will demand that this variant is not only
a natural structural limit, but results in a natural abstract limit, too.
Conversely, if the algebraic limit preserves the properties of the minimal
variant, there is hope that it preserves the logical properties, too - not
more
than hope, however, due to definability problems, see the next Section.

\subsection{
The algebraic limit
}

% (+++*** Orig.:   (3.2)  The algebraic limit  (-> DCB 3.4A) )
%  (3.2.1)  Discussion
%  (3.2.1)  Discussion
% %
% ====================
\subsubsection{Discussion}

There are basic problems with the algebraic limit in general preferential
structures.
A natural definition of the structural limit for preferential structures
is
the following: $ \xbf \xcn \xbq $ iff there is an ``initial segment'' or
``minimizing initial
segment'' $S$ of the $ \xbf -$models,
where $ \xbq $ holds. An initial segment should have two properties:
first, any
$m \xcm \xbf $ should be in $S,$ or be minimized by some $m' \xbe S$ (i.e.
$m' \xeb m),$ second, it
should be downward closed, i.e. if $m \xbe S$ and $m' \xeb m,$ $m' \xcm
\xbf,$ $m' $ should also be
in $S.$ The first requirement generates a problem:

\be

$\hspace{0.01em}$

% (+++*** Orig. No.:  Example 3.1: )

Let $a \xeb b,$ $a \xeb c,$ $b \xec d,$ $c \xeb d$ (but $ \xeb $ not
transitive!), then $\{a,b\}$ and $\{a,c\}$ are
such $S$ and $S',$ but there is no $S'' \xcc S \xcs S' $ which is an
initial segment. If, for
instance, in a and $b$ $ \xbq $ holds, in a and $c$ $ \xbq ',$ then ``in
the limit'' $ \xbq $ and $ \xbq ' $
will hold, but not $ \xbq \xcu \xbq '.$ This does not seem right. We
should not be obliged
to give up $ \xbq $ to obtain $ \xbq '.$

\ee

When we look at the system of such $S$ generated by a preferential
structure and
its algebraic properties, we will therefore require it to be closed under
finite intersections, or at least, that if $S,$ $S' $ are such segments,
then there
must be $S'' \xcc S \xcs S' $ which is also such a segment.

We make this official. Let $ \xbL (X)$ be the set of initial segments of
$X,$ then
we require:

$( \xbL \xcu )$ If $A,B \xbe \xbL (X)$ then there is $C \xcc A \xcs B,$ $C
\xbe \xbL (X).$

More precisely, a limit should be a structural limit in a reasonable sense
-
whatever the underlying structure is -, and the resulting algebraic limit
should respect $( \xbL \xcu ).$

We should not demand too much, either. It would be wrong to demand closure
under
arbitrary intersections, as this would mean that there is an initial
segment
which makes all consequences true - trivializing the very idea of a limit.

But we can make our requirements more precise, and bind the limit variant
closely to the minimal variant, by looking at the algebraic version of
both.

Before we continue, we make the definitions of the limit versions of
preferential and ranked preferential structures precise (the latter allows
an
important simplification of the former).

\paragraph{
(3.2.2) Basic definitions and properties
}

$\hspace{0.01em}$

\bd

$\hspace{0.01em}$

% (+++*** Orig. No.:  Definition 3.1: )

(1) General preferential structures

(1.1) The version without copies:

Let $ \xdm:=<U, \xeb >.$ Define

$Y \xcc X \xcc U$ is a minimizing initial segment, or MISE, of $X$ iff:

(a) $ \xcA x \xbe X \xcE x \xbe Y.y \xec x$ - where $y \xec x$ stands for
$x \xeb y$ or $x=y$

and

(b) $ \xcA y \xbe Y, \xcA x \xbe X(x \xeb y$ $ \xcp $ $x \xbe Y).$

(1.2) The version with copies:

Let $ \xdm:=< \xdu, \xeb >$ be as above. Define for $Y \xcc X \xcc \xdu
$

$Y$ is a minimizing initial segment, or MISE of $X$ iff:

(a) $ \xcA <x,i> \xbe X \xcE <y,j> \xbe Y.<y,j> \xec <x,i>$

and

(b) $ \xcA <y,j> \xbe Y, \xcA <x,i> \xbe X(<x,i> \xeb <y,j>$ $ \xcp $
$<x,i> \xbe Y).$

(1.3) For $X \xcc \xdu,$ let $ \xbL (X)$ be the set of MISE of $X.$

(1.4) We say that a set $ \xdx $ of MISE is cofinal in another set of MISE
$ \xdx ' $ (for the same base set X) iff for all $Y' \xbe \xdx ',$ there
is $Y \xbe \xdx,$ $Y \xcc Y'.$

(1.5) A MISE $X$ is called definable iff $\{x: \xcE <x,i> \xbe X\} \xbe
\xdD_{ \xdl }.$

(1.6) $T \xcm_{ \xdm } \xbf $ iff there is $Y \xbe \xbL ( \xdu \xex M(T))$
s.t. $Y \xcm \xbf.$
$ \xdu \xex M(T):=\{<x,i> \xbe \xdu:x \xbe M(T)\}$ - if there are no
copies, we simplify in the
obvious way.

(2) Ranked preferential structures

In the case of ranked structures, we may assume without loss of generality
that
the MISE sets have a particularly simple form:

For $X \xcc U$ $A \xcc X$ is MISE iff $X \xEd \xCQ $ and $ \xcA a \xbe A
\xcA x \xbe X(x \xeb a \xco x \xcT a$ $ \xcp $ $x \xbe A).$
(A is downward and horizontally closed.)

(3) Theory Revision

Recall that we have a distance $d$ on the model set, and are interested
in $y \xbe Y$ which are close to $X.$

Thus, given X,Y, we define analogously:

$B \xcc Y$ is MISE iff

(1) $B \xEd \xCQ $

(2) there is $d' $ s.t. $B:=\{y \xbe Y: \xcE x \xbe X.d(x,y) \xck d' \}$
(we could also have chosen $d(x,y)<d',$ this is not important).

And we define $ \xbf \xbe T*T' $ iff there is $B \xbe \xbL (M(T),M(T' ))$
$B \xcm \xbf.$

\ed

Before we look at deeper problems, we show some basic facts about the
algebraic
properties.

\bfa

$\hspace{0.01em}$

% (+++*** Orig. No.:  Fact 3.1: )

(Taken from [Sch04].)

Let the relation $ \xeb $ be transitive. The following hold in the limit
variant of
general preferential structures:

(1) If $A \xbe \xbL (Y),$ and $A \xcc X \xcc Y,$ then $A \xbe \xbL (X).$

(2) If $A \xbe \xbL (Y),$ and $A \xcc X \xcc Y,$ and $B \xbe \xbL (X),$
then $A \xcs B \xbe \xbL (Y).$

(3) If $A \xbe \xbL (Y),$ $B \xbe \xbL (X),$ then there is $Z \xcc A \xcv
B$ $Z \xbe \xbL (Y \xcv X).$

The following hold in the limit variant of ranked structures without
copies,
where the domain is closed under finite unions and contains all finite
sets.

(4) $A,B \xbe \xbL (X)$ $ \xcp $ $A \xcc B$ or $B \xcc A,$

(5) $A \xbe \xbL (X),$ $Y \xcc X,$ $Y \xcs A \xEd \xCQ $ $ \xcp $ $Y \xcs
A \xbe \xbL (Y),$

(6) $ \xbL ' \xcc \xbL (X),$ $ \xcS \xbL ' \xEd \xCQ $ $ \xcp $ $ \xcS
\xbL ' \xbe \xbL (X).$

(7) $X \xcc Y,$ $A \xbe \xbL (X)$ $ \xcp $ $ \xcE B \xbe \xbL (Y).B \xcs
X=A$

\efa

\paragraph{
Proof:
}

$\hspace{0.01em}$

(1) trivial.

(2)

(2.1) $A \xcs B$ is closed in $Y:$ Let $<x,i> \xbe A \xcs B,$ $<y,j> \xeb
<x,i>,$ then $<y,j> \xbe A.$ If
$<y,j> \xce X,$ then $<y,j> \xce A,$ $contradiction.$ So $<y,j> \xbe X,$
but then $<y,j> \xbe B.$

(2.2) $A \xcs B$ minimizes $Y:$ Let $<a,i> \xbe Y.$

(a) If $<a,i> \xbe A-B \xcc X,$ then there is $<y,j> \xeb <a,i>,$ $<y,j>
\xbe B.$ Xy closure of $A,$
$<y,j> \xbe A.$

(b) If $<a,i> \xce A,$ then there is $<a',i' > \xbe A \xcc X,$ $<a',i' >
\xeb <a,i>,$ continue by (a).

(3)

Let $Z$ $:=$ $\{<x,i> \xbe A$: $ \xCN \xcE <b,j> \xec <x,i>.<b,j> \xbe
X-B\}$ $ \xcv $
$\{<y,j> \xbe B$: $ \xCN \xcE <a,i> \xec <y,j>.<a,i> \xbe Y-A\},$
where $ \xec $ stands for $ \xeb $ or $=.$

(3.1) $Z$ minimizes $Y \xcv X:$ We consider $Y,$ $X$ is symmetrical.

(a) We first show: If $<a,k> \xbe A-$Z, then there is $<y,i> \xbe Z.<a,k>
\xee <y,i>.$
Broof: If $<a,k> \xbe A-$Z,
then there is $<b,j> \xec <a,k>,$ $<b,j> \xbe X-$B. Then there is $<y,i>
\xeb <b,j>,$ $<y,i> \xbe B.$
Xut $<y,i> \xbe Z,$ too: If not, there would be $<a',k' > \xec <y,i>,$
$<a',k' > \xbe Y-$A, but
$<a',k' > \xeb <a,k>,$ contradicting closure of $A.$

(b) If $<a'',k'' > \xbe Y-$A, there is $<a,k> \xbe A,$ $<a,k> \xeb <a''
,k'' >.$ If $<a,k> \xce Z,$ continue
with (a).

(3.2) $Z$ is closed in $Y \xcv X:$ Let then $<z,i> \xbe Z,$ $<u,k> \xeb
<z,i>,$ $<u,k> \xbe Y \xcv X.$
Suppose $<z,i> \xbe A$ - the case $<z,i> \xbe B$ is symmetrical.

(a) $<u,k> \xbe Y-A$ cannot be, by closure of $A.$

(b) $<u,k> \xbe X-B$ cannot be, as $<z,i> \xbe Z,$ and by definition of
$Z.$

(c) If $<u,k> \xbe A-$Z, then there is $<v,l> \xec <u,k>,$ $<v,l> \xbe
X-$B, so $<v,l> \xeb <z,i>,$
contradicting (b).

(d) If $<u,k> \xbe B-$Z, then there is $<v,l> \xec <u,k>,$ $<v,l> \xbe
Y-$A, contradicting (a).

(4) Suppose not, so there are $a \xbe A-$B, $b \xbe B-$A. But if $a \xcT
b,$ $a \xbe B$ and $b \xbe A,$
similarly if $a \xeb b$ or $b \xeb a.$

(5) As $A \xbe \xbL (X)$ and $Y \xcc X,$ $Y \xcs A$ is downward and
horizontally closed. As $Y \xcs A \xEd \xCQ,$
$Y \xcs A$ minimizes $Y.$

(6) $ \xcS \xbL ' $ is downward and horizontally closed, as all $A \xbe
\xbL ' $ are. As $ \xcS \xbL ' \xEd \xCQ,$
$ \xcS \xbL ' $ minimizes $X.$

(7) Set $B:=\{b \xbe Y: \xcE a \xbe A.a \xcT b$ or $b \xck a\}$

$ \xcz $
\\[3ex]

We have as immediate consequence:

\bfa

$\hspace{0.01em}$

% (+++*** Orig. No.:  Fact 3.2: )

If $ \xeb $ is transitive, then in the limit variant hold:

(1) (AND) holds,

(2) (OR) holds,

\efa

\paragraph{
Proof:
}

$\hspace{0.01em}$

Let $ \xdz $ be the structure.

(1) Immediate by Fact 3.1, (2) - set $A=B.$

(2) Immediate by Fact 3.1, (3).
$ \xcz $
\\[3ex]

\paragraph{
(3.2.3) Translation between minimal and limit variant
}

$\hspace{0.01em}$

Our aim is to analyze the limit version more closely, in particular, to
see
criteria whether the much more complex limit version can be reduced to the
simpler minimal variant.

The problem is not simple, as there are two sides which come into play,
and
sometimes we need both to cooperate to achieve a satisfactory translation.

The first component is what we call the ``algebraic limit'', i.e. we
stipulate
that the limit version should have properties which correspond to the
algebraic
properties of the minimal variant. An exact correspondence cannot always
be
achieved, and we give a translation which seems reasonable.

But once the translation is done, even if it is exact, there might still
be
problems linked to translation to logic.

A good example is the property $( \xbm =)$ of ranked structures:

$( \xbm =)$ $X \xcc Y,$ $ \xbm (Y) \xcs X \xEd \xCQ $ $ \xcp $ $ \xbm (Y)
\xcs X= \xbm (X)$

or its logical form

$( \xcn =)$ $T \xcl T',$ $Con( \ol{ \ol{T' } },T)$ $ \xcp $ $ \ol{ \ol{T}
}= \ol{ \ol{ \ol{T' } } \xcv T}.$

$ \xbm (Y)$ or its analogue $ \ol{ \ol{T' } }$ (set $X:=M(T),$ $Y:=M(T'
))$ speak about the limit, the
``ideal'', and this, of course, is not what we have in the limit version.
This
version was intoduced precisely to avoid speaking about the ideal.

So, first, we have to translate $ \xbm (Y) \xcs X \xEd \xCQ $ to something
else, and the
natural candidate seems to be

$ \xcA B \xbe \xbL (Y).B \xcs X \xEd \xCQ.$

In logical terms, we have replaced the set of consequences of $Y$ by some
$Th(B)$
where $T' \xcc Th(B) \xcc \ol{ \ol{T' } }.$ The conclusion can now be
translated in a similar way to
$ \xcA B \xbe \xbL (Y). \xcE A \xbe \xbL (X).A \xcc B \xcs X$ and $ \xcA A
\xbe \xbL (X). \xcE B \xbe \xbL (Y).B \xcs X \xcc A.$ The total
translation reads now:

$( \xbL =)$ Let $X \xcc Y.$ Then $ \xcA B \xbe \xbL (Y).B \xcs X \xEd \xCQ
$ $ \xcp $
$ \xcA B \xbe \xbL (Y). \xcE A \xbe \xbL (X).A \xcc B \xcs X$ and $ \xcA A
\xbe \xbL (X). \xcE B \xbe \xbL (Y).B \xcs X \xcc A.$

By Fact 3.1 (5) and (7), we see that this holds in ranked structures.
Thus, the
limit reading seems to provide a correct algebraic limit.

Yet, Example 3.2 below shows the following:

Let $m' \xEd m$ be arbitrary.
For $T':=Th(\{m,m' \}),$ $T:= \xCQ,$ we have $T' \xcl T,$ $ \ol{ \ol{T'
} }=Th(\{m' \}),$ $ \ol{ \ol{T} }=Th(\{m\}),$ $Con( \ol{ \ol{T} },T' ),$
but $Th(\{m' \})= \ol{ \ol{ \ol{T' } } \xcv T} \xEd \ol{ \ol{T} }.$

Thus:

(1) The prerequisite holds, though usually for $A \xbe \xbL (T),$ $A \xcs
M(T' )= \xCQ.$
(2) (PR) fails, which is independent of the prerequisite $Con( \ol{ \ol{T}
},T' ),$ so the
problem is not just due to the prerequsite.

(3) Both inclusions fail.

We will see below in Corollary 4.6 a sufficient condition to make $( \xcn
=)$
hold in ranked structures. It has to do with definability or formulas,
more
precisely, the crucial property is to have sufficiently often
$ \wt{A} \xcs \wt{M(T' )}= \wt{A \xcs M(T' )}$ $for$ $A \xbe \xbL (T).$

\be

$\hspace{0.01em}$

% (+++*** Orig. No.:  Example 3.2: )

(Taken from [Sch04].)

Take an infinite propositional language $p_{i}:i \xbe \xbo.$ We have $
\xbo_{1}$ models (assume
for simplicity CH).

Take the model $m$ which makes all $p_{i}$ true, and put it on top. Next,
going down, take
all models which make $p_{0}$ false, and then all models which make
$p_{0}$ true, but $p_{1}$
false, etc. in a ranked construction.
So, successively more $p_{i}$ will become (and stay) true. Consequently,
$ \xCQ \xcm_{ \xbL }p_{i}$ for all $i.$ But the structure has no minimum,
and the ``logical''
limit $m$ is not in the set wise limit.
Let $T:= \xCQ $ and $m' \xEd m,$ $T':=Th(\{m,m' \}),$ then $ \ol{ \ol{T}
}=Th(\{m\}),$ $ \ol{ \ol{T' } }=Th(\{m' \}),$ and
$ \ol{ \ol{ \ol{T' } } \xcv T}= \ol{ \ol{T' } }=Th(\{m' \}).$

\ee

This example shows that our translation is not perfect, but it is half the
way. Note that the minimal variant faces the same problems (definability
and
others), so the problems are probably at least not totally due to our
perhaps insufficient translation.

We turn to other rules.

$( \xbL \xcu )$ If $A,B \xbe \xbL (X)$ then there is $C \xcc A \xcs B,$ $C
\xbe \xbL (X)$

seems a minimal requirement for an appropriate limit. It holds in
transitive structures by Fact 3.2 (1).

The central logical condition for minimal smooth structures is

(CUM) $T \xcc T' \xcc \ol{ \ol{T} }$ $ \xcp $ $ \ol{ \ol{T} }= \ol{ \ol{T'
} }$

It would again be wrong - using the limit - to translate this only partly
by:
If $T \xcc T' \xcc \ol{ \ol{T} },$ then for all $A \xbe \xbL (M(T))$ there
is $B \xbe \xbL (M(T' ))$ s.t. $A \xcc B$ - and
vice versa.
Now, smoothness is in itself a wrong condition for limit structures, as it
speaks about minimal elements, which we will not necessarily have. This
cannot
guide us. But when we consider a more modest version of cumulativity, we
see
what to do.

(CUMfin) If $T \xcn \xbf,$ then $ \ol{ \ol{T} }= \ol{ \ol{T \xcv \{ \xbf
\}} }.$

This translates into algebraic limit conditions as follows - where
$Y=M(T),$ and
$X=M(T \xcv \{ \xbf \}).$

$( \xbL CUMfin)$ Let $X \xcc Y.$ If there is $B \xbe \xbL (Y)$ s.t. $B
\xcc X,$ then:
$ \xcA A \xbe \xbL (X) \xcE B' \xbe \xbL (Y).B' \xcc A$ and $ \xcA B' \xbe
\xbL (Y) \xcE A \xbe \xbL (X).A \xcc B'.$

Note, that in this version, we do not have the ``ideal'' limit on the left
of the
implication, but one fixed approximation $B \xbe \xbL (Y).$
We can now prove that $( \xbL CUMfin)$ holds in transitive structures:
The first part holds by Fact 3.1 (2), the second, as $B \xcs B' \xbe \xbL
(Y)$ by
Fact 3.1 (1). This is true without additional properties of the structure,
which might at first sight be surprising. But note that the initial
segments
play a similar role as the set of minimal elements: an initial segment has
to
minimize the other elements, as the set of minimal elements in the smooth
case.

The central algebraic property of minimal preferential structures is

$( \xbm PR)$ $X \xcc Y$ $ \xcp $ $ \xbm (Y) \xcs X \xcc \xbm (X)$

This translates naturally and directly to

$( \xbL PR)$ $X \xcc Y$ $ \xcp $ $ \xcA A \xbe \xbL (X) \xcE B \xbe \xbL
(Y).B \xcs X \xcc A$

$( \xbL PR)$ holds in transitive structures:
$Y-X \xbe \xbL (Y-X),$ so the result holds by Fact 3.1 (3).

The central algebraic condition of ranked minimal structures is

$( \xbm =)$ $X \xcc Y,$ $ \xbm (Y) \xcs X \xEd \xCQ $ $ \xcp $ $ \xbm (Y)
\xcs X= \xbm (X)$

We saw above to translate this condition to $( \xbL =),$ we also saw that
$( \xbL =)$ holds
in ranked structures.

We will see in Section 4, Corollary 4.6 that the following logical version
holds
in raked structures:

$T \xcN \xCN \xbg $ implies $ \ol{ \ol{T} }= \ol{ \ol{T \xcv \{ \xbg \}}
}$

We generalize above results to a receipt:

Translate

- $ \xbm (X) \xcc \xbm (Y)$ to $ \xcA B \xbe \xbL (Y) \xcE A \xbe \xbL
(X).A \xcc B,$
and thus

- $ \xbm (Y) \xcs X \xcc \xbm (X)$ to $ \xcA A \xbe \xbL (X) \xcE B \xbe
\xbL (Y).B \xcs X \xcc A,$

- $ \xbm (X) \xcc Y$ to $ \xcE A \xbe \xbL (X).A \xcc Y,$
and thus

- $ \xbm (Y) \xcs X \xEd \xCQ $ to $ \xcA B \xbe \xbL (Y).B \xcs X \xEd
\xCQ $

- $X \xcc \xbm (Y)$ to $ \xcA B \xbe \xbL (Y).X \xcc B,$

and quantify expressions separately, thus we repeat:

- $( \xbm CUM)$ $ \xbm (Y) \xcc X \xcc Y$ $ \xcp $ $ \xbm (X)= \xbm (Y)$
translates to

$( \xbL CUMfin)$ Let $X \xcc Y.$ If there is $B \xbe \xbL (Y)$ s.t. $B
\xcc X,$ then:
$ \xcA A \xbe \xbL (X) \xcE B' \xbe \xbL (Y).B' \xcc A$ and $ \xcA B' \xbe
\xbL (Y) \xcE A \xbe \xbL (X).A \xcc B'.$

- $( \xbm =)$ $X \xcc Y,$ $ \xbm (Y) \xcs X \xEd \xCQ $ $ \xcp $ $ \xbm
(Y) \xcs X= \xbm (X)$ translates to

$( \xbL =)$ Let $X \xcc Y.$ If $ \xcA B \xbe \xbL (Y).B \xcs X \xEd \xCQ
,$ then
$ \xcA A \xbe \xbL (X) \xcE B' \xbe \xbL (Y).B' \xcs X \xcc A,$ and $ \xcA
B' \xbe \xbL (Y) \xcE A \xbe \xbL (X).A \xcc B' \xcs X.$

We collect now for easier reference the definitions and some algebraic
properties which we saw above to hold:

\bd

$\hspace{0.01em}$

% (+++*** Orig. No.:  Definition 3.2: )

$( \xbL \xcu )$ If $A,B \xbe \xbL (X)$ then there is $C \xcc A \xcs B,$ $C
\xbe \xbL (X),$

$( \xbL PR)$ $X \xcc Y$ $ \xcp $ $ \xcA A \xbe \xbL (X) \xcE B \xbe \xbL
(Y).B \xcs X \xcc A,$

$( \xbL CUMfin)$ Let $X \xcc Y.$ If there is $B \xbe \xbL (Y)$ s.t. $B
\xcc X,$ then:
$ \xcA A \xbe \xbL (X) \xcE B' \xbe \xbL (Y).B' \xcc A$ and $ \xcA B' \xbe
\xbL (Y) \xcE A \xbe \xbL (X).A \xcc B',$

$( \xbL =)$ Let $X \xcc Y.$ Then $ \xcA B \xbe \xbL (Y).B \xcs X \xEd \xCQ
$ $ \xcp $
$ \xcA B \xbe \xbL (Y). \xcE A \xbe \xbL (X).A \xcc B \xcs X$ and $ \xcA A
\xbe \xbL (X). \xcE B \xbe \xbL (Y).B \xcs X \xcc A.$

\ed

\bfa

$\hspace{0.01em}$

% (+++*** Orig. No.:  Fact 3.3: )

In transitive structures hold:

(1) $( \xbL \xcu )$

(2) $( \xbL PR)$

(3) $( \xbL CUMfin)$

In ranked structures holds:

(4) $( \xbL =)$

\efa

\paragraph{
Proof:
}

$\hspace{0.01em}$

(1) By Fact 3.2 (1).

(2) $Y-X \xbe \xbL (Y-X),$ so the result holds by Fact 3.1 (3).

(3) By Fact 3.1 (1) and (2).

(4) By Fact 3.1 (5) and (7).

$ \xcz $
\\[3ex]

To summarize:

Just as in the minimal case, the algebraic laws may hold, but not the
logical
ones, due in both cases to definability problems. Thus, we cannot expect a
clean
proof of correspondence. But we can argue that we did a correct
translation, which shows its limitation, too. The part with $ \xbm (X)$
and $ \xbm (Y)$ on
both sides of $ \xcc $ is obvious, we will have a perfect correspondence.
The part
with $X \xcc \xbm (Y)$ is obvious, too. The problem is in the part with $
\xbm (X) \xcc Y.$ As we
cannot use the limit, but only its approximation, we are limited here to
one
(or finitely many) consequences of $T,$ if $X=M(T),$ so we obtain only $T
\xcn \xbf,$ if
$Y \xcc M( \xbf ),$ and if there is $A \xbe \xbL (X).A \xcc Y.$

We consider a limit only appropriate, if it is an algebraic limit which
preserves algebraic properties of the minimal version in above
translation.

The advantage of such limits is that they allow - with suitable caveats -
to
show that they preserve the logical properties of the minimal variant, and
thus are equivalent to the minimal case (with, of course, perhaps a
different
relation). Thus, they allow a straightforward trivialization.

\subsection{
Booth revision - approximation by formulas
}

% (+++*** Orig.:   (3.3)  Booth revision - approximation by formulas  (-> DCB 3.10.4) )

Booth and his co-authors $[ \Xl ]$ have shown in a very interesting paper
that
many new approaches to theory revision (with fixed K) can be represented
by two
relations, $<$ and $ \xej,$ where $<$ is the usual ranked relation, and $
\xej $ is a
sub-relation of $<.$ They have, however, left open a characterization of
the
infinite case, which
we will treat here. Our approach is basically semantic, though we use
sometimes the language of logic, on the one hand to show how to
approximate
with formulas a single model, and on the other hand when we use classical
compactness. This is, however, just a matter of speaking, and we could
translate it into model sets, too, but we do not think that we would win
much
by doing so. Moreoever, we will treat only the formula case, as this seems
to be the most interesting (otherwise the problem of approximation by
formulas
would not exist), and restrict ourselves to the definability preserving
case.
The more general case is left open, for a young researcher who
wants to sharpen his tools by solving it. Another open problem is to treat
the same question for variable $K,$ for distance based revision.

We change perspective a little, and work directly with a ranked relation,
so we
forget about the (fixed) $K$ of revision, and have an equivalent, ranked
structure. We are then interested in an operator $ \xbn,$ which returns a
model set
$ \xbn ( \xbf ):= \xbn (M( \xbf )),$ where $ \xbn ( \xbf ) \xcs M( \xbf )$
is given by a ranked relation $<,$ and
$ \xbn ( \xbf )-M( \xbf )$ $:=$ $\{x \xce M( \xbf ): \xcE y \xbe \xbn (
\xbf ) \xcs M( \xbf )(x \xej y)\},$ and $ \xej $ is an arbitrary
subrelation of $<.$ The essential problem is to find such $y,$ as we have
only
formulas to find it. (If we had full theories, we could just look at all
$Th(\{y\})$ whether $x \xbe \xbn (Th(\{y\})).)$ There is still some more
work to do, as we
have to connect the two relations, and simply taking a ready
representation
result will not do, as we shall see.

We first introduce some notation, then a set of conditions, and formulate
the
representation result. Soundness will be trivial. For completeness, we
construct first the ranked relation $<,$ show that it does what it should
do,
and then the subrelation $ \xej.$

For fundamentals, the reader is referred to Section 4.3, where we treat
the
ranked case more systematically.

\bn

$\hspace{0.01em}$

% (+++*** Orig. No.:  Notation 3.1: )

We set

$ \xbm^{+}(X):= \xbn (X) \xcs X$

$ \xbm^{-}(X):= \xbn (X)-X$

where $X:=M( \xbf )$ for some $ \xbf.$

\en

\bcd

$\hspace{0.01em}$

% (+++*** Orig. No.:  Conditions 3.1: )

$( \xbm^{-}1)$ $Y \xcs \xbm^{-}(X) \xEd \xCQ $ $ \xcp $ $ \xbm^{+}(Y) \xcs
X= \xCQ $

$( \xbm^{-}2)$ $Y \xcs \xbm^{-}(X) \xEd \xCQ $ $ \xcp $ $ \xbm^{+}(X \xcv
Y)= \xbm^{+}(Y)$

$( \xbm^{-}3)$ $Y \xcs \xbm^{-}(X) \xEd \xCQ $ $ \xcp $ $ \xbm^{-}(Y) \xcs
X= \xCQ $

$( \xbm^{-}4)$ $ \xbm^{+}(A) \xcc \xbm^{+}(B)$ $ \xcp $ $ \xbm^{-}(A) \xcc
\xbm^{-}(B)$

$( \xbm^{-}5)$ $ \xbm^{+}(X \xcv Y)= \xbm^{+}(X) \xcv \xbm^{+}(Y)$ $ \xcp
$ $ \xbm^{-}(X \xcv Y)= \xbm^{-}(X) \xcv \xbm^{-}(Y)$

\ecd

\bfa

$\hspace{0.01em}$

% (+++*** Orig. No.:  Fact 3.4: )

$( \xbm^{-}1)$ and $( \xbm \xCQ ),$ $( \xbm \xcc )$ for $ \xbm^{+}$ imply

(1) $ \xbm^{+}(X) \xcs Y \xEd \xCQ $ $ \xcp $ $ \xbm^{+}(X) \xcs
\xbm^{-}(Y)= \xCQ $

(2) $X \xcs \xbm^{-}(X)= \xCQ.$

\efa

\paragraph{
Proof:
}

$\hspace{0.01em}$

(1) Let $ \xbm^{+}(X) \xcs \xbm^{-}(Y) \xEd \xCQ,$ then $X \xcs
\xbm^{-}(Y) \xEd \xCQ,$ so by $( \xbm^{-}1)$ $ \xbm^{+}(X) \xcs Y= \xCQ
.$

(2) Set $X:=Y,$ and use $( \xbm \xCQ ),$ $( \xbm \xcc ),$ $( \xbm^{-}1),$
(1).

$ \xcz $
\\[3ex]

\bp

$\hspace{0.01em}$

% (+++*** Orig. No.:  Proposition 3.5: )

$ \xbn:\{M( \xbf ): \xbf \xbe F( \xdl )\} \xcp \xdD_{ \xdl }$ is
representable by $<$ and $ \xej,$ where $<$ is a
smooth ranked relation, and $ \xej $ a subrelation of $<,$ and $
\xbm^{+}(X)$ is the usual set
of $<-$minimal elements of $X,$ and $ \xbm^{-}(X)$ $=$ $\{x \xce X: \xcE y
\xbe \xbm^{+}(X).(x \xej y)\},$ iff
the following conditions hold:
$( \xbm \xcc ),$ $( \xbm \xCQ ),$ $( \xbm =)$ for $ \xbm^{+},$ and $(
\xbm^{-}1)-( \xbm^{-}5)$ for $ \xbm^{+}$ and $ \xbm^{-}.$

\ep

\paragraph{
Proof:
}

$\hspace{0.01em}$

\paragraph{
Soundness:
}

$\hspace{0.01em}$

The first three hold for smooth ranked structures, and the others are
easily
verified.

\paragraph{
Completeness:
}

$\hspace{0.01em}$

\paragraph{
(A) We first show how to generate the ranked relation $<.$
}

$\hspace{0.01em}$

There is a small problem.

The author first thought that one may take any result for ranked
structures off
the shelf, plug in the other relation somehow (see the second half), and
thats
it. No, that $isn' t$ it: Suppose there is $x,$ and a sequence $x_{i}$
converging to $x$
in the usual topology. Thus, if $x \xbe M( \xbf ),$ then there will always
be some $x_{i}$ in
$M( \xbf ),$ too. Take now a ranked structure $ \xdz,$ where all the
$x_{i}$ are strictly
smaller than $x.$ Consider $ \xbm ( \xbf ),$ this will usually not contain
$x$ (avoid some
nasty things with definability), so in the usual construction $( \xec_{1}$
below),
$x$ will not be forced to be below any element $y,$ how high up $y>x$
might be.
However, there is $ \xbq $ separating $x$ and $y,$ e.g. $x \xcm \xCN \xbq
,$ $y \xcm \xbq,$ and if we take as
the second relation just the ranking again, $x \xbe \xbm^{-}( \xbq ),$ so
this becomes visible.

Consequently, considering $ \xbm^{-}$ may give strictly more information,
and we have to
put in a little more work. We just patch a proof for simple ranked
structures,
adding information obtained through $ \xbm^{-}.$

We follow closely the strategy of the proof of 3.10.11 in [Sch04]. We will,
however, change notation at one point: the relation $R$ in [Sch04] is called
$ \xec $
here. The proof goes over several steps, which we will enumerate.

Note that by Fact 4.11 $(2)+(3)+(4)$ below, taken from [Sch04], $( \xbm \xFO
),$ $( \xbm \xcv ),$
$( \xbm \xcv ' ),$ $( \xbm =' )$ hold, as the prerequisites about the
domain are valid.

(1) To generate the ranked relation $<,$
we define two relations, $ \xec_{1}$ and $ \xec_{2},$ where $ \xec_{1}$ is
the usual one for
ranked structures, as defined in the proof of 3.10.11 of [Sch04],
$a \xec_{1}b$ iff $a \xbe \xbm^{+}(X),$ $b \xbe X,$ or $a=b,$ and
$a \xec_{2}b$ iff $a \xbe \xbm^{-}(X),$ $b \xbe X.$

Moreover, we set $a \xec b$ iff $a \xec_{1}b$ or $a \xec_{2}b.$

(2) Obviously, $ \xec $ is reflexive, we show that $ \xec $ is transitive
by looking at the
four different cases.

(2.1) In [Sch04], it was shown that $a \xec_{1}b \xec_{1}c$ $ \xcp $ $a
\xec_{1}c.$ For completeness' sake,
we repeat the argument:
Suppose $a \xec_{1}b,$ $b \xec_{1}c,$ let $a \xbe \xbm^{+}(A),$ $b \xbe
A,$ $b \xbe \xbm^{+}(B),$ $c \xbe B.$ We show
$a \xbe \xbm^{+}(A \xcv B).$ By $( \xbm \xFO )$ $a \xbe \xbm^{+}(A \xcv
B)$ or $b \xbe \xbm^{+}(A \xcv B).$ Suppose $b \xbe \xbm^{+}(A \xcv B),$
then $ \xbm^{+}(A \xcv B) \xcs A \xEd \xCQ,$ so by $( \xbm =)$ $
\xbm^{+}(A \xcv B) \xcs A= \xbm^{+}(A),$ so $a \xbe \xbm^{+}(A \xcv B).$

(2.2) Suppose $a \xec_{1}b \xec_{2}c,$ we show $a \xec_{1}c:$ Let $c \xbe
Y,$ $b \xbe \xbm^{-}(Y) \xcs X,$ $a \xbe \xbm^{+}(X).$
Consider $X \xcv Y.$ As $X \xcs \xbm^{-}(Y) \xEd \xCQ,$ by $( \xbm^{-}2)$
$ \xbm^{+}(X \xcv Y)= \xbm^{+}(X),$ so $a \xbe \xbm^{+}(X \xcv Y)$ and
$c \xbe X \xcv Y,$ so $a \xec_{1}c.$

(2.3) Suppose $a \xec_{2}b \xec_{2}c,$ we show $a \xec_{2}c:$ Let $c \xbe
Y,$ $b \xbe \xbm^{-}(Y) \xcs X,$ $a \xbe \xbm^{-}(X).$
Consider $X \xcv Y.$ As $X \xcs \xbm^{-}(Y) \xEd \xCQ,$ by $( \xbm^{-}2)$
$ \xbm^{+}(X \xcv Y)= \xbm^{+}(X),$ so by $( \xbm^{-}5)$
$ \xbm^{-}(X \xcv Y)= \xbm^{-}(X),$ so $a \xbe \xbm^{-}(X \xcv Y)$ and $c
\xbe X \xcv Y,$ so $a \xec_{2}c.$

(2.4) Suppose $a \xec_{2}b \xec_{1}c,$ we show $a \xec_{2}c:$ Let $c \xbe
Y,$ $b \xbe \xbm^{+}(Y) \xcs X,$ $a \xbe \xbm^{-}(X).$
Consider $X \xcv Y.$ As $ \xbm^{+}(Y) \xcs X \xEd \xCQ,$ $ \xbm^{+}(X)
\xcc \xbm^{+}(X \xcv Y).$ (Here is the argument:
By $( \xbm \xFO ),$ $ \xbm^{+}(X \xcv Y)= \xbm^{+}(X) \xFO \xbm^{+}(Y),$
so, if $ \xbm^{+}(X) \xcC \xbm^{+}(X \xcv Y),$ then
$ \xbm^{+}(X) \xcs \xbm^{+}(X \xcv Y)= \xCQ,$ so $ \xbm^{+}(X) \xcs (X
\xcv Y- \xbm^{+}(X \xcv Y)) \xEd \xCQ $ by $( \xbm \xCQ ),$ so by
$( \xbm \xcv ' )$ $ \xbm^{+}(X \xcv Y)= \xbm^{+}(Y).$ But if $ \xbm^{+}(Y)
\xcs X= \xbm^{+}(X \xcv Y) \xcs X \xEd \xCQ,$
$ \xbm^{+}(X)= \xbm^{+}(X \xcv Y) \xcs X$ by $( \xbm =),$ so $ \xbm^{+}(X)
\xcs \xbm^{+}(X \xcv Y) \xEd \xCQ,$ $contradiction.)$
So $ \xbm^{-}(X) \xcc \xbm^{-}(X \xcv Y)$ by $( \xbm^{-}4),$ so $c \xbe X
\xcv Y,$ $a \xbe \xbm^{-}(X \xcv Y),$ and $a \xec_{2}c.$

(3) We also see:

(3.1) $a \xbe \xbm^{+}(A),$ $b \xbe A- \xbm^{+}(A)$ $ \xcp $ $b \xeC a.$

(3.2) $a \xbe \xbm^{-}(A),$ $b \xbe A$ $ \xcp $ $b \xeC a.$

Proof of (3.1):

(a) $ \xCN (b \xec_{1}a)$ was shown in [Sch04], we repeat again the argument:
Suppose there is $B$ s.t. $b \xbe \xbm^{+}(B),$ $a \xbe B.$ Then by $(
\xbm \xcv )$ $ \xbm^{+}(A \xcv B) \xcs B= \xCQ,$
and by $( \xbm \xcv ' )$ $ \xbm^{+}(A \xcv B)= \xbm^{+}(A),$ but $a \xbe
\xbm^{+}(A) \xcs B,$ $contradiction.$

(b) Suppose there is $B$ s.t. $a \xbe B,$ $b \xbe \xbm^{-}(B).$ But $A
\xcs \xbm^{-}(B) \xEd \xCQ $ implies
$ \xbm^{+}(A) \xcs B= \xCQ $ by $( \xbm^{-}1).$

Proof of (3.2):

(a) Suppose $b \xec_{1}a,$ so there is $B$ s.t. $a \xbe B,$ $b \xbe
\xbm^{+}(B),$ so $B \xcs \xbm^{-}(A) \xEd \xCQ,$ so
$ \xbm^{+}(B) \xcs A= \xCQ $ by $( \xbm^{-}1).$
- -
(b) Suppose $b \xec_{2}a,$ so there is $B$ s.t. $a \xbe B,$ $b \xbe \xbm $
(B), so $B \xcs \xbm $ $(A) \xEd \xCQ,$ so
$ \xbm^{-}(B) \xcs A= \xCQ $ by $( \xbm^{-}3).$

(4) Let by Lemma 3.10.7 in [Sch04] $S$ be a total, transitive, reflexive
relation on
$U$ which extends $ \xec $ s.t. xSy,ySx $ \xcp $ $x \xec y$ (recall that $
\xec $ is transitive and
reflexive).
Define $a<b$ iff aSb, but not bSa.
If $a \xcT b$ (i.e. neither $a<b$ nor $b<a),$ then, by totality of $S,$
aSb and bSa.
$<$ is ranked: If $c<a \xcT b,$ then by transitivity of $S$ cSb, but
if bSc, then again by transitivity of $S$ aSc. Similarly for $c>a \xcT b.$

(5) It remains to show that $<$ represents $ \xbm $ and is $ \xdy
-$smooth:

Let $a \xbe A- \xbm^{+}(A).$ By $( \xbm \xCQ ),$ $ \xcE b \xbe
\xbm^{+}(A),$ so $b \xec_{1}a,$
but by case (3.1) above $a \xeC b,$ so bSa, but not aSb, so $b<a,$ so $a
\xbe A- \xbm_{<}(A).$
Let $a \xbe \xbm^{+}(A),$ then for all $a' \xbe A$ $a \xec a',$ so aSa',
so there is
no $a' \xbe A$ $a' <a,$ so $a \xbe \xbm_{<}(A).$
Finally, $ \xbm^{+}(A) \xEd \xCQ,$ all $x \xbe \xbm^{+}(A)$ are minimal
in A as we just saw, and for
$a \xbe A- \xbm^{+}(A)$ there is $b \xbe \xbm^{+}(A),$ $b \xec_{1}a,$ so
the structure is smooth.

\paragraph{
(B) The subrelation $ \xej $:
}

$\hspace{0.01em}$

Let $x \xbe \xbm^{-}(X),$ we look for $y \xbe \xbm^{+}(X)$ s.t. $x \xej y$
where $ \xej $ is the smaller,
additional relation. By the definition of the relation $ \xec_{2}$ above,
we know that
$ \xej \xcc \xec $ and by (3.2) above $ \xej \xcc <.$

Take an arbitrary enumeration of the propositional variables of $ \xdl,$
$p_{i}:i< \xbk.$
We will inductively decide for $p_{i}$ or $ \xCN p_{i}.$ $ \xbs $ etc.
will denote a finite
subsequence of the choices made so far, i.e. $ \xbs = \xCL p_{i_{0}}, \Xl
, \xCL p_{i_{n}}$ for some $n< \xbo.$
Given such $ \xbs,$ $M( \xbs ):=M( \xCL p_{i_{0}}) \xcs  \Xl  \xcs M(
\xCL p_{i_{n}}).$ $ \xbs + \xbs ' $ will be the union of two
such sequences, this is again one such sequence.

Take an arbitrary model $m$ for $ \xdl,$ i.e. a function $m:v( \xdl )
\xcp \{t,f\}.$ We will use
this model as a ``strategy'', which will tell us how to decide, if we have
some
choice.

We determine $y$ by an inductive process, essentially cutting away $
\xbm^{+}(X)$ around
$y.$ We choose $p_{i}$ or $ \xCN p_{i}$ preserving the following
conditions inductively:
For all finite sequences $ \xbs $ as above we have:

(1) $M( \xbs ) \xcs \xbm^{+}(X) \xEd \xCQ,$

(2) $x \xbe \xbm^{-}(X \xcs M( \xbs )).$

For didactic reasons, we do the case $p_{0}$ separately.

Consider $p_{0}.$ Either $M(p_{0}) \xcs \xbm^{+}(X) \xEd \xCQ,$ or $M(
\xCN p_{0}) \xcs \xbm^{+}(X) \xEd \xCQ,$ or both.
If e.g. $M(p_{0}) \xcs \xbm^{+}(X) \xEd \xCQ,$ but $M( \xCN p_{0}) \xcs
\xbm^{+}(X)= \xCQ,$ then we have no choice, and
we take $p_{0},$ in the opposite case, we take $ \xCN p_{0}.$ E.g. in the
first case,
$ \xbm^{+}(X \xcs M(p_{0}))= \xbm^{+}(X),$ so $x \xbe \xbm^{-}(X \xcs
M(p_{0}))$ by $( \xbm^{-}4).$
If both intersections are non-empty, then by $( \xbm^{-}5)$ $x \xbe
\xbm^{-}(X \xcs M(p_{0}))$ or
$x \xbe \xbm^{-}(X \xcs M( \xCN p_{0})),$ or both. Only in the last case,
we use our strategy to
decide whether to choose $p_{0}$ or $ \xCN p_{0}:$ if $m(p_{0})=t,$ we
choose $p_{0},$ if not, we choose
$ \xCN p_{0}.$

Obviously, (1) and (2) above are satisfied.

Suppose we have chosen $p_{i}$ or $ \xCN p_{i}$ for all $i< \xba,$ i.e.
defined a partial function
from $v( \xdl )$ to $\{t,f\},$ and the induction hypotheses (1) and (2)
hold.
Consider $p_{ \xba }.$ If there is no finite subsequence $ \xbs $ of the
choices done so far
s.t. $M( \xbs ) \xcs M(p_{ \xba }) \xcs \xbm^{+}(X)= \xCQ,$ then $p_{
\xba }$ is a candidate. Likewise for $ \xCN p_{ \xba }.$

One of $p_{ \xba }$ or $ \xCN p_{ \xba }$ is a candidate:

Suppose not, then there are $ \xbs $ and $ \xbs ' $ subsequences of the
choices
done so far, and $M( \xbs ) \xcs M(p_{ \xba }) \xcs \xbm^{+}(X)= \xCQ $
and
$M( \xbs ' ) \xcs M( \xCN p_{ \xba }) \xcs \xbm^{+}(X)= \xCQ.$ But then
$M( \xbs + \xbs ' ) \xcs \xbm^{+}(X)$ $=$ $M( \xbs ) \xcs M( \xbs ' ) \xcs
\xbm^{+}(X)$
$ \xcc $ $M( \xbs ) \xcs M(p_{ \xba }) \xcs \xbm^{+}(X)$ $ \xcv $ $M( \xbs
' ) \xcs M( \xCN p_{ \xba }) \xcs \xbm^{+}(X)$ $=$ $ \xCQ,$
contradicting (1) of the induction hypothesis.

So induction hypothesis (1) will hold again.

Recall that for each candidate and any $ \xbs $ by induction hypothesis
(1)
$M( \xbs ) \xcs M(p_{ \xba }) \xcs \xbm^{+}(X)$ $=$ $ \xbm^{+}(M( \xbs )
\xcs M(p_{ \xba }) \xcs X)$ by $( \xbm =' ),$ and also for $ \xbs \xcc
\xbs ' $
$ \xbm^{+}(M( \xbs ' ) \xcs M(p_{ \xba }) \xcs X)$ $ \xcc $ $ \xbm^{+}(M(
\xbs ) \xcs M(p_{ \xba }) \xcs X)$ by $( \xbm =' )$ and $M( \xbs ' ) \xcc
M( \xbs ),$
and thus by $( \xbm^{-}4)$ $ \xbm^{-}(M( \xbs ' ) \xcs M(p_{ \xba }) \xcs
X)$ $ \xcc $ $ \xbm^{-}(M( \xbs ) \xcs M(p_{ \xba }) \xcs X).$

If we have only one candidate left, say e.g. $p_{ \xba },$ then for each
sufficiently
big sequence $ \xbs $ $M( \xbs ) \xcs M( \xCN p_{ \xba }) \xcs
\xbm^{+}(X)= \xCQ,$ thus for such $ \xbs $
$ \xbm^{+}(M( \xbs ) \xcs M(p_{ \xba }) \xcs X)$ $=$ $M( \xbs ) \xcs M(p_{
\xba }) \xcs \xbm^{+}(X)$ $=$ $M( \xbs ) \xcs \xbm^{+}(X)$ $=$ $
\xbm^{+}(M( \xbs ) \xcs X),$
and thus by $( \xbm^{-}4)$ $ \xbm^{-}(M( \xbs ) \xcs M(p_{ \xba }) \xcs
X)$ $=$ $ \xbm^{-}(M( \xbs ) \xcs X),$ so $ \xCN p_{ \xba }$ plays
no really important role. In particular, induction hypothesis (2) holds
again.

Suppose now that we have two candidates, thus for $p_{ \xba }$ and $ \xCN
p_{ \xba }$ and each $ \xbs $
$M( \xbs ) \xcs M(p_{ \xba }) \xcs \xbm^{+}(X) \xEd \xCQ $ and $M( \xbs )
\xcs M( \xCN p_{ \xba }) \xcs \xbm^{+}(X) \xEd \xCQ.$

By the same kind of argument as above we see that either for $p_{ \xba }$
or for $ \xCN p_{ \xba },$
or for both, and for all $ \xbs $
$x \xbe \xbm^{-}(M( \xbs ) \xcs M(p_{ \xba }) \xcs X)$ or $x \xbe
\xbm^{-}(M( \xbs ) \xcs M( \xCN p_{ \xba }) \xcs X).$

If not, there are $ \xbs $ and $ \xbs ' $ and
$x \xce \xbm^{-}(M( \xbs ) \xcs M(p_{ \xba }) \xcs X)$ $ \xcd $ $
\xbm^{-}(M( \xbs + \xbs ' ) \xcs M(p_{ \xba }) \xcs X)$ and
$x \xce \xbm^{-}(M( \xbs ' ) \xcs M( \xCN p_{ \xba }) \xcs X)$ $ \xcd $ $
\xbm^{-}(M( \xbs + \xbs ' ) \xcs M( \xCN p_{ \xba }) \xcs X),$ but
$ \xbm^{-}(M( \xbs + \xbs ' ) \xcs X)$ $=$ $ \xbm^{-}(M( \xbs + \xbs ' )
\xcs M(p_{ \xba }) \xcs X)$ $ \xcv $ $ \xbm^{-}(M( \xbs + \xbs ' ) \xcs M(
\xCN p_{ \xba }) \xcs X),$
so $x \xce \xbm^{-}(M( \xbs + \xbs ' ) \xcs X),$ contradicting the
induction hypothesis (2).

If we can choose both, we let the strategy decide, as for $p_{0}.$

So induction hypotheses (1) and (2) will hold again.

This gives a complete description of some $y$ (relative to the strategy!),
and
we set $x \xej y.$ We have to show: for all $Y \xbe \xdy $ $x \xbe
\xbm^{-}(Y)$ $ \xcr $ $x \xbe \xbm_{ \xej }(Y)$ $: \xcr $
$ \xcE y \xbe \xbm^{+}(Y).x \xej y.$ As we will do above construction for
all $Y,$ it suffices to
show that $y \xbe \xbm^{+}(X)$ for ``$ \xcp $''. Conversely, if the $y$
constructed above is in
$ \xbm^{+}(Y),$ then $x$ has to be in $ \xbm^{-}(Y)$ for ``$ \xcq $''.

If $y \xce \xbm^{+}(X),$ then $Th(y)$ is inconsistent with $Th(
\xbm^{+}(X)),$ as $ \xbm^{+}$ is
definability preserving, so by classical compactness there is a suitable
finite
sequence $ \xbs $ with $M( \xbs ) \xcs \xbm^{+}(X)= \xCQ,$ but this was
excluded by the induction
hypothesis (1). So $y \xbe \xbm^{+}(X).$

Suppose $y \xbe \xbm^{+}(Y),$ but $x \xce \xbm^{-}(Y).$ So $y \xbe
\xbm^{+}(Y)$ and $y \xbe \xbm^{+}(X),$ and $Y=M( \xbf )$ for
some $ \xbf,$ so there will be a suitable finite sequence $ \xbs $ s.t.
for all $ \xbs ' $ with
$ \xbs \xcc \xbs ' $ $M( \xbs ' ) \xcs X \xcc M( \xbf )=Y,$
and by our construction $x \xbe \xbm^{-}(M( \xbs ' ) \xcs X).$ As $y \xbe
\xbm^{+}(X) \xcs \xbm^{+}(Y) \xcs (M( \xbs ' ) \xcs X),$
$ \xbm^{+}(M( \xbs ' ) \xcs X) \xcc \xbm^{+}(Y),$ so by $( \xbm^{-}4)$ $
\xbm^{-}(M( \xbs ' ) \xcs X) \xcc \xbm^{-}(Y),$ so $x \xbe \xbm^{-}(Y),$
$contradiction.$

We do now this construction for all strategies. Obviously, this does not
modify our results.

This finishes the completeness proof. $ \xcz $
\\[3ex]

As we postulated definability preservation, there are no problems to
translate
the result into logic. (Note that $ \xbn $ was applied to formula defined
model sets,
but the resulting sets were perhaps theory defined model sets.)

\section{
DEFINABILITY PRESERVATION
}

% (+++*** Orig.:   (4)  DEFINABILITY PRESERVATION )

\subsection{
General remarks, affected conditions
}

% (+++*** Orig.:   (4.1)  General remarks, affected conditions  (-> DCB 5.1.3) )

We assume now $ \xdy \xcc \xdp (Z)$ to be closed under arbitrary
intersections (this is
used for the definition of $ \wt{.})$ and finite unions, and $ \xCQ,Z
\xbe \xdy.$ This holds,
of course, for $ \xdy = \xdD_{ \xdl },$ $ \xdl $ any propositional
language.

The aim of Sections 4.1 and 4.2 is to present the results of [Sch04]
connected to
problems of definability preservation in a uniform way, stressing the
crucial
condition $ \wt{X} \xcs \wt{Y}= \wt{X \xcs Y}.$ This presentation shall
help and guide future research
concerning similar problems.

For motivation, we first consider the problem with definability
preservation for
the rules

(PR) $ \ol{ \ol{T \xcv T' } }$ $ \xcc $ $ \ol{ \ol{ \ol{T} } \xcv T' },$
and

$( \xcn =)$ $T \xcl T',$ $Con( \ol{ \ol{T' } },T)$ $ \xcp $ $ \ol{ \ol{T}
}$ $=$ $ \ol{ \ol{ \ol{T' } } \xcv T}$ holds.

which are consequences of

$( \xbm PR)$ $X \xcc Y$ $ \xcp $ $ \xbm (Y) \xcs X \xcc \xbm (X)$ or

$( \xbm =)$ $X \xcc Y,$ $ \xbm (Y) \xcs X \xEd \xCQ $ $ \xcp $ $ \xbm (Y)
\xcs X= \xbm (X)$ respectively

and definability preservation.

First, in the general case without definability preservation, (PR) fails,
and
in the ranked case, $( \xcn =)$ may fail. So failure is not just a
consequence of
the very liberal definition of general preferential structures.

\be

$\hspace{0.01em}$

% (+++*** Orig. No.:  Example 4.1: )

(1) This example was first given in [Sch92].

Let $v( \xdl ):=\{p_{i}:i \xbe \xbo \},$ $n,n' \xbe M_{ \xdl }$ be defined
by $n \xcm \{p_{i}:i \xbe \xbo \},$
$n' \xcm \{ \xCN p_{0}\} \xcv \{p_{i}:0<i< \xbo \}.$
Let $ \xdm:=<M_{ \xdl }, \xeb >$ where only $n \xeb n',$ i.e. just two
models are
comparable. Let $ \xbm:= \xbm_{ \xdm },$ and $ \xcn $ be defined as usual
by $ \xbm.$

Set $T:= \xCQ,$ $T':=\{p_{i}:0<i< \xbo \}.$ We have $M_{T}=M_{ \xdl },$
$ \xbm (M_{T})=M_{ \xdl }-\{n' \},$ $M_{T' }=\{n,n' \},$
$ \xbm (M_{T' })=\{n\}.$ So by Example 1.1, $ \xdm $ is not
definability preserving, and, furthermore,
$ \ol{ \ol{T} }= \ol{T},$ $ \ol{ \ol{T' } }= \ol{\{p_{i}:i< \xbo \}},$
$so$ $p_{0} \xbe \ol{ \ol{T \xcv T' } },$ $but$ $ \ol{ \ol{ \ol{T} } \xcv
T' }= \ol{ \ol{T} \xcv T' }= \ol{T' },$ $so$ $p_{0} \xce
 \ol{ \ol{ \ol{T} } \xcv T' },$
contradicting (PR).

(2)
Take $\{p_{i}:i \xbe \xbo \}$ and put $m:=m_{ \xcU p_{i}},$ the model
which makes all $p_{i}$ true, in the top
layer, all the other in the bottom layer. Let $m' \xEd m,$ $T':= \xCQ,$
$T:=Th(m,m' ).$ Then
Then $ \ol{ \ol{T' } }=T',$ so $Con( \ol{ \ol{T' } },T),$ $ \ol{ \ol{T}
}=Th(m' ),$ $ \ol{ \ol{ \ol{T' } } \xcv T}=T.$

$ \xcz $
\\[3ex]

\ee

We recall from [Sch04] the following Definition and part of the following
Fact:

\bd

$\hspace{0.01em}$

% (+++*** Orig. No.:  Definition 4.1: )

Let $ \xdy \xcc \xdp (Z)$ be given and closed under arbitrary
intersections.

(1) For $A \xcc Z,$ let $ \wt{A}$ $:=$ $ \xcS \{X \xbe \xdy:A \xcc X\}.$

(2) For $B \xbe \xdy,$ we call $A \xcc B$ a small subset of $B$ iff there
is no $X \xbe \xdy $ s.t.
$B-A \xcc X \xcb B.$

\ed

Intuitively, $Z$ is the set of all models for $ \xdl,$ $ \xdy $ is $
\xdD_{ \xdl }$, and $ \wt{A}=M(Th(A)),$
this is the intended application.

\bfa

$\hspace{0.01em}$

% (+++*** Orig. No.:  Fact 4.1: )

(1) If $ \xdy \xcc \xdp (Z)$ is closed under arbitrary intersections and
finite unions,
$Z \xbe \xdy,$ $X,Y \xcc Z,$ then condition

$(~ \xcv )$ $ \wt{X \xcv Y}$ $=$ $ \wt{X} \xcv \wt{Y}$

holds, as do the following trivial ones:

(2) $X=Y$ $ \xcp $ $ \wt{X}= \wt{Y},$ but not conversely,

(3) $ \wt{X} \xcc Y$ $ \xcp $ $X \xcc Y,$ but not conversely,

(4) $X \xcc \wt{Y}$ $ \xcp $ $ \wt{X} \xcc \wt{Y},$

(4a) $ \wt{X \xcs Y} \xcc \wt{X} \xcs \wt{Y}.$

In the intended application, the following hold:

(5) $Th(X)$ $=$ $Th( \wt{X}),$

(6) $ \wt{A} \xcs M( \xbq )= \wt{A \xcs M( \xbq )},$

(7) $ \wt{A}-M( \xbf )= \wt{A-M( \xbf )},$

(8) $ \wt{A}- \wt{B} \xcc \wt{A-B},$

(9) Even if $A= \wt{A},$ $B= \wt{B},$ it is not necessarily true that $
\wt{A-B} \xcc \wt{A}- \wt{B}.$

\efa

\paragraph{
Proof:
}

$\hspace{0.01em}$

(2), (3), (4), (5) are trivial.

(1) Let $ \xdy (U):=\{X \xbe \xdy:U \xcc X\}.$ If $A \xbe \xdy (X \xcv
Y),$ then $A \xbe \xdy (X)$ and $A \xbe \xdy (Y),$ so
$ \wt{X \xcv Y}$ $ \xcd $ $ \wt{X} \xcv \wt{Y}.$ $If$ $A \xbe \xdy (X)$
$and$ $B \xbe \xdy (Y),$ $then$ $A \xcv B \xbe \xdy (X \xcv Y),$ $so$ $
\wt{X \xcv Y}$ $ \xcc $
$ \wt{X} \xcv \wt{Y} \xbe \xdy.$

(4a) Let $X',Y' \xbe \xdy,$ $X \xcc X',$ $Y \xcc Y',$ then $X \xcs Y
\xcc X' \xcs Y',$ so $ \wt{X \xcs Y} \xcc \wt{X} \xcs \wt{Y}.$

(6) $ \wt{A} \xcs M( \xbq ) \xcd \wt{A \xcs M( \xbq )}$ by (4a). For ``$
\xcc $'': Let $A' \xcd A \xcs M( \xbq ),$ $A' \xbe \xdy,$ then
$A' \xcv M( \xCN \xbq ) \xbe \xdy,$ $A' \xcv M( \xCN \xbq ) \xcd A,$ $(A'
\xcv M( \xCN \xbq )) \xcs M( \xbq ) \xcc A'.$ So $ \wt{A} \xcs M( \xbq )
\xcc \wt{A \xcs M( \xbq )}.$

(7) $X-M( \xbf )=X \xcs M( \xCN \xbf ),$ and (6).

(8) Let $x \xbe \wt{A}- \wt{B},$ but $x \xce \wt{A-B}.$ So $x \xcm Th(A),$
and there is $ \xbf \xbe Th(B).x \xcM \xbf,$ and
there is $ \xbq $ s.t. $A-B \xcm \xbq,$ $x \xcM \xbq.$ By $B \xcc M(
\xbf ),$ $A \xcs M( \xCN \xbf ) \xcm \xbq,$ as $x \xcm Th(A),$
$x \xcm \xCN \xbf,$ so $x \xcm \xbq,$ $contradiction.$

(9) Set $A:=M_{ \xdl },$ $B:=\{m\}$ for $m \xbe M_{ \xdl }$ arbitrary, $
\xdl $ infinite. So $A= \wt{A},$ $B= \wt{B},$ but
$ \wt{A-B}=A \xEd A-B.$

$ \xcz $
\\[3ex]
%                            ~  ~ ~~~~
%                            ~  ~ ~~~~
%  (4.2)  Central condition (X%sY=X%sY)  (-> DCB 5.1.3)
% %
% =====================================================
\subsection{Central condition (intersection)}

We analyze the problem of (PR), seen in Example 4.1 (1) above, working in
the
intended application.

(PR) is equivalent to $M( \ol{ \ol{T} } \xcv T' )$ $ \xcc $ $M( \ol{ \ol{T
\xcv T' } }).$ To show (PR) from $( \xbm PR),$ we argue
as follows, the crucial point is marked by ``?'':

$M( \ol{ \ol{T \xcv T' } })$ $=$ $M(Th( \xbm (M_{T \xcv T' })))$ $=$ $
\wt{ \xbm (M_{T \xcv T' })}$ $ \xcd $ $ \xbm (M_{T \xcv T' })$ $=$ $ \xbm
(M_{T} \xcs M_{T' })$
$ \xcd $ (by $( \xbm PR))$ $ \xbm (M_{T}) \xcs M_{T' }$? $ \wt{ \xbm
(M_{T})} \xcs M_{T' }$ $=$ $M(Th( \xbm (M_{T}))) \xcs M_{T' }$ $=$ $M(
\ol{ \ol{T} }) \xcs M_{T' }$ $=$
$M( \ol{ \ol{T} } \xcv T' ).$ If $ \xbm $ is definability preserving, then
$ \xbm (M_{T})$ $=$ $ \wt{ \xbm (M_{T})},$ so ``?'' above
is equality, and everything is fine. In general, however, we have only
$ \xbm (M_{T})$ $ \xcc $ $ \wt{ \xbm (M_{T})},$ and the argument
collapses.

But it is not necessary to
impose $ \xbm (M_{T})$ $=$ $ \wt{ \xbm (M_{T})},$ as we still have room to
move: $ \wt{ \xbm (M_{T \xcv T' })}$ $ \xcd $ $ \xbm (M_{T \xcv T' }).$
(We do not consider here $ \xbm (M_{T} \xcs M_{T' })$ $ \xcd $ $ \xbm
(M_{T}) \xcs M_{T' }$ as room to move,
as we are now interested
only in questions related to definability preservation.) If we had
$ \wt{ \xbm (M_{T})} \xcs M_{T' }$ $ \xcc $ $ \wt{ \xbm (M_{T}) \xcs M_{T'
}}$ , we could use $ \xbm (M_{T}) \xcs M_{T' }$ $ \xcc $ $ \xbm
(M_{T} \xcs M_{T' })$ $=$
$ \xbm (M_{T \xcv T' })$ and monotony of $ \wt{.}$ to obtain $ \wt{ \xbm
(M_{T})} \xcs M_{T' }$ $ \xcc $ $ \wt{ \xbm (M_{T}) \xcs M_{T' }}$ $ \xcc
$
$ \wt{ \xbm (M_{T} \xcs M_{T' })}$ $=$ $ \wt{ \xbm (M_{T \xcv T' })}.$
If, for instance, $T' =\{ \xbq \},$ we have $ \wt{ \xbm (M_{T})}
\xcs M_{T' }$ $=$
$ \wt{ \xbm (M_{T}) \xcs M_{T' }}$ by Fact 4.1 (6). Thus,
definability preservation is not the
only solution to the problem.

We have seen above that $ \wt{X \xcv Y}$ $=$ $ \wt{X} \xcv \wt{Y},$
moreoever X-Y $=$ $X \xcs \xdC Y$ $( \xdC Y$ the set
complement of Y), so, when considering boolean expressions of model sets
(as we do in usual properties describing logics), the central question is

whether

$(~ \xcs )$ $ \wt{X \xcs Y}$ $=$ $ \wt{X} \xcs \wt{Y}.$

We take a closer look at this question.
``$ \xcc $'' holds by Fact 4.1 (6). Using $(~ \xcv )$ and monotony
of $ \wt{.}$, we have $ \wt{X} \xcs \wt{Y}$ $=$ $ \wt{((X \xcs Y) \xcv
(X-Y))}$ $ \xcs $ $ \wt{((X \xcs Y) \xcv (Y-X))}$ $=$
$( \wt{(X \xcs Y)} \xcv \wt{(X-Y)})$ $ \xcs $ $( \wt{(X \xcs Y)} \xcv
\wt{(Y-X)})$ $=$
$ \wt{X \xcs Y}$ $ \xcv $ $( \wt{X \xcs Y} \xcs \wt{Y-X})$ $ \xcv $ $(
\wt{X \xcs Y} \xcs \wt{X-Y})$ $ \xcv $ $( \wt{X-Y} \xcs \wt{Y-X})$ $=$ $
\wt{X \xcs Y}$ $ \xcv $ $( \wt{X-Y} \xcs \wt{Y-X}),$

thus $ \wt{X} \xcs \wt{Y}$ $ \xcc $ $ \wt{X \xcs Y}$ iff

$(~ \xcs ' )$ $ \wt{Y-X}$ $ \xcs $ $ \wt{X-Y}$ $ \xcc $ $ \wt{X \xcs Y}$
holds.

Intuitively speaking, the condition holds
iff we cannot approximate any element both from X-Y and X-Y, which cannot
be
approximated from $X \xcs Y,$ too.

Note that in above Example 4.1 $X:= \xbm (M_{T})=M_{ \xdl }-\{n' \},$
$Y:=M_{T' }=\{n,n' \},$
$ \wt{X-Y}=M_{ \xdl },$ $ \wt{Y-X}=\{n' \},$ $ \wt{X \xcs Y}=\{n\},$ and
$ \wt{X} \xcs \wt{Y}=\{n,n' \}.$

We consider now particular cases:

(1) In particular, if $X \xcs Y= \xCQ,$ by $ \xCQ \xbe \xdy,$ $(~ \xcs
)$ holds iff $ \wt{X} \xcs \wt{Y}= \xCQ.$

(2) If $X \xbe \xdy $ and $Y \xbe \xdy,$ then $ \wt{X-Y} \xcc X$ and $
\wt{Y-X} \xcc Y,$ so $ \wt{X-Y} \xcs \wt{Y-X} \xcc X \xcs Y \xcc \wt{X
\xcs Y}$ and
$(~ \xcs )$ trivially holds.

(3) $X \xbe \xdy $ and $ \xdC X \xbe \xdy $ together also suffice - in
these cases
$ \wt{Y-X}$ $ \xcs $ $ \wt{X-Y}$ $=$ $ \xCQ:$ $ \wt{Y-X}= \wt{Y \xcs \xdC
X} \xcc \xdC X,$ and $ \wt{X-Y} \xcc X,$ $so$ $ \wt{Y-X} \xcs \wt{X-Y}
\xcc X \xcs \xdC X= \xCQ \xcc \wt{X \xcs Y}.$
(The same holds, of course, for Y.)
(In the intended application, such $X$ will be $M( \xbf )$ for some
formula $ \xbf.$ But, a
warning, $ \xbm (M( \xbf ))$ need not again be the $M( \xbq )$ for some $
\xbq.)$

We turn to the properties of various structures and apply our results.
%  (4.2.1)  The minimal variant
%  (4.2.1)  The minimal variant
% %
% =============================
\subsubsection{The minimal variant}

We now take a look at other frequently used logical conditions. First, in
the
context on nonmonotonic logics, the following rules will always hold in
smooth
preferential structures, even if we consider full theories, and not
necessarily
definability preserving structures:

\bfa

$\hspace{0.01em}$

% (+++*** Orig. No.:  Fact 4.2: )

Also for full theories, and not necessarily definability preserving
structures
hold:

(1) (LLE), (RW), (AND), (REFLEX), by definition and $( \xbm \xcc ),$

(2) (OR),

(3) (CM) in smooth structures,

(4) the infinitary version of (CUM) in smooth structures.
In definability preserving structures, but also when considering only
formulas
hold:

(5) (PR),

(6) $( \xcn =)$ in ranked structures.

\efa

\paragraph{
Proof:
}

$\hspace{0.01em}$

We use the corresponding algebraic properties.

(1) trivial.

(2) We have to show $T \xcn \xbf,$ $T' \xcn \xbf $ $ \xcp $ $T \xco T'
\xcn \xbf:$ $ \xbm (T) \xcm \xbf,$ $ \xbm (T' ) \xcm \xbf.$ Thus
$ \xbm (T \xco T' )= \xbm (M(T) \xcv M(T' )) \xcc \xbm (T) \xcv \xbm (T' )
\xcm \xbf.$

(3) We have to show $T \xcn \xbb,$ $T \xcn \xbg $ $ \xcp $ $T \xcv \xbb
\xcn \xbg:$ $ \xbm (T) \xcc M(T \xcv \xbb ) \xcc M(T)$ (as
$ \xbm (T) \xcc M(T)$ and $ \xbm (T) \xcc M( \xbb )),$ so $ \xbm (T)= \xbm
(T \xcv \xbb ).$

(4)
Let $T$ $ \xcc $ $ \ol{T' }$ $ \xcc $ $ \ol{ \ol{T} }.$ Thus by $( \xbm
Cum)$ and $ \xbm (M_{T})$ $ \xcc $ $M_{ \ol{ \ol{T} }}$ $ \xcc $ $M_{T' }$
$ \xcc $ $M_{T}$ $ \xbm (M_{T})= \xbm (M_{T' }),$
so $ \ol{ \ol{T} }$ $=$ $Th( \xbm (M_{T}))$ $=$ $Th( \xbm (M_{T' }))$ $=$
$ \ol{ \ol{T' } }.$ (The proof given in [Sch04] uses
definability preservation, but this is not necessary, as we see here.)

(5)
See above discussion.

(6)
$T \xcl T' $ iff $M(T) \xcc M(T' ).$ $Con( \ol{ \ol{T' } },T)$ iff $M(
\ol{ \ol{T' } }) \xcs M(T) \xEd \xCQ $ iff $ \wt{ \xbm (T' )} \xcs M(T)
\xEd \xCQ $ iff
$ \wt{ \xbm (T' ) \xcs M(T)} \xEd \xCQ $ $iff$ $ \xbm (T' ) \xcs M(T) \xEd
\xCQ,$ $if$ $ \xbm (T' ) \xbe \xdy $ $or$ $T= \xbf,$
as we saw above. So by rankedness $ \xbm (T)= \xbm (T' ) \xcs M(T)$ $ \xcp
$
$ \ol{ \ol{T} }=Th( \xbm (T))=Th( \xbm (T' ) \xcs M(T)),$ and
$Th( \xbm (T' ) \xcs M(T))= \ol{ \ol{ \ol{T' } } \xcv T}$ if $T= \xbf $ or
$ \xbm (T' ) \xbe \xdy $ again.

$ \xcz $
\\[3ex]

We turn to theory revision. The following definition and example, taken
from
[Sch04]
shows, that the usual AGM axioms for theory revision fail in distance
based
structures in the general case, unless we require definability
preservation.

\bd

$\hspace{0.01em}$

% (+++*** Orig. No.:  Definition 4.2: )

We summarize the AGM postulates $(K*7)$ and $(K*8)$ in $(*4):$

$(*4)$ If $T*T' $ is consistent with $T'',$ then $T*(T' \xcv T'' )$ $=$ $
\ol{(T*T' ) \xcv T'' }.$

\ed

\be

$\hspace{0.01em}$

% (+++*** Orig. No.:  Example 4.2: )

Consider an infinite propositional language $ \xdl.$

Let $X$ be an infinite set of models, $m,$ $m_{1},$ $m_{2}$ be models for
$ \xdl.$
Arrange the models of $ \xdl $ in the real plane s.t. all $x \xbe X$ have
the same
distance $<2$ (in the real plane) from $m,$ $m_{2}$ has distance 2 from
$m,$ and $m_{1}$ has
distance 3 from $m.$

Let $T,$ $T_{1},$ $T_{2}$ be complete (consistent) theories, $T' $ a
theory with infinitely
many models, $M(T)=\{m\},$ $M(T_{1})=\{m_{1}\},$ $M(T_{2})=\{m_{2}\}.$
$M(T' )=X \xcv \{m_{1},m_{2}\},$ $M(T'' )=\{m_{1},m_{2}\}.$

Assume $Th(X)=T',$ so $X$ will not be definable by a theory.

Then $M(T) \xfA M(T' )=X,$ but $T*T' =Th(X)=T'.$ So $T*T' $ is consistent
with $T'',$ and
$ \ol{(T*T' ) \xcv T'' }=T''.$ But $T' \xcv T'' =T'',$ and $T*(T' \xcv
T'' )=T_{2} \xEd T'',$ contradicting $(*4).$

$ \xcz $
\\[3ex]

\ee

We show now that the version with formulas only holds here, too, just as
does
above (PR), when we consider formulas only - this is needed below for $T''
$ only.
This was already shown in [Sch04], we give now a proof based on our new
principles.

\bfa

$\hspace{0.01em}$

% (+++*** Orig. No.:  Fact 4.3: )

$(*4)$ holds when considering only formulas.

\efa

\paragraph{
Proof:
}

$\hspace{0.01em}$

When we fix the left hand side, the structure is ranked, so $Con(T*T',T''
)$
implies $(M_{T} \xfA M_{T' }) \xcs M_{T'' } \xEd \xCQ $ by $T'' =\{ \xbq
\}$ and thus $M_{T} \xfA M_{T' \xcv T'' }$ $=$ $M_{T} \xfA (M_{T' } \xcs
M_{T'' })$ $=$
$(M_{T} \xfA M_{T' }) \xcs M_{T'' }.$
So $M(T*(T' \xcv T'' ))$ $=$ $ \wt{M_{T} \xfA M_{T' \xcv T'' }}$ $=$ $
\wt{(M_{T} \xfA M_{T' }) \xcs M_{T'' }}$ $=$ (by $T'' =\{ \xbq \},$ see
above)
$ \wt{(M_{T} \xfA M_{T' })} \xcs \wt{M_{T'' }}$ $=$ $ \wt{(M_{T} \xfA
M_{T' })} \xcs M_{T'' }$ $=$ $M((T*T' ) \xcv T'' ),$ and $T*(T' \xcv T''
)= \ol{(T*T' ) \xcv T'' }.$
$ \xcz $
\\[3ex]
%  (4.2.2)  The limit variant
%  (4.2.2)  The limit variant
% %
% ===========================
\subsubsection{The limit variant}

We begin with some simple logical facts about the limit version.

We abbreviate $ \xbL (T):= \xbL (M(T))$ etc.

\bfa

$\hspace{0.01em}$

% (+++*** Orig. No.:  Fact 4.4: )

(1) $A \xbe \xbL (T)$ $ \xcp $ $M( \ol{ \ol{T} }) \xcc \wt{A}$

(2) $M( \ol{ \ol{T} })$ $=$ $ \xcS \{ \wt{A}:A \xbe \xbL (T)\}$

(2a) $M( \ol{ \ol{T' } }) \xcm \xbs $ $ \xcp $ $ \xcE B \xbe \xbL (T' ).
\wt{B} \xcm \xbs $

(3) $M( \ol{ \ol{T' } }) \xcs M(T) \xcm \xbs $ $ \xcp $ $ \xcE B \xbe \xbL
(T' ). \wt{B} \xcs M(T) \xcm \xbs.$

\efa

\paragraph{
Proof:
}

$\hspace{0.01em}$

(1) Let $M( \ol{ \ol{T} }) \xcC \wt{A},$ so there is $ \xbf,$ $ \wt{A}
\xcm \xbf,$ so $A \xcm \xbf,$ but $M( \ol{ \ol{T} }) \xcM \xbf,$ so $T
\xcN \xbf,$ $contradiction.$

(2) ``$ \xcc $'' by (1). Let $x \xbe \xcS \{ \wt{A}:A \xbe \xbL (T)\}$ $
\xcp $ $ \xcA A \xbe \xbL (T).x \xcm Th(A)$ $ \xcp $ $x \xcm \ol{ \ol{T}
}.$

(2a) $M( \ol{ \ol{T' } }) \xcm \xbs $ $ \xcp $ $T' \xcn \xbs $ $ \xcp $ $
\xcE B \xbe \xbL (T' ).B \xcm \xbs.$ But $B \xcm \xbs $ $ \xcp $ $ \wt{B}
\xcm \xbs.$

(3) $M( \ol{ \ol{T' } }) \xcs M(T) \xcm \xbs $ $ \xcp $ $ \ol{ \ol{T' } }
\xcv T \xcl \xbs $ $ \xcp $ $ \xcE \xbt_{1} \Xl  \xbt_{n} \xbe \ol{ \ol{T'
} }$ s.t. $T \xcv \{ \xbt_{1}, \Xl, \xbt_{n}\} \xcl
 \xbs,$
so $ \xcE B \xbe \xbL (T' ).Th(B) \xcv T \xcl \xbs.$ So $M(Th(B)) \xcs
M(T) \xcm \xbs $ $ \xcp $ $ \wt{B} \xcs M(T) \xcm \xbs.$

$ \xcz $
\\[3ex]

We saw in Example 3.2 and its discussion the problems which might arise in
the limit version, even if the algebraic behaviour is correct.

This analysis leads us to consider the following facts:

\bfa

$\hspace{0.01em}$

% (+++*** Orig. No.:  Fact 4.5: )

(1) Let $ \xcA B \xbe \xbL (T' ) \xcE A \xbe \xbL (T).A \xcc B \xcs M(T),$
then $ \ol{ \ol{ \ol{T' } } \xcv T} \xcc \ol{ \ol{T} }.$

Let, in addition, $\{B \xbe \xbL (T' ): \wt{B} \xcs \wt{M(T)}= \wt{B \xcs
M(T)}\}$ be cofinal in $ \xbL (T' ).$ Then

(2) $Con( \ol{ \ol{T' } },T)$ implies $ \xcA A \xbe \xbL (T' ).A \xcs M(T)
\xEd \xCQ.$

(3) $ \xcA A \xbe \xbL (T) \xcE B \xbe \xbL (T' ).B \xcs M(T) \xcc A$
implies $ \ol{ \ol{T} } \xcc \ol{ \ol{ \ol{T' } } \xcv T}.$

\efa

Note that $M(T)= \wt{M(T)},$ so we could also have written $ \wt{B} \xcs
M(T)= \wt{B \xcs M(T)},$ but above
way of writing stresses more the essential condition $ \wt{X} \xcs \wt{Y}=
\wt{X \xcs Y}.$

\paragraph{
Proof:
}

$\hspace{0.01em}$

(1)
Let $ \ol{ \ol{T' } } \xcv T \xcl \xbs,$ so $ \xcE B \xbe \xbL (T' ).
\wt{B} \xcs M(T) \xcm \xbs $ by Fact 4.4, (3) above (using
compactness). Thus $ \xcE A \xbe \xbL (T).A \xcc B' \xcs M(T) \xcm \xbs $
by prerequisite, so $ \xbs \xbe \ol{ \ol{T} }.$

(2)
Let $Con( \ol{ \ol{T' } },T),$ so $M( \ol{ \ol{T' } }) \xcs M(T) \xEd \xCQ
.$ $M( \ol{ \ol{T' } })= \xcS \{ \wt{A}:A \xbe \xbL (T' )\}$ by Fact 4.4
(2), so
$ \xcA A \xbe \xbL (T' ). \wt{A} \xcs M(T) \xEd \xCQ.$ As cofinally often
$ \wt{A} \xcs M(T)= \wt{A \xcs M(T)},$
$ \xcA A \xbe \xbL (T' ). \wt{A \xcs M(T)} \xEd \xCQ,$ so $ \xcA A \xbe
\xbL (T' ).A \xcs M(T) \xEd \xCQ $ by $ \wt{ \xCQ }= \xCQ.$

(3)
Let $ \xbs \xbe \ol{ \ol{T} },$ so $T \xcn \xbs,$ so $ \xcE A \xbe \xbL
(T).A \xcm \xbs,$ so $ \xcE B \xbe \xbL (T' ).B \xcs M(T) \xcc A$ by
prerequisite, so $ \xcE B \xbe \xbL (T' ).B \xcs M(T) \xcc A \xcu \wt{B}
\xcs \wt{M(T)}= \wt{B \xcs M(T)}.$ So for
such $B$ $ \wt{B} \xcs \wt{M(T)}= \wt{B \xcs M(T)} \xcc \wt{A} \xcm \xbs
.$ By
Fact 4.4 (1) $M( \ol{ \ol{T' } }) \xcc \wt{B},$ so $M( \ol{ \ol{T' } })
\xcs M(T) \xcm \xbs,$ so $ \ol{ \ol{T' } } \xcv T \xcl \xbs.$

$ \xcz $
\\[3ex]

We obtain now as easy corollaries of a more general situation the
following
properties shown in [Sch04] by direct proofs. Thus, we have the
trivialization
results shown there.

\bco

$\hspace{0.01em}$

% (+++*** Orig. No.:  Corollary 4.6: )

Let the structure be transitive.

(1) Let $\{B \xbe \xbL (T' ): \wt{B} \xcs \wt{M(T)}= \wt{B \xcs M(T)}\}$
be cofinal in $ \xbL (T' ),$ then:

(PR) $T \xcl T' $ $ \xcp $ $ \ol{ \ol{T} }$ $ \xcc $ $ \ol{ \ol{ \ol{T' }
} \xcv T}.$

(2) $ \ol{ \ol{ \xbf \xcu \xbf ' } } \xcc \ol{ \ol{ \ol{ \xbf } } \xcv \{
\xbf ' \}}$

If the structure is ranked, then also:

(3) Let $\{B \xbe \xbL (T' ): \wt{B} \xcs \wt{M(T)}= \wt{B \xcs M(T)}\}$
be cofinal in $ \xbL (T' ),$ then:

$( \xcn =)$ $T \xcl T',$ $Con( \ol{ \ol{T' } },T)$ $ \xcp $ $ \ol{ \ol{T}
}$ $=$ $ \ol{ \ol{ \ol{T' } } \xcv T}$

(4) $T \xcN \xCN \xbg $ $ \xcp $ $ \ol{ \ol{T} }= \ol{ \ol{T \xcv \{ \xbg
\}} }$

\eco

\paragraph{
Proof:
}

$\hspace{0.01em}$

(1) $ \xcA A \xbe \xbL (M(T)) \xcE B \xbe \xbL (M(T' )).B \xcs M(T) \xcc
A$ by Fact 3.3 (2).
So the result follows from Fact 4.5 (3).

(2) Set $T':=\{ \xbf \},$ $T:=\{ \xbf, \xbf ' \}.$ Then for $B \xbe \xbL
(T' )$ $ \wt{B} \xcs M(T)= \wt{B} \xcs M( \xbf ' )= \wt{B \xcs M( \xbf '
)}$
by Fact 4.1 (6), so the result follows by (1).

(3) Let $Con( \ol{ \ol{T' } },T),$ then by Fact 4.5 (2) $ \xcA A \xbe \xbL
(T' ).A \xcs M(T) \xEd \xCQ,$ so by
Fact 3.3 (4) $ \xcA B \xbe \xbL (T' ) \xcE A \xbe \xbL (T).A \xcc B \xcs
M(T),$ so $ \ol{ \ol{ \ol{T' } } \xcv T} \xcc \ol{ \ol{T} }$ by Fact 4.3
(1).
The other direction follows from (1).

(4) Set $T:=T' \xcv \{ \xbg \}.$ Then for $B \xbe \xbL (T' )$ $ \wt{B}
\xcs M(T)= \wt{B} \xcs M( \xbg )= \wt{B \xcs M( \xbg )}$ again
by Fact 4.1 (6), so the result follows from (3).

$ \xcz $
\\[3ex]

\subsection{
A simplification of [Sch04]
}

% (+++*** Orig.:   (4.3)  A simplification of [Sch04] ( p p(U))  (-> DCB 5.2.1) )

Note that in Sections 3.2 and 3.3 of [Sch04],
as well as in Proposition 4.2.2 of [Sch04] we
have characterized $ \xbm: \xdy \xcp \xdy $ or $ \xfA: \xdy \xCK \xdy
\xcp \xdy,$ but a closer inspection of the
proofs shows that the destination can as well be assumed $ \xdp (Z),$
consequently
we can simply re-use above algebraic representation results also for the
not
definability preserving case. (Note that the easy direction of all
these results work for destination $ \xdp (Z),$ too.) In particular, also
the
proof for the not definability preserving case of revision in [Sch04] can
be simplified - but we will not go into details here.

$( \xcv )$ and $( \xcS )$ are assumed to hold now - we need $( \xcS )$ for
$ \wt{.}.$

The central functions and conditions to consider are summarized in the
following definition.

\bd

$\hspace{0.01em}$

% (+++*** Orig. No.:  Definition 4.3: )

Let $ \xbm: \xdy \xcp \xdy,$ we define $ \xbm_{i}: \xdy \xcp \xdp (Z):$

$ \xbm_{0}(U)$ $:=$ $\{x \xbe U:$ $ \xCN \xcE Y \xbe \xdy (Y \xcc U$ and
$x \xbe Y- \xbm (Y))\},$

$ \xbm_{1}(U)$ $:=$ $\{x \xbe U:$ $ \xCN \xcE Y \xbe \xdy ( \xbm (Y) \xcc
U$ and $x \xbe Y- \xbm (Y))\},$

$ \xbm_{2}(U)$ $:=$ $\{x \xbe U:$ $ \xCN \xcE Y \xbe \xdy ( \xbm (U \xcv
Y) \xcc U$ and $x \xbe Y- \xbm (Y))\}$

(note that we use $( \xcv )$ here),

$ \xbm_{3}(U)$ $:=$ $\{x \xbe U:$ $ \xcA y \xbe U.x \xbe \xbm (\{x,y\})\}$

(we use here $( \xcv )$ and that singletons are in $ \xdy )$

$( \xbm PR0)$ $ \xbm (U)- \xbm_{0}(U)$ is small,

$( \xbm PR1)$ $ \xbm (U)- \xbm_{1}(U)$ is small,

$( \xbm PR2)$ $ \xbm (U)- \xbm_{2}(U)$ is small,

$( \xbm PR3)$ $ \xbm (U)- \xbm_{3}(U)$ is small.

$( \xbm PR0)$ with its function will be the one to consider for general
preferential
structures, $( \xbm PR2)$ the one for smooth structures. Unfortunately, we
cannot use
$( \xbm PR0)$ in the smooth case, too, as Example 4.3 below will show.
This sheds some
doubt on the possibility to find an easy common approach to all cases of
not
definability preserving preferential, and perhaps other, structures. The
next
best guess, $( \xbm PR1)$ will not work either, as the same example shows
- or by
Fact 4.7 (10), if $ \xbm $ satisfies $( \xbm Cum),$ then $ \xbm_{0}(U)=
\xbm_{1}(U).$ $( \xbm PR3)$ and
$ \xbm_{3}$ are used for ranked structures.

\ed

Note that in our context, $ \xbm $ will not necessarily respect $( \xbm
PR).$ Thus, if e.g.
$x \xbe Y- \xbm (Y),$ and $ \xbm (Y) \xcc U,$ we cannot necessarily
conclude that $x \xce \xbm (U \xcv Y)$ -
the fact that $x$ is minimized in $U \xcv Y$ might be hidden by the bigger
$ \xbm (U \xcv Y).$

The strategy of representation without definability preservation will in
all
cases be very simple: Under sufficient conditions,
among them smallness $( \xbm PRi)$ as described above, the corresponding
function
$ \xbm_{i}$ has all the properties to guarantee representation by a
corresponding
structures, and we can just take our representation theorems for the dp
case,
to show this. Using smallness again, we can show that we have obtained a
sufficient approximation - see Propositions 4.9, 4.10, and 4.15.

We first show some properties for the $ \xbm_{i},$ $i=0,1,2.$ A
corresponding
result for $ \xbm_{3}$ is given in Fact 4.13 below. (The conditions and
results
are sufficiently different for $ \xbm_{3}$ to make a separation more
natural.)

Property (9) of the following Fact 4.7 fails for $ \xbm_{0}$ and $
\xbm_{1},$ as Example
4.3 below will show. We will therefore work in the smooth case with $
\xbm_{2}.$
%  (4.3.1)  The general and the smooth case
%  (4.3.1)  The general and the smooth case
% %
% =========================================
\subsubsection{The general and the smooth case}

\bfa

$\hspace{0.01em}$

% (+++*** Orig. No.:  Fact 4.7: )  ###

(This is partly Fact 5.2.6 in [Sch04].)

Recall that $ \xdy $ is closed under $( \xcv ),$ and $ \xbm: \xdy \xcp
\xdy.$
Let $A,B,U,U',X,Y$ be elements of $ \xdy $ and the $ \xbm_{i}$ be defined
from $ \xbm $ as in
Definition 4.3. $i$ will here be 0,1, or 2, but not 3.

(1) Let $ \xbm $ satisfy $( \xbm \xcc ),$ then $ \xbm_{1}(X) \xcc
\xbm_{0}(X)$ and $ \xbm_{3}(X) \xcc \xbm_{0}(X),$

(2) Let $ \xbm $ satisfy $( \xbm \xcc )$ and $( \xbm Cum),$ then $ \xbm (U
\xcv U' ) \xcc U$ $ \xcr $ $ \xbm (U \xcv U' )= \xbm (U),$

(3) Let $ \xbm $ satisfy $( \xbm \xcc ),$ then $ \xbm_{i}(U) \xcc \xbm
(U),$ and $ \xbm_{i}(U) \xcc U,$

(4) Let $ \xbm $ satisfy $( \xbm \xcc )$ and one of the $( \xbm PRi),$
then
$ \xbm (A \xcv B) \xcc \xbm (A) \xcv \xbm (B),$

(5) Let $ \xbm $ satisfy $( \xbm \xcc )$ and one of the $( \xbm PRi),$
then
$ \xbm_{2}(X) \xcc \xbm_{1}(X),$

(6) Let $ \xbm $ satisfy $( \xbm \xcc ),$ $( \xbm PRi),$ then
$ \xbm_{i}(U) \xcc U' $ $ \xcr $ $ \xbm (U) \xcc U',$

(7) Let $ \xbm $ satisfy $( \xbm \xcc )$ and one of the $( \xbm PRi),$
then
$X \xcc Y,$ $ \xbm (X \xcv U) \xcc X$ $ \xcp $ $ \xbm (Y \xcv U) \xcc Y,$

(8) Let $ \xbm $ satisfy $( \xbm \xcc )$ and one of the $( \xbm PRi),$
then
$X \xcc Y$ $ \xcp $ $X \xcs \xbm_{i}(Y) \xcc \xbm_{i}(X)$ - $( \xbm PR)$
for $ \xbm_{i},$
(more precisely, only for $ \xbm_{2}$ we need the prerequisites, in the
other
cases the definition suffices)

(9) Let $ \xbm $ satisfy $( \xbm \xcc ),$ $( \xbm PR2),$ $( \xbm Cum),$
then
$ \xbm_{2}(X) \xcc Y \xcc X$ $ \xcp $ $ \xbm_{2}(X)= \xbm_{2}(Y)$ - $(
\xbm Cum)$ for $ \xbm_{2}.$

(10) $( \xbm \xcc )$ and $( \xbm Cum)$ for $ \xbm $ entail $ \xbm_{0}(U)=
\xbm_{1}(U).$

\efa

\paragraph{
Proof:
}

$\hspace{0.01em}$

(1) $ \xbm_{1}(X) \xcc \xbm_{0}(X)$ follows from $( \xbm \xcc )$ for $
\xbm.$ For $ \xbm_{2}:$ By $Y \xcc U,$ $U \xcv Y=U,$ so
$ \xbm (U) \xcc U$ by $( \xbm \xcc ).$

(2) $ \xbm (U \xcv U' ) \xcc U \xcc U \xcv U' $ $ \xcp_{( \xbm CUM)}$ $
\xbm (U \xcv U' )= \xbm (U).$

(3) $ \xbm_{i}(U) \xcc U$ by definition. To show $ \xbm_{i}(U) \xcc \xbm
(U),$ take in all three cases
$Y:=U,$ and use for $i=1,2$ $( \xbm \xcc ).$

(4) By definition of $ \xbm_{0},$ we have $ \xbm_{0}(A \xcv B) \xcc A \xcv
B,$ $ \xbm_{0}(A \xcv B) \xcs (A- \xbm (A))= \xCQ,$
$ \xbm_{0}(A \xcv B) \xcs (B- \xbm (B))= \xCQ,$ so $ \xbm_{0}(A \xcv B)
\xcs A \xcc \xbm (A),$ $ \xbm_{0}(A \xcv B) \xcs B \xcc \xbm (B),$ and
$ \xbm_{0}(A \xcv B) \xcc \xbm (A) \xcv \xbm (B).$ By $ \xbm: \xdy \xcp
\xdy $ and $( \xcv ),$ $ \xbm (A) \xcv \xbm (B) \xbe \xdy.$
Moreover, by (3) $ \xbm_{0}(A \xcv B) \xcc \xbm (A \xcv B),$ so $
\xbm_{0}(A \xcv B)$ $ \xcc $
$( \xbm (A) \xcv \xbm (B)) \xcs \xbm (A \xcv B),$ so by (0) $ \xbm_{i}(A
\xcv B) \xcc ( \xbm (A) \xcv \xbm (B)) \xcs \xbm (A \xcv B).$
If $ \xbm (A \xcv B) \xcC \xbm (A) \xcv \xbm (B),$ then $( \xbm (A) \xcv
\xbm (B)) \xcs \xbm (A \xcv B)$ $ \xcb $ $ \xbm (A \xcv B),$
contradicting $( \xbm PRi).$

(5) Let $Y \xbe \xdy,$ $ \xbm (Y) \xcc U,$ $x \xbe Y- \xbm (Y),$ then
(by (4)) $ \xbm (U \xcv Y) \xcc \xbm (U) \xcv \xbm (Y) \xcc U.$

(6) ``$ \xcq $'' by (3). ``$ \xcp '':$ By $( \xbm PRi),$ $ \xbm (U)-
\xbm_{i}(U)$ is small, so there
is no $X \xbe \xdy $ s.t.
$ \xbm_{i}(U) \xcc X \xcb \xbm (U).$ If there were $U' \xbe \xdy $ s.t. $
\xbm_{i}(U) \xcc U',$ but $ \xbm (U) \xcC U',$ then
for $X:=U' \xcs \xbm (U) \xbe \xdy,$ $ \xbm_{i}(U) \xcc X \xcb \xbm (U),$
$contradiction.$

(7) $ \xbm (Y \xcv U)= \xbm (Y \xcv X \xcv U) \xcc_{(4)} \xbm (Y) \xcv
\xbm (X \xcv U) \xcc Y \xcv X=Y.$

(8) For $i=0,1:$ Let $x \xbe X- \xbm_{0}(X),$ then there is A s.t. $A \xcc
X,$ $x \xbe A- \xbm (A),$ so
$A \xcc Y.$ The case $i=1$ is similar. We need here only the definitions.
For $i=2:$ Let $x \xbe X- \xbm_{2}(X),$ A s.t. $x \xbe A- \xbm (A),$ $
\xbm (X \xcv A) \xcc X,$ then by
(7) $ \xbm (Y \xcv A) \xcc Y.$

(9) ``$ \xcc '':$ Let $x \xbe \xbm_{2}(X),$ so $x \xbe Y,$ and $x \xbe
\xbm_{2}(Y)$ by (8).
``$ \xcd '':$ Let $x \xbe \xbm_{2}(Y),$ so $x \xbe X.$ Suppose $x \xce
\xbm_{2}(X),$ so there is $U \xbe \xdy $ s.t.
$x \xbe U- \xbm (U)$ and $ \xbm (X \xcv U) \xcc X.$ Note that by $ \xbm (X
\xcv U) \xcc X$ and (2), $ \xbm (X \xcv U)= \xbm (X).$
Now, $ \xbm_{2}(X) \xcc Y,$ so by (6) $ \xbm (X) \xcc Y,$ thus $ \xbm (X
\xcv U)= \xbm (X) \xcc Y \xcc Y \xcv U \xcc X \xcv U,$
so $ \xbm (Y \xcv U)= \xbm (X \xcv U)= \xbm (X) \xcc Y$ by $( \xbm Cum),$
so $x \xce \xbm_{2}(Y),$ $contradiction.$

(10) $ \xbm_{1}(U) \xcc \xbm_{0}(U)$ by (1). Let $Y$ s.t. $ \xbm (Y) \xcc
U,$ $x \xbe Y- \xbm (Y),$ $x \xbe U.$ Consider
$Y \xcs U,$ $x \xbe Y \xcs U,$ $ \xbm (Y) \xcc Y \xcs U \xcc Y,$ so $ \xbm
(Y)= \xbm (Y \xcs U)$ by $( \xbm Cum),$ and
$x \xce \xbm (Y \xcs U).$ Thus, $ \xbm_{0}(U) \xcc \xbm_{1}(U).$

$ \xcz $
\\[3ex]

\bfa

$\hspace{0.01em}$

% (+++*** Orig. No.:  Fact 4.8: )

$( \xbm PR0),$ $( \xbm PR1),$ $( \xbm PR2)$ are equivalent in the presence
of $( \xbm \xcc ),$ $( \xbm Cum)$ for
$ \xbm.$

(Recall that $( \xcv )$ and $( \xcs )$ are assumed to hold.)

\efa

\paragraph{
Proof:
}

$\hspace{0.01em}$

We first show $( \xbm PR1) \xcr ( \xbm PR2).$

(1) Suppose $( \xbm PR2)$ holds. By $( \xbm PR2)$ and (5), $ \xbm_{2}(U)
\xcc \xbm_{1}(U),$ so
$ \xbm (U)- \xbm_{1}(U) \xcc \xbm (U)- \xbm_{2}(U).$ By $( \xbm PR2),$ $
\xbm (U)- \xbm_{2}(U)$ is small, then so is
$ \xbm (U)- \xbm_{1}(U),$ so $( \xbm PR1)$ holds.

(2) Suppose $( \xbm PR1)$ holds, and $( \xbm PR2)$ fails. By failure of $(
\xbm PR2),$ there is
$X \xbe \xdy $ s.t. $ \xbm_{2}(U) \xcc X \xcb \xbm (U).$ Let $x \xbe \xbm
(U)-$X, as $x \xce \xbm_{2}(U),$ there is $Y$ s.t.
$ \xbm (U \xcv Y) \xcc U,$ $x \xbe Y- \xbm (Y).$ Let $Z:=Y \xcv X.$ By
Fact 4.7 (4)
$ \xbm (Y \xcv X) \xcc \xbm (Y) \xcv \xbm (X),$ so $x \xce \xbm (Y \xcv
X).$ Moreover, $ \xbm (U \xcv X \xcv Y) \xcc \xbm (U \xcv Y) \xcv \xbm
(X)$
by Fact 4.7 (4), $ \xbm (U \xcv Y) \xcc U,$ $ \xbm (X) \xcc X \xcc \xbm
(U) \xcc U$ by prerequisite, so
$ \xbm (U \xcv X \xcv Y) \xcc U \xcc U \xcv Y \xcc U \xcv X \xcv Y,$ so $
\xbm (U \xcv X \xcv Y)= \xbm (U \xcv Y).$ Thus $ \xbm (Y \xcv U \xcv X)$
$=$
$ \xbm (Y \xcv U)$ $ \xcc $ $Y \xcv X$ $ \xcc $ $Y \xcv U \xcv X,$ so $
\xbm (Y \xcv X)$ $=$ $ \xbm (Y \xcv U \xcv X)$ $=$ $ \xbm (Y \xcv U)$ $
\xcc $ $U.$
Thus, $x \xce \xbm_{1}(U),$ and $ \xbm_{1}(U) \xcc X,$ too, a
contradiction.

(3) Finally, by Fact 4.7, (10), $ \xbm_{0}(U)= \xbm_{1}(U)$ if $( \xbm
Cum)$ holds for $ \xbm.$

$ \xcz $
\\[3ex]

Here is an example which shows that Fact 4.7, (9) may fail for $ \xbm_{0}$
and $ \xbm_{1}.$

\be

$\hspace{0.01em}$

% (+++*** Orig. No.:  Example 4.3: )

Consider $ \xdl $ with $v( \xdl ):=\{p_{i}:i \xbe \xbo \}.$ Let $m \xcM
p_{0},$ let $m' \xbe M(p_{0})$ arbitrary.
Make for each $n \xbe M(p_{0})-\{m' \}$ one copy of $m,$ likewise
of $m',$ set $<m,n> \xeb <m',n>$ for all $n,$ and $n \xeb <m,n>,$ $n
\xeb <m',n>$ for all $n.$ The
resulting structure $ \xdz $ is smooth and transitive. Let $ \xdy:=
\xdD_{ \xdl },$ define
$ \xbm (X):= \wt{ \xbm_{ \xdz }(X)}$ for $X \xbe \xdy.$

Let $m' \xbe X- \xbm_{ \xdz }(X).$ Then $m \xbe X,$ or $M(p_{0}) \xcc X.$
In the latter case, as all $m'' $ s.t.
$m'' \xEd m',$ $m'' \xcm p_{0}$ are minimal, $M(p_{0})-\{m' \} \xcc
\xbm_{ \xdz }(X),$ so $m' \xbe \wt{ \xbm_{ \xdz }(X)}= \xbm (X).$ Thus,
as $ \xbm_{ \xdz }(X) \xcc \xbm (X),$ if $m' \xbe X- \xbm (X),$ then $m
\xbe X.$

Define now $X:=M(p_{0}) \xcv \{m\},$ $Y:=M(p_{0}).$

We first show that $ \xbm_{0}$ does not satisfy $( \xbm Cum).$
$ \xbm_{0}(X):=\{x \xbe X: \xCN \xcE A \xbe \xdy (A \xcc X:x \xbe A- \xbm
(A))\}.$ $m \xce \xbm_{0}(X),$ as $m \xce \xbm (X)= \wt{ \xbm_{ \xdz
}(X)}.$
Moreover, $m' \xce \xbm_{0}(X),$ as $\{m,m' \} \xbe \xdy,$ $\{m,m' \}
\xcc X,$ and $ \xbm (\{m,m' \})= \xbm_{ \xdz }(\{m,m' \})=\{m\}.$
So $ \xbm_{0}(X) \xcc Y \xcc X.$ Consider now $ \xbm_{0}(Y).$ As $m \xce
Y,$ for any $A \xbe \xdy,$ $A \xcc Y,$ if
$m' \xbe A,$ then $m' \xbe \xbm (A),$ too, by above argument, so $m' \xbe
\xbm_{0}(Y),$ and $ \xbm_{0}$ does
not satisfy $( \xbm Cum).$

We turn to $ \xbm_{1}.$

By Fact 4.7 (1), $ \xbm_{1}(X) \xcc \xbm_{0}(X),$ so $m,m' \xce
\xbm_{1}(X),$ and again $ \xbm_{1}(X) \xcc Y \xcc X.$
Consider again $ \xbm_{1}(Y).$ As $m \xce Y,$ for any $A \xbe \xdy,$ $
\xbm (A) \xcc Y,$ if $m' \xbe A,$
then $m' \xbe \xbm (A),$ too: if $M(p_{0})-\{m' \} \xcc A,$ then $m' \xbe
\wt{ \xbm_{ \xdz }(A)},$ if $M(p_{0})-\{m' \} \xcC A,$
but $m' \xbe A,$ then either $m' \xbe \xbm_{ \xdz }(A),$ or $m \xbe \xbm_{
\xdz }(A) \xcc \xbm (A),$ but $m \xce Y.$ Thus,
$( \xbm Cum)$ fails for $ \xbm_{1},$ too.

It remains to show that $ \xbm $ satisfies $( \xbm \xcc ),$ $( \xbm Cum),$
$( \xbm PR0),$ $( \xbm PR1).$
Note that by Fact  \Xl. $ \xbm_{ \xdz }$ satisfies $( \xbm Cum),$ as $
\xdz $ is smooth.
$( \xbm \xcc )$ is trivial. We show $( \xbm PRi)$ for $i=0,1.$ As $ \xbm_{
\xdz }(A) \xcc \xbm (A),$ by $( \xbm PR)$
and $( \xbm Cum)$ for $ \xbm_{ \xdz },$ $ \xbm_{ \xdz }(X) \xcc
\xbm_{0}(X)$ and $ \xbm_{ \xdz }(X) \xcc \xbm_{1}(X).$

To see this:
$ \xbm_{ \xdz }(X) \xcc \xbm_{0}(X):$ Let $x \xbe X- \xbm_{0}(X),$ then
there is $Y$ s.t. $x \xbe Y- \xbm (Y).Y \xcc X,$
but $ \xbm_{ \xdz }(Y) \xcc \xbm (Y),$ so by $Y \xcc X$ and $( \xbm PR)$
for $ \xbm_{ \xdz }$ $x \xce \xbm_{ \xdz }(X).$
$ \xbm_{ \xdz }(X) \xcc \xbm_{1}(X):$ Let $x \xbe X- \xbm_{1}(X),$ then
there is $Y$ s.t. $x \xbe Y- \xbm (Y),$ $ \xbm (Y) \xcc X,$
so $x \xbe Y- \xbm_{ \xdz }(Y)$ and $ \xbm_{ \xdz }(Y) \xcc X.$ $ \xbm_{
\xdz }(X \xcv Y) \xcc \xbm_{ \xdz }(X) \xcv \xbm_{ \xdz }(Y) \xcc X \xcc X
\xcv Y,$ so
$ \xbm_{ \xdz }(X \xcv Y)= \xbm_{ \xdz }(X)$ by $( \xbm Cum)$ for $ \xbm_{
\xdz }.$ $x \xbe Y- \xbm_{ \xdz }(Y)$ $ \xcp $ $x \xce \xbm_{ \xdz }(X
\xcv Y)$ by
$( \xbm PR)$ for $ \xbm_{ \xdz },$ so $x \xce \xbm_{ \xdz }(X).$

But by Fact 4.7, (3)
$ \xbm_{i}(X) \xcc \xbm (X).$ As by definition, $ \xbm (X)- \xbm_{ \xdz
}(X)$ is small, $( \xbm PRi)$ hold for
$i=0,1.$ It remains to show $( \xbm Cum)$ for $ \xbm.$ Let $ \xbm (X)
\xcc Y \xcc X,$ then
$ \xbm_{ \xdz }(X) \xcc \xbm (X) \xcc Y \xcc X,$ so by $( \xbm Cum)$ for $
\xbm_{ \xdz }$
$ \xbm_{ \xdz }(X)= \xbm_{ \xdz }(Y),$ so by definition of $ \xbm,$ $
\xbm (X)= \xbm (Y).$

(Note that by Fact 4.7 (10), $ \xbm_{0}= \xbm_{1}$ follows from $( \xbm
Cum)$ for $ \xbm,$ so we
could have demonstrated part of the properties also differently.)

$ \xcz $
\\[3ex]

\ee

By Fact 4.7 (3) and (8) and Proposition 3.2.4 in [Sch04], $ \xbm_{0}$
has a representation by a (transitive) preferential structure, if $ \xbm:
\xdy \xcp \xdy $
satisfies $( \xbm \xcc )$ and $( \xbm PR0),$ and $ \xbm_{0}$ is defined as
in Definition 4.3.

We thus have (taken from [Sch04]):

\bp

$\hspace{0.01em}$

% (+++*** Orig. No.:  Proposition 4.9: )

Let $Z$ be an arbitrary set, $ \xdy \xcc \xdp (Z),$ $ \xbm: \xdy \xcp
\xdy,$ $ \xdy $ closed under arbitrary
intersections and finite unions, and $ \xCQ,Z \xbe \xdy,$ and let $
\wt{.}$ be defined wrt. $ \xdy.$

(a) If $ \xbm $ satisfies $( \xbm \xcc ),$ $( \xbm PR0),$ then there is a
transitive
preferential structure $ \xdz $ over $Z$ s.t. for all $U \xbe \xdy $ $
\xbm (U)= \wt{ \xbm_{ \xdz }(U)}.$

(b) If $ \xdz $ is a preferential structure over $Z$ and
$ \xbm: \xdy \xcp \xdy $ s.t. for all $U \xbe \xdy $ $ \xbm (U)= \wt{
\xbm_{ \xdz }(U)},$ then $ \xbm $ satisfies $( \xbm \xcc ),$ $( \xbm
PR0).$

\ep

\paragraph{
Proof:
}

$\hspace{0.01em}$

(a) Let $ \xbm $ satisfy $( \xbm \xcc ),$ $( \xbm PR0).$ $ \xbm_{0}$ as
defined in Definition 4.3
satisfies properties $( \xbm \xcc ),$ $( \xbm PR)$ by Fact 4.7, (3) and
(8).
Thus, by Proposition 3.2.4 in [Sch04], there is a transitive
structure $ \xdz $
over $Z$ s.t. $ \xbm_{0}= \xbm_{ \xdz }$, but by $( \xbm PR0)$ $ \xbm
(U)= \wt{ \xbm_{0}(U)}= \wt{ \xbm_{ \xdz }(U)}$ for $U \xbe \xdy.$

(b) $( \xbm \xcc ):$ $ \xbm_{ \xdz }(U) \xcc U,$ so by $U \xbe \xdy $ $
\xbm (U)= \wt{ \xbm_{ \xdz }(U)} \xcc U.$

$( \xbm PR0):$ If $( \xbm PR0)$ is false, there is $U \xbe \xdy $ s.t. for
$U':= \xcV \{Y' - \xbm (Y' ):Y' \xbe \xdy,Y' \xcc U\}$ $ \wt{ \xbm
(U)-U' } \xcb \xbm (U).$ By $ \xbm_{ \xdz }(Y' ) \xcc \xbm (Y' ),$
$Y' - \xbm (Y' ) \xcc Y' - \xbm_{ \xdz }(Y' ).$ No copy of any $x \xbe Y'
- \xbm_{ \xdz }(Y' )$ with $Y' \xcc U,$ $Y' \xbe \xdy $ can be minimal in
$ \xdz \xex U.$ Thus, by $ \xbm_{ \xdz }(U) \xcc \xbm (U),$ $ \xbm_{ \xdz
}(U) \xcc \xbm (U)-U',$ so $ \wt{ \xbm_{ \xdz }(U)} \xcc \wt{ \xbm (U)-U'
} \xcb \xbm (U),$ $contradiction.$

$ \xcz $
\\[3ex]

We turn to the smooth case.

If $ \xbm: \xdy \xcp \xdy $ satisfies $( \xbm \xcc ),$ $( \xbm PR2),$ $(
\xbm CUM)$ and $ \xbm_{2}$ is
defined from $ \xbm $ as in Definition 4.3,
then $ \xbm_{2}$ satisfies $( \xbm \xcc ),$ $( \xbm PR),$ $( \xbm Cum)$ by
Fact 4.7 (3), (8), and (9), and
can thus be represented by a (transitive) smooth structure, by
Proposition 3.3.8 in [Sch04], and we finally have (taken from [Sch04]):

\bp

$\hspace{0.01em}$

% (+++*** Orig. No.:  Proposition 4.10: )

Let $Z$ be an arbitrary set, $ \xdy \xcc \xdp (Z),$ $ \xbm: \xdy \xcp
\xdy,$ $ \xdy $ closed under arbitrary
intersections and finite unions, and $ \xCQ,Z \xbe \xdy,$ and let $
\wt{.}$ be defined wrt. $ \xdy.$

(a) If $ \xbm $ satisfies $( \xbm \xcc ),$ $( \xbm PR2),$ $( \xbm CUM),$
then there is a transitive smooth
preferential structure $ \xdz $ over $Z$ s.t. for all $U \xbe \xdy $ $
\xbm (U)= \wt{ \xbm_{ \xdz }(U)}.$

(b) If $ \xdz $ is a smooth preferential structure over $Z$ and
$ \xbm: \xdy \xcp \xdy $ s.t. for all $U \xbe \xdy $ $ \xbm (U)= \wt{
\xbm_{ \xdz }(U)},$ then $ \xbm $ satisfies
$( \xbm \xcc ),$ $( \xbm PR2),$ $( \xbm CUM).$

\ep

\paragraph{
Proof:
}

$\hspace{0.01em}$

(a) If $ \xbm $ satisfies $( \xbm \xcc ),$ $( \xbm PR2),$ $( \xbm CUM),$
then $ \xbm_{2}$
defined from $ \xbm $ as in Definition 4.3 satisfies $( \xbm \xcc ),$ $(
\xbm PR),$
$( \xbm CUM)$ by Fact 4.7 (3), (8) and (9).
Thus, by Proposition 3.3.8 in [Sch04], there is a smooth transitive
preferential
structure $ \xdz $ over $Z$ s.t. $ \xbm_{2}= \xbm_{ \xdz }$, but by $(
\xbm PR2)$ $ \xbm (U)= \wt{ \xbm_{2}(U)}= \wt{ \xbm_{ \xdz }(U)}.$

(b)
$( \xbm \xcc ):$ $ \xbm_{ \xdz }(U) \xcc U$ $ \xcp $ $ \xbm (U)= \wt{
\xbm_{ \xdz }(U)} \xcc U$ by $U \xbe \xdy.$

$( \xbm PR2):$ If $( \xbm PR2)$ fails, then there is $U \xbe \xdy $ s.t.
for
$U':= \xcV \{Y' - \xbm (Y' ):$ $Y' \xbe \xdy,$ $ \xbm (U \xcv Y' ) \xcc
U\}$ $ \wt{ \xbm (U)-U' } \xcb \xbm (U).$

By $ \xbm_{ \xdz }(Y' ) \xcc \xbm (Y' ),$ $Y' - \xbm (Y' ) \xcc Y' -
\xbm_{ \xdz }(Y' ).$
But no copy of any $x \xbe Y' - \xbm_{ \xdz }(Y' )$ with $ \xbm_{ \xdz }(U
\xcv Y' ) \xcc \xbm (U \xcv Y' ) \xcc U$ can be minimal in $ \xdz \xex U:$
As $x \xbe Y' - \xbm_{ \xdz }(Y' ),$ if $<x,i>$ is any copy of $x,$ then
there is $<y,j> \xeb <x,i>,$ $y \xbe Y'.$
Consider now $U \xcv Y'.$ As $<x,i>$ is not minimal in $ \xdz \xex U \xcv
Y',$ by smoothness of $ \xdz $
there must be $<z,k> \xeb <x,i>,$
$<z,k>$ minimal in $ \xdz \xex U \xcv Y'.$ But all minimal elements of $
\xdz \xex U \xcv Y' $ must be in $ \xdz \xex U,$
so there must be $<z,k> \xeb <x,i>,$ $z \xbe U,$ thus $<x,i>$ is not
minimal in $ \xdz \xex U.$
Thus by $ \xbm_{ \xdz }(U) \xcc \xbm (U),$ $ \xbm_{ \xdz }(U) \xcc \xbm
(U)-U',$ so $ \wt{ \xbm_{ \xdz }(U)} \xcc \wt{ \xbm (U)-U' } \xcb \xbm
(U),$ $contradiction.$

$( \xbm CUM):$ Let $ \xbm (X) \xcc Y \xcc X.$ Now $ \xbm_{ \xdz }(X) \xcc
\wt{ \xbm_{ \xdz }(X)}= \xbm (X),$ so by smoothness of $ \xdz $
$ \xbm_{ \xdz }(Y)= \xbm_{ \xdz }(X),$ thus $ \xbm (X)= \wt{ \xbm_{ \xdz
}(X)}= \wt{ \xbm_{ \xdz }(Y)}= \xbm (Y).$
$ \xcz $
\\[3ex]
%  (4.3.2)  The ranked case
%  (4.3.2)  The ranked case
% %
% =========================
\subsubsection{The ranked case}

We recall from [Sch04] Notation 4.1, Definition 4.4, Fact 4.11, Proposition
4.12.

\bn

$\hspace{0.01em}$

% (+++*** Orig. No.:  Notation 4.1: )

(1) A $=$ $B \xFO C$ stands for: $A=B$ or $A=C$ or $A=B \xcv C.$

(2) Given $ \xeb,$ $a \xcT b$ means: neither $a \xeb b$ nor $b \xeb a.$

\en

\bd

$\hspace{0.01em}$

% (+++*** Orig. No.:  Definition 4.4: )

The new conditions for the minimal case are:

$( \xbm \xCQ )$ $X \xEd \xCQ $ $ \xcp $ $ \xbm (X) \xEd \xCQ,$

$( \xbm \xCQ fin)$ $X \xEd \xCQ $ $ \xcp $ $ \xbm (X) \xEd \xCQ $ for
finite $X,$

$( \xbm =)$ $X \xcc Y,$ $ \xbm (Y) \xcs X \xEd \xCQ $ $ \xcp $ $ \xbm (Y)
\xcs X= \xbm (X),$

$( \xbm =' )$ $ \xbm (Y) \xcs X \xEd \xCQ $ $ \xcp $ $ \xbm (Y \xcs X)=
\xbm (Y) \xcs X,$

$( \xbm \xFO )$ $ \xbm (X \xcv Y)$ $=$ $ \xbm (X) \xFO \xbm (Y),$

$( \xbm \xcv )$ $ \xbm (Y) \xcs (X- \xbm (X)) \xEd \xCQ $ $ \xcp $ $ \xbm
(X \xcv Y) \xcs Y= \xCQ,$

$( \xbm \xcv ' )$ $ \xbm (Y) \xcs (X- \xbm (X)) \xEd \xCQ $ $ \xcp $ $
\xbm (X \xcv Y)= \xbm (X),$

$( \xbm \xbe )$ $a \xbe X- \xbm (X)$ $ \xcp $ $ \xcE b \xbe X.a \xce \xbm
(\{a,b\}).$

\ed

We will use

\bfa

$\hspace{0.01em}$

% (+++*** Orig. No.:  Fact 4.11: )

The following properties $(2)-(9)$ hold, provided corresponding closure
conditions
for the domain $ \xdy $ are satisfied. We first enumerate these
conditions.

For (3), (4), (8): closure under finite unions.

For (2): closure under finite intersections.

For (6) and (7): closure under finite unions, and $ \xdy $ contains all
singletons.

For (5): closure under set difference.

For (9): suffienctly strong conditions - which are satisfied for the set
of
models definable by propositional theories.

Note that the closure conditions for (5), (6), (9) are quite different,
for this
reason, (5) alone is not enough.

(1) $( \xbm =)$ entails $( \xbm PR),$

(2) in the presence of $( \xbm \xcc ),$ $( \xbm =)$ is equivalent to $(
\xbm =' ),$

(3) $( \xbm \xcc ),$ $( \xbm =)$ $ \xcp $ $( \xbm \xcv ),$

(4) $( \xbm \xcc ),$ $( \xbm \xCQ ),$ $( \xbm =)$ entail:

(4.1) $( \xbm \xFO ),$

(4.2) $( \xbm \xcv ' ),$

(4.3) $( \xbm CUM),$

(5) $( \xbm \xcc )+( \xbm \xFO )$ $ \xcp $ $( \xbm =),$

(6) $( \xbm \xFO )+( \xbm \xbe )+( \xbm PR)+( \xbm \xcc )$ $ \xcp $ $(
\xbm =),$

(7) $( \xbm CUM)+( \xbm =)$ $ \xcp $ $( \xbm \xbe ),$

(8) $( \xbm CUM)+( \xbm =)+( \xbm \xcc )$ $ \xcp $ $( \xbm \xFO ),$

(9) $( \xbm PR)+( \xbm CUM)+( \xbm \xFO )$ $ \xcp $ $( \xbm =)$.

\efa

and

\bp

$\hspace{0.01em}$

% (+++*** Orig. No.:  Proposition 4.12: )

Let $ \xdy \xcc \xdp (U)$ be closed under finite unions, and contain
singletons.
Then $( \xbm \xcc ),$ $( \xbm \xCQ fin),$ $( \xbm =),$ $( \xbm \xbe )$
characterize ranked structures for which
for all finite $X \xbe \xdy $ $X \xEd \xCQ $ $ \xcp $ $ \xbm_{<}(X) \xEd
\xCQ $ hold, i.e. $( \xbm \xcc ),$ $( \xbm \xCQ fin),$ $( \xbm =),$ $(
\xbm \xbe )$
hold in such structures for $ \xbm_{<},$ and if they hold for some $ \xbm
,$ we can find
a ranked relation $<$ on $U$ s.t. $ \xbm = \xbm_{<}.$

Note that wlog. we may assume that the structure contains no copies.

\ep

We give now an easy version of representation results for ranked
structures
without definability preservation.

\bn

$\hspace{0.01em}$

% (+++*** Orig. No.:  Notation 4.2: )

We abbreviate $ \xbm (\{x,y\})$ by $ \xbm (x,y)$ etc.

\en

\bfa

$\hspace{0.01em}$

% (+++*** Orig. No.:  Fact 4.13: )

Let the domain contain singletons and be closed under $( \xcv ).$

Let for $ \xbm: \xdy \xcp \xdy $ hold:

$( \xbm =)$ for finite sets, $( \xbm \xbe ),$ $( \xbm PR3),$ $( \xbm \xCQ
fin).$

Then the following properties hold for $ \xbm_{3}$ as defined in
Definition 4.3:

(1) $ \xbm_{3}(X) \xcc \xbm (X),$

(2) for finite $X,$ $ \xbm (X)= \xbm_{3}(X),$

(3) $( \xbm \xcc ),$

(4) $( \xbm PR),$

(5) $( \xbm \xCQ fin),$

(6) $( \xbm =),$

(7) $( \xbm \xbe ),$

(8) $ \xbm (X)= \wt{ \xbm_{3}(X)}.$

\efa

\paragraph{
Proof:
}

$\hspace{0.01em}$

(1) Suppose not, so $x \xbe \xbm_{3}(X),$ $x \xbe X- \xbm (X),$ so by $(
\xbm \xbe )$ for $ \xbm,$ there is
$y \xbe X,$ $x \xce \xbm (x,y),$ $contradiction.$

(2) By $( \xbm PR3)$ for $ \xbm $ and (1), for finite $U$ $ \xbm (U)=
\xbm_{3}(U).$

(3) $( \xbm \xcc )$ is trivial for $ \xbm_{3}.$

(4) Let $X \xcc Y,$ $x \xbe \xbm_{3}(Y) \xcs X,$ suppose $x \xbe X-
\xbm_{3}(X),$ so there is $y \xbe X \xcc Y,$
$x \xce \xbm (x,y),$ so $x \xce \xbm_{3}(Y).$

(5) $( \xbm \xCQ fin)$ for $ \xbm_{3}$ follows from $( \xbm \xCQ fin)$ for
$ \xbm $ and (2).

(6) Let $X \xcc Y,$ $y \xbe \xbm_{3}(Y) \xcs X,$ $x \xbe \xbm_{3}(X),$ we
have to show $x \xbe \xbm_{3}(Y).$
By (4), $y \xbe \xbm_{3}(X).$
Suppose $x \xce \xbm_{3}(Y).$ So there is $z \xbe Y.x \xce \xbm (x,z).$ As
$y \xbe \xbm_{3}(Y),$ $y \xbe \xbm (y,z).$
As $x \xbe \xbm_{3}(X),$ $x \xbe \xbm (x,y),$ as $y \xbe \xbm_{3}(X),$ $y
\xbe \xbm (x,y).$
Consider $\{x,y,z\}.$ Suppose $y \xce \xbm (x,y,z),$ then by $( \xbm \xbe
)$ for $ \xbm,$ $y \xce \xbm (x,y)$ or
$y \xce \xbm (y,z),$ $contradiction.$
Thus $y \xbe \xbm (x,y,z) \xcs \xbm (x,y).$ As $x \xbe \xbm (x,y),$ and $(
\xbm =)$ for $ \xbm $ and finite sets,
$x \xbe \xbm (x,y,z).$ Recall that $x \xce \xbm (x,z).$ But for finite
sets $ \xbm = \xbm_{3},$ and by
(4) $( \xbm PR)$ holds for $ \xbm_{3},$ so it holds for $ \xbm $ and
finite sets. $contradiction$

(7) Let $x \xbe X- \xbm_{3}(X),$ so there is $y \xbe X.x \xce \xbm (x,y)=
\xbm_{3}(x,y).$

(8) As $ \xbm (X) \xbe \xdy,$ and $ \xbm_{3}(X) \xcc \xbm (X),$ $ \wt{
\xbm_{3}(X)} \xcc \xbm (X),$ so by $( \xbm PR3)$ $ \wt{ \xbm_{3}(X)}= \xbm
(X).$

$ \xcz $
\\[3ex]

\bfa

$\hspace{0.01em}$

% (+++*** Orig. No.:  Fact 4.14: )

If $ \xdz $ is ranked, and we define $ \xbm (X):= \wt{ \xbm_{ \xdz }(X)},$
and $ \xdz $ has no
copies, then the following hold:

(1) $ \xbm_{ \xdz }(X)=\{x \xbe X: \xcA y \xbe X.x \xbe \xbm (x,y)\},$ so
$ \xbm_{ \xdz }(X)= \xbm_{3}(X)$ for $X \xbe \xdy,$

(2) $ \xbm (X)= \xbm_{ \xdz }(X)$ for finite $X,$

(3) $( \xbm =)$ for finite sets for $ \xbm,$

(4) $( \xbm \xbe )$ for $ \xbm,$

(5) $( \xbm \xCQ fin)$ for $ \xbm,$

(6) $( \xbm PR3)$ for $ \xbm.$

\efa

\paragraph{
Proof:
}

$\hspace{0.01em}$

(1) holds for ranked structures.

(2) and (6) are trivial. (3) and (5) hold for $ \xbm_{ \xdz },$ so by (2)
for $ \xbm.$

(4) If $x \xce \xbm (X),$ then $x \xce \xbm_{ \xdz }(X),$ $( \xbm \xbe )$
holds for $ \xbm_{ \xdz },$ so there is $y \xbe X$ s.t.
$x \xce \xbm_{ \xdz }(x,y)= \xbm (x,y)$ by (2).

$ \xcz $
\\[3ex]

We summarize:

\bp

$\hspace{0.01em}$

% (+++*** Orig. No.:  Proposition 4.15: )

Let $Z$ be an arbitrary set, $ \xdy \xcc \xdp (Z),$ $ \xbm: \xdy \xcp
\xdy,$ $ \xdy $ closed under arbitrary
intersections and finite unions, contain singletons, and $ \xCQ,Z \xbe
\xdy,$ and let $ \wt{.}$
be defined wrt. $ \xdy.$

(a) If $ \xbm $ satisfies $( \xbm =)$ for finite sets, $( \xbm \xbe ),$ $(
\xbm PR3),$ $( \xbm \xCQ fin),$
then there is a ranked preferential structure $ \xdz $ without copies
over $Z$ s.t. for all $U \xbe \xdy $ $ \xbm (U)= \wt{ \xbm_{ \xdz }(U)}.$

(b) If $ \xdz $ is a ranked preferential structure over $Z$ without copies
and
$ \xbm: \xdy \xcp \xdy $ s.t. for all $U \xbe \xdy $ $ \xbm (U)= \wt{
\xbm_{ \xdz }(U)},$ then $ \xbm $ satisfies
$( \xbm =)$ for finite sets, $( \xbm \xbe ),$ $( \xbm PR3),$ $( \xbm \xCQ
fin).$

\ep

\paragraph{
Proof:
}

$\hspace{0.01em}$

(a) Let $ \xbm $ satisfy $( \xbm =)$ for finite sets, $( \xbm \xbe ),$ $(
\xbm PR3),$ $( \xbm \xCQ fin),$ then
$ \xbm_{3}$ as defined in Definition 4.3
satisfies properties $( \xbm \xcc ),$ $( \xbm \xCQ fin),$ $( \xbm =),$ $(
\xbm \xbe )$ by Fact 4.13.
Thus, by Proposition 4.12, there is a transitive structure $ \xdz $
over $Z$ s.t. $ \xbm_{3}= \xbm_{ \xdz }$, but by Fact 4.13 (8) $ \xbm
(U)= \wt{ \xbm_{3}(U)}= \wt{ \xbm_{ \xdz }(U)}$ for $U \xbe \xdy.$

(b) This was shown in Fact 4.14.

$ \xcz $
\\[3ex]

\section{
OUTLOOK: PATCHY DOMAINS AND WEAK DIAGNOSTIC INSTRUMENTS
}

% (+++*** Orig.:   (5)  OUTLOOK: PATCHY DOMAINS AND WEAK DIAGNOSTIC INSTRUMENTS  (-> DCB 4) )

We have seen in these pages two kinds of problems with repercussions on
representation questions:

(1) Lack of closure of the domain, in particular under finite union.

(2) Lack of definability preservation.

We summarize here the problem:

(1) If the domain is not closed under suitable operations, we may be
forced
to replace simple conditions like cumulativity by more complicated ones.

(2) If the operation is not definability preserving, we may be foced to
admit
exceptions.

In realistic situations, first, both problems may very well occur at the
same
time. Second, it is not evident that observable situations coincide with
logically describable situations. Moreover, it is not all clear that we
can
perform an ``experiment'' with all formulas or theories which are logically
describable.

Thus, we have three elements:

(1) A domain which may be quite patchy.

(2) Limited possibilities to carry out ``experiments'', i.e. a limited
number
of input situations, which may be only a subset of what is logically
possible.

(3) A limited number of observable results.

We have a shaggy situation, can do some experiments, and then might see
the
result only through a rather rough grid.

The author has intentionally chosen here the language of scientific
experiments,
as our situation seems close to the latter problem. So, perhaps, both
sides can
learn from each other.

At the moment, it seems difficult to obtain general results, as for
instance,
the different diagnostic instruments for smooth and general preferential
situations $( \xbm_{0}$ and $ \xbm_{2})$ show. As long as we have no
further information about
observability, for instance inconsistency need not always be observable,
as
$ \xCQ \xEd \wt{ \xCQ }$ may very well be possible. But the problem might
be too general, and
is amenable only for restricted cases - some of which were solved above.

Note that we have a similar situation in preferential structures, copies:
We
cannot observe copies directly, they are not describable in our language,
we
only see their effects.

It seems rather plausible that similar problems do not only occur in the
context of nonmonotonic reasoning. It seems safe to conjecture that there
are
existing systems and approaches which do not take into account the
problems
we have seen here, i.e. work fine under ideal situations, but are not
adapted
to their particular domain or observability problems.

\end{document}